\documentclass{article}

\usepackage[margin=25truemm]{geometry}
\usepackage{bm}
\usepackage{graphicx}
\usepackage{here}
\usepackage{bbm}
\usepackage{mathptmx}
\usepackage{amsmath}
\usepackage{amssymb}
\usepackage{mathrsfs}
\usepackage{amsthm}
\usepackage[authoryear]{natbib}
\usepackage{hyperref}
\usepackage{color}
\usepackage{url}

\theoremstyle{plain}
\newtheorem{theorem}{Theorem}[section]
\newtheorem{corollary}[theorem]{Corollary}
\newtheorem{proposition}[theorem]{Proposition}
\newtheorem{lemma}{Lemma}
\newtheorem{assumption}{Assumption}
\newtheorem{remark}{Remark}
\theoremstyle{definition}

\hypersetup{colorlinks=true, citecolor=blue, urlcolor=magenta, linkcolor=blue}

\allowdisplaybreaks 

\title{Higher-Order Asymptotic Properties of Kernel Density Estimator with Global Plug-In and Its Accompanying Pilot Bandwidth}

\author{Shunsuke Imai \thanks{Graduate School of Economics, Kyoto University, Yoshidahonmachi, Sakyoku, Kyoto, 606–8501, JAPAN, imai.shunsuke.57n@st.kyoto-u.ac.jp} \and Yoshihiko Nishiyama \thanks{Institute of Economic Research, Kyoto University, Yoshidahonmachi, Sakyoku, Kyoto, 606–8501, JAPAN, nishiyama@kier.kyoto-u.ac.jp}}

\begin{document}

\maketitle

\begin{abstract}
This study investigates the effect of bandwidth selection via a plug-in method on the asymptotic structure of the nonparametric kernel density estimator. We generalise the result of \cite{HK01} and find that the plug-in method has no effect on the asymptotic structure of the estimator up to the order of $O\{(nh_0)^{-1/2}+h_0^L\}=O(n^{-L/(2L+1)})$ for a bandwidth $h_0$ and any kernel order $L$ when the kernel order for pilot estimation $L_p$ is high enough. We also provide the valid Edgeworth expansion up to the order of $O\{(nh_0)^{-1}+h_0^{2L}\}$ and find that, as long as the $L_p$ is high enough , the plug-in method has an effect from on the term whose convergence rate is $O\{(nh_0)^{-1/2}h_0+h_0^{L+1}\}=O(n^{-(L+1)/(2L+1)})$. In other words, we derive the exact achievable convergence rate of the deviation between the distribution functions of the estimator with a deterministic bandwidth and with the plug-in bandwidth. In addition, we weaken the conditions on kernel order $L_p$ for pilot estimation by considering the effect of pilot bandwidth associated with the plug-in bandwidth. We also show that the bandwidth selection via the global plug-in method possibly has an effect on the asymptotic structure even up to the order of $O\{(nh_0)^{-1/2}+h_0^L\}$. Finally, Monte Carlo experiments are conducted to see whether our approximation improves previous results.

\end{abstract}

\numberwithin{equation}{section}

\newpage
\section{Introduction}

    In nonparametric statistics, the target of statistical inference is a function or an infinite dimensional vector $f$ that is not specifically modelled itself (See \cite{wasserman06} for introductive overviews, \cite{GN16} for mathematically unified understanding and \cite{IT07} and \cite{Chen07} for overviews especially in the context of economic literature). One of the important components of the function $f$ is the density function because, in statistics and its related fields, there are cases where we are interested in the distribution as a wage distribution (See e.g. \cite{DFL96}) or where a target of statistical inference depends on the density function as a conditional expectation function. Although there are different methods for estimating a density function, we focus on the estimator based on the kernel method, namely \textit{kernel density estimator} (KDE), also called Rosenblatt estimator or Rosenblatt-Parzen estimator after their pioneering works (\cite{Rosenblatt56} and \cite{Parzen62}).

    The first-order asymptotic properties of KDE have been studied over a long period and it has been proven that, under certain conditions, KDE has pointwise consistency and asymptotic normality (see e.g. \cite{Parzen62}, and the monograph by \cite[pp.28-30]{LR07}). As we will review in Section \ref{section:review}, the rate of convergence of KDE is slower than the parametric rate, and furthermore, becomes slower as the dimension increases. This property is called the curse of dimensionality. We can understand this as being the cost of using local data to avoid misspecification. \cite{Hall91} has clarified the higher-order asymptotic properties of the estimator in both non-Studentised and Studentised cases. The asymptotic expansion of  KDE is no longer a series of $n^{-1/2}$ as parametric estimators, but a series of $(nh)^{-1/2}$, even in the non-Studentised case; it is a more complicated series in the Studentised case, where $n$ and $h$ are the sample size and bandwidth, respectively.

    Bandwidth $h$ specifies the flexibility of statistic models and is adjusted between the bias and variance trade-offs in the sense that creating flexible models and consequently decreasing the bias results in increasing variance while creating non-flexible models and decreasing the variance results in increasing bias. It is well known that the performance of the kernel-based estimators depends greatly on the bandwidth, not so much on the kernel function. By defining a loss function, one can compute the theoretically optimal bandwidth $h_0$ that minimises loss. For example, mean integrated squared error (MISE) is the most commonly used global loss measure. However, in practice, such a bandwidth is typically infeasible because it depends on the unknown density. Therefore, one has to choose the bandwidth in a data-driven manner. Among the many bandwidth selection methods, two famous ones are cross-validation and plug-in method. In this paper, we focus on the latter.

    It is natural to ask whether the choice of bandwidth affects the asymptotic structure of the estimator. \cite{Ichimura00} and \cite{LL10} have considered the asymptotic distribution of kernel-based non/semiparametric estimators with data-driven bandwidth. They argue that, under certain conditions, the bandwidth selection has no effect on the first-order asymptotic structure of the estimators. \cite{HK01} showed that the bandwidth selection by the global plug-in method also has no effect on the asymptotic structure of KDE up to the order of $O(n^{-2/5})$ for $L=2$ and $L_p=6$, where $L$ and $L_p$ are kernel orders for the density estimation and estimation of an unknown part of the optimal bandwidth, respectively.

    Our contributions are fivefold. 
    First, we provide the Edgeworth expansion of KDE with global plug-in bandwidth up to the order of $O\{(nh_0)^{-1}+h_0^{2L}\}=O(n^{\frac{-2L}{2L+1}})$ and show that the bandwidth selection by the plug-in method begins to affect the term whose convergence rate is $O\{(nh_0)^{-1/2}h_0+h_0^{L+1}\}=O(n^{\frac{-(L+1)}{2L+1}})$ under the condition that $L_p$ is large enough. 
    Second, we generalise Theorem 3.2 of \cite{HK01}, which states that bandwidth selection via the global plug-in method has no effect on the asymptotic structure of KDE up to the order of $O\{(nh_0)^{-1/2}+h_0^{L}\}=O(n^{\frac{-L}{2L+1}})$.
    Their results limit the order of kernel functions $K(u)$ and $H(u)$ to $L=2,L_p=6$, respectively, but we show that they are valid for general orders $L$  as well under the condition that $L_p$ is large enough. 
    Third, we explore Edgeworth expansion of KDE with deterministic bandwidth in more detail than \cite{Hall91}. We show that Edgeworth expansion of Standardised KDE with deterministic bandwidth has the term of order $O\{(nh_0)^{-1/2}+h_0^{L}\}=O(n^{\frac{-L}{2L+1}})$ right after the term $\Phi(z)$ with a gap between them. After that however, the terms decrease at the rate of $O(h_0)=O(n^{\frac{-1}{2L+1}})$.
    However, the result of \cite{HK01} and our results above need the kernel order $L_p$ for the estimation of unknown parts of the optimal bandwidth to be high enough. We have two motivations to avoid imposing this condition on $L_p$. One is that although the higher-order kernel is theoretically justified, in terms of implementation using a computer, it has undesirable properties. The other is that the condition forces pilot bandwidth to be relatively large but the range is restrictive especially in multidimensional settings. For details of the latter motivation, see the seminal works of \cite{CCJ10,CCJ13,CCJ14a,CCJ14b} and \cite{CCJ18}.  Then, as a fourth contribution , we weaken this condition on $L_p$ assumed by \cite{HK01} and our Theorem \ref{T:Main Theorem} and provide the Edgeworth expansion including the effect of pilot bandwidth up to the order of $O\{(nh_0)^{-1}+h_0^{2L}\}$. In this situation, the bandwidth selection via the global plug-in method possibly has an effect on the asymptotic structure of KDE even up to the order of $O\{(nh_0)^{-1/2}+h_0^L\}$ (for example, when $L=2$ and $L_p=2$).
    Finally, we consider the intersectional effect of the bandwidth selection via the global plug-in method, its accompanying pilot bandwidth, and Studentisation.
    The proof of our main theorem owes much to \cite{NR00}. 
    They have established the valid Edgeworth expansion for the semiparametric density-weighted averaged derivatives estimator of the single index model, which has an exact second-order $U$-statistic form. Although the higher-order asymptotic structure of $U$-statistics had been studied before \cite{NR00} (See e.g. \cite{CJV80}), the estimator is different from \textit{standard} $U$-statistics in that it is $U$-statistics whose kernel depends on the sample size $n$ through the bandwidth. Since KDE with plug-in bandwidth can also be approximated by a sum of first- and second-order $U$-statistics whose kernel depends on the sample size $n$ through the bandwidth, we can benefit from their proof.

    The remainder of this paper is organised as follows. In the next section, we introduce KDE and review its known properties. Section \ref{section:main result} provides the main results, namely the Edgeworth expansion of the estimator with the global plug-in bandwidth. In section \ref{section:simulation study}, we employ Monte Carlo studies to compare our results with those of previous works. Section \ref{Section:Discussion} concludes and discusses future research directions.

\section{Review of the Estimator's Properties} \label{section:review}
\subsection{Estimator and Its First Order Properties}
\begin{assumption} \label{A:DGP}
Let $\{X_i\}_{i=1}^n$ be a random sample with an absolutely continuous distribution with
Lebesgue density $f$.
\end{assumption}

First, we introduce nonparametric KDE $\hat{f}$ for unknown density $f$. Estimator $\hat{f}$ at a point $x$ with a bandwidth $h$ is defined as follows:
\begin{equation}
    \hat{f}_h(x) \equiv \frac{1}{nh}\sum_{i=1}^n K\left(\frac{X_i-x}{h}\right) \equiv \frac{1}{nh}\sum_{i=1}^n K_{i,h}(x)\nonumber,
\end{equation}
where $K$ is a kernel function, and we say that $K$ is a $L$-th order kernel, for a positive integer $L$, if
\begin{equation}
    \int u^lK(u)du = 
    \begin{cases}1 & (l=0) \\
    0 & (1\leqq l\leqq L-1)\\
    C \neq 0,<\infty& (l=L). \nonumber\end{cases}
\end{equation}
\begin{assumption} \label{A: L times differentiable}
In a neighbourhood of $x$, $f$ is $L$ times continuously differentiable and its first $L$ derivatives are bounded.
\end{assumption}
\begin{assumption} \label{A: Kernel}
Kernel function $K$ is a bounded, even function with a compact support, of order $L\geqq 2$ and $\int  K(u)du=1$. 
\end{assumption}
\begin{assumption} \label{A: interior point}
$x$ is an interior point in the support of $X$.
\end{assumption}
\begin{assumption} \label{A: nh infty}
    $h\rightarrow 0,~~~nh\rightarrow \infty~~~\text{as}~~n\rightarrow\infty$
\end{assumption}

KDE has pointwise consistency and asymptotic normality for an interior point in the support of $X$. Although it also converges uniformly for an interior point in the support of $X$, we only review pointwise properties because we investigate the pointwise higher-order asymptotics of KDE with global plug-in bandwidth. Under Assumption \ref{A:DGP}--\ref{A: Kernel}, we can expand mean squared error (MSE) of $\hat{f}_h(x)$ as follows:
\begin{align}
    MSE[\hat{f}(x)]&\equiv\mathbb{E}[\{\hat{f}(x)-f(x)\}^2]=\Bigl(C_Lf^{(L)}(x)h^{L}\Bigl)^2+\frac{R(K)f(x)}{nh}+o\{h^{2L}+(nh)^{-1}\}, \label{MSE}
\end{align}
where $R(K)=\int K(u)^2du, C_L=\frac{1}{L!}\int u^LK(u)du$.
Therefore, Markov's inequality, Assumptions \ref{A:DGP}--\ref{A: nh infty}, and (\ref{MSE}) imply pointwise consistency $\hat{f}_h(x)\xrightarrow{p}f(x)$. Moreover, we can show that KDE has asymptotic normality by applying Lindberg-Feller's central limit theorem:
\begin{equation}
\sqrt{nh}\left(\hat{f}_h(x)-\mathbb{E}\hat{f}_h(x)\right)\xrightarrow{d} N\Bigl(0,R(K)f(x)\Bigl). \nonumber
\end{equation}

\begin{remark} \label{R:Asymptotic Bias}
Since $\mathbb{E}[\hat{f}_h(x)]\approx f(x)+C_Lf^{(L)}(x)h^L$, the statistics centred by $f(x)$ asymptotically follows a zero-mean normal distribution if  $nh^{2L+1}\rightarrow 0$ holds. However, the theoretically optimal bandwidth does not satisfy this condition, as we will discuss later. Therefore, we consider the statistics centred by $\mathbb{E}[\hat{f}_h(x)]$, not $f(x)$. For recent studies on asymptotic bias of KDE, see, for example, \cite{HH13} and \cite{CCF18}. For other nonparametric estimators, recent related studies are those by \cite{AK18}, \cite{CCT14},\cite{CCF20,CCF22} and \cite{Schennach20}.
\end{remark}

\subsection{Plug-In Method}
Bandwidth $h$ is a parameter that analysts need to choose in advance. One of the criteria for bandwidth selection is the mean integrated squared error (MISE):
\begin{align}
    MISE(h)=\int \mathbb{E}[\{\hat{f}_h(x)-f(x)\}^2]dx. \nonumber
\end{align}
The theoretically optimal bandwidth is the one that minimises MISE and, from the MISE expansion, this bandwidth is defined as follows:
\begin{align}
  h_0&=\left(\frac{R(K)}{2LC_L^2I_L}\right)^{\frac{1}{2L+1}}n^{-\frac{1}{2L+1}}, \nonumber
\end{align}
where $I_L=\int f^{(L)}(x)^2dx$. Although $h_0$ would perform the best, it is infeasible because $I_L$ is unknown, so one has to select the bandwidth from the available data. We examine the effect of a certain plug-in method on the distribution of the estimator.

Several plug-in methods have been proposed so far (see e.g. \cite{HSJM91}, \cite{SJ91}). 
In this paper, we adopt as \cite{HK01}, a simple plug-in method that estimates $I_L$ directly and nonparametrically using the estimator proposed by \cite{HM87}. 
Their estimator, $\hat{I}_L$ for $I_L$, is given as follows:
\begin{equation}
    \hat{I}_L=\binom{n}{2}^{-1}\sum_{i=1}^{n-1}\sum_{j=i+1}^nb^{-(2L+1)}H^{(2L)}\left(\frac{X_i-X_j}{b}\right)\equiv\binom{n}{2}^{-1}\sum_{i=1}^{n-1}\sum_{j=i+1}^n \hat{I}_{Lij}, \nonumber
\end{equation}
where $b$ (called \textit{pilot bandwidth}) is a bandwidth for estimation of $I_L$, different from $h$, and $H$ is a kernel function of order $L_p$.

Another estimator for $I_L$ proposed by \cite{HM87} is 
\begin{align}
    \int \left\{\hat{f}^{(L)}(x)\right\}^2dx = \frac{1}{nb^{2L+1}}\Bar{H}^{(L)}(0) + \frac{1}{n^2b^{2L+1}}\sum_{i=1}^n\sum_{j\neq i}^n\Bar{H}^{(L)}\left(\frac{X_i-X_j}{b}\right) \label{eq:I hat convo}
    \end{align}
where $\hat{f}^{(L)}(x) \equiv \frac{1}{nb^{L+1}}\sum_{i=1}^n K^{(L)}\left(\frac{X_i-x}{b}\right)$  and $\Bar{H}^{(L)}(v) \equiv \int H^{(L)}(u)H^{(L)}(v-u)du$. \cite{HM87} state that 'the first term does not make use of the data, and hence may be thought of as adding a type of bias in the estimator. This motivates the estimator'.
\begin{align}
    \hat{I}_L^{convo} \equiv \frac{1}{n(n-1)b^{2L+1}}\sum_{i=1}^n\sum_{j\neq i}^n\Bar{H}^{(L)}\left(\frac{X_i-X_j}{b}\right). \label{eq:I hat convo 2}
\end{align}

\begin{remark}
     $\hat{I}_L$ and $\hat{I}_L^{convo}$ can be negative in small samples. Although they are asymptotically justified, it can cause problems in empirical applications. \cite{HK01} avoid this problem by using $|\hat{I}_L|$ instead of $\hat{I}_L$. Another way is to use $\int\left\{\hat{f}^{(L)}(x)\right\}^2dx$ instead of $\hat{I}_L^{convo}$. In Section \ref{section:simulation study}, we employ the Monte Carlo Study in these two ways. 
\end{remark}

\begin{assumption} \label{A:rate b}
$b=cn^{-2/(4L+2L_p+1)}$
\end{assumption}

Proposition \ref{p:Linearization of Plug-In Bandwidth} provides the expansion of the plug-in bandwidth (defined as $\hat{h}$) and plays an essential role in the derivation of the asymptotic expansion of KDE with the plug-in bandwidth. We assume additional conditions for Proposition \ref{p:Linearization of Plug-In Bandwidth}:
\begin{assumption} \label{A:2L+L_p times differentiable}
In a neighbourhood of $x$, $f$ is $(2L+L_p)$-times continuously differentiable and its first $(2L+L_p)$ derivatives are bounded.
\end{assumption}
\begin{assumption} \label{A:pilot kernel}
Kernel function $H$ is a bounded, even function with compact support, of order $L_p\geqq 2$, (2L)-times continuously differentiable and for all integers $k$ such that $1\leqq k \leqq 2L-1$, $\lim_{u\rightarrow\pm\infty}|H^{(k)}(u)|\rightarrow0$.
\end{assumption}
Assumption \ref{A:2L+L_p times differentiable} gives regularity conditions on the smoothness of the estimand, which implies Assumption \ref{A: L times differentiable}. Assumption \ref{A:pilot kernel} is on the kernel function $H$ for the estimation of $I_L$, and the condition at the infinity of $u$ is necessary for integration by parts in the expanding process of $(\hat{h}-h_0)/h_0$. These assumptions can be interpreted as a generalisation of assumption $(A_{gpi})$ of \cite{HK01} to $K$ of order $L$ and $H$ of order $L_p$. 

\begin{proposition}[Expansion of Plug-In Bandwidth]\label{p:Linearization of Plug-In Bandwidth}
 Under Assumptions \ref{A:DGP}, \ref{A: Kernel}, \ref{A: interior point}, \ref{A:rate b}, \ref{A:2L+L_p times differentiable} and \ref{A:pilot kernel}, and additionally \ref{A:L_p order 1} for Theorem \ref{T: Hall and Kang}, \ref{A:L_p order 2} for Theorem \ref{T:Main Theorem}, and \ref{A:L_p order 3} and \ref{A:L_p order 4} for Theorem \ref{T:Edgeworth Expansion Including Pilot Bandwidth} and \ref{T:Edgeworth Expansion Including Pilot Bandwidth and the effect of Studentisation}, we can expand $(\hat{h}-h_0)/h_0$ as follows:
\begin{align}
    \frac{\hat{h}-h_0}{h_0}= \frac{-C_{PI}}{n}\sum_{i=1}^nV_i-\frac{C_{PI}}{2}\binom{n}{2}^{-1}\sum_{i=1}^{n-1}\sum_{j\neq i}^nW_{ij} + o_p\{(nh_0)^{-1}+h_0^{2L}\} \label{Linearized PI}
\end{align}
where
\begin{align}
    & C_{PI} = \frac{2}{2L+1}I_L^{-1},  \nonumber \\
    & V_i \equiv  \{f^{(2L)}(X_i)-\mathbb{E}f^{(2L)}(X_i)\}+\frac{\int u^{L_p}H(u)du}{(L_p)!}b^{L_p} \Bigl\{f^{(2L+L_{p})}(X_i)-\mathbb{E}f^{(2L+L_{p})}(X_i)\Bigl\} + o_p(n^{-1/2}b^{L_p}), \nonumber\\
    & W_{ij} \equiv \Bigl\{\hat{I}_{Lij} - \mathbb{E}\left[\hat{I}_{Lij}|X_i\right] - \mathbb{E}\left[\hat{I}_{Lij}|X_j\right] + \mathbb{E}\left[\hat{I}_{Lij}\right]\Bigl\}. \nonumber
\end{align}
\end{proposition}
\noindent The proof is in \ref{Proof: Linearization of PI bandwidth}.

\begin{remark}\label{R:projection of h^}
The first term on the right-hand side of (\ref{Linearized PI}) reflects the projection term of the Hoeffding-decomposition of $\hat{I}_L$, whose convergence rate is $O_p(n^{-1/2})$. The second term reflects the quadratic term of the decomposed $\hat{I}_L$, whose convergence rate is $O_p(n^{-1}b^{-(4L+1)/2})$.

\end{remark}
\begin{remark}\label{R:pilot}
Since the MSE optimal rate of $b$ is $O_p(n^{\frac{-2}{4L+2L_p+1}})$ from \cite{HM87}, for example, when one chooses the pilot bandwidth via the rule of thumb (see \cite{Silverman86}) or second-stage plug-in method, the convergence rate of the second term in (\ref{Linearized PI}) is $O_p(n^{-1/2}b^{L_p})=O_p(n^{\frac{-4L-6L_p-1}{2(4L+2L_p+1)}})$. We can make the second term in $V_i$ as small as we like up to the order of $O(n^{-3/2})$ by letting kernel order $L_p$ be large enough. This is not an unrealistic statement; for example, when one uses a second order kernel function $K$, adopting a second order kernel function is sufficient to make the effect of the second order term negligible in the sense that they do not affect on the asymptotic structure of KDE up to the order of $O\{(nh_0)^{-1}+h_0^{2L}\}$. 
\end{remark}

\begin{remark}
Since the MSE optimal rate of $b$ is $O_p(n^{\frac{-2}{4L+2L_p+1}})$ from \cite{HM87}, for example, when one choose the pilot bandwidth via rule of thumb (see \cite{Silverman86}), the convergence rate of the second term in (\ref{Linearized PI}) is $O_p(n^{-1}b^{-(4L+1)/2})=O_p(n^{\frac{-2L_p}{4L+2L_p+1}})$. This implies that we can also make the third term of (\ref{Linearized PI})  as small as we like up to the order of $O(n^{-1})$ by letting kernel order $L_p$ be large enough. Although we cannot immediately identify how large $L_p$ needs to be to make the effect of pilot bandwidth negligible without deriving the Edgeworth expansion with pilot bandwidth, as we will see later, one has to adopt a considerably large $L_p$.
\end{remark}

\begin{remark} \label{R:degenerate I_L}
     Since the convergence rate of the second term is $O_p(n^{-1}b^{-(4L+1)/2})=O_p(n^{\frac{-2L_p}{4L+2L_p+1}})$, if not $L_p > (4L+1)/2$, the convergence rate of the second term is slower than that of the first term. In order to ignore the effect of the second term, \cite{HK01} provide the expansion under the condition that $L=2$ and $L_p=6$. The generalised version of this assumption is provided as Assumption \ref{A:L_p order 1}. In addition, we weaken the condition by considering the effect of pilot bandwidth. We provide such results as Theorem \ref{T:Edgeworth Expansion Including Pilot Bandwidth} and \ref{T:Edgeworth Expansion Including Pilot Bandwidth and the effect of Studentisation}. 
\end{remark}

\subsection{Review of Previous Studies}

Theorem 2.1 of \cite{Hall91} established the Edgeworth expansion for KDE with a deterministic bandwidth, which we replicate in Proposition \ref{p:Edgeworth Expansion for KDE with optimal bandwidth}. Let $S_h(x)$ be the Standardised version of KDE with a bandwidth $h$:
\begin{align}
     S_h(x)\equiv\frac{\sqrt{nh}\{\hat{f}_h(x)-\mathbb{E}\hat{f}_h(x)\}}{\mu_{20}(h)^{1/2}}, \nonumber
\end{align}
where $K_{i,h}(x)=K\left(\frac{X_i-x}{h}\right)$ and
\begin{equation}
    \mu_{kl}(h)\equiv h^{-1}\mathbb{E}\left[\left\{K_{i,h}(x)-\mathbb{E}[K_{i,h}(x)]\right\}^k\left\{K_{i,h}(x)^2-\mathbb{E}[K_{i,h}(x)^2]\right\}^l\right]. \nonumber
\end{equation}

\begin{assumption} \label{A: nh/long infty}
$h\rightarrow 0,~~~nh/\log n\rightarrow\infty~~~\text{as}~~n\rightarrow\infty$
\end{assumption}
\begin{assumption}[\textbf{Cram\'er Condition}] \label{A:Cramer}
For a sufficiently small $h$:
\begin{equation}
    \sup_{t\in\mathbb{R}}\left|\int_{-\infty}^{\infty}\exp\left\{itK(u)\right\}f(x-uh)du\right|< 1 \nonumber.
\end{equation}
\end{assumption}
\begin{remark}
    Assumption \ref{A:Cramer} is a high-level condition. Lemma 4.1 in \cite{Hall91} shows that primitive condition (2.1) in \cite{Hall91} implies Assumption \ref{A:Cramer}. Moreover, Assumption \ref{A:Cramer} is weaker than the Cram\'er condition in Lemma 4.1 of \cite{Hall91}. This is because Theorem \ref{T:Main Theorem} only deal with the Standardised case, while \cite{Hall91} also deals with the Studentised case. Our Theorem \ref{T:Edgeworth Expansion Including Pilot Bandwidth and the effect of Studentisation} needs the same Cram\'er condition as \cite{Hall91}.
\end{remark}
\begin{remark} \label{R: cramer condition}
    Assumption \ref{A:Cramer} rules out the uniform kernel, but many kernels which are practically used will satisfy this condition. However, as stated in \cite{Hall91}, one can also derive the Edgeworth expansion in the case of the uniform kernel by routine methods for lattice-valued random variables. 
\end{remark}

\begin{proposition}[\cite{Hall91}, Expansion with a Deterministic Bandwidth]\label{p:Edgeworth Expansion for KDE with optimal bandwidth}
Under Assumptions \ref{A:DGP}, \ref{A: Kernel}, \ref{A: interior point}, \ref{A: nh/long infty}, and \ref{A:Cramer}, the following expansions are valid:
    \begin{align}
        &\sup_{z \in \mathbb{R}} \left|\mathbb{P}(S_h(x)\leqq z)-\Phi(z)-\phi(z)\Biggl[(nh)^{-1/2}p_1(z)\Biggl]\right|=o\{(nh)^{-1/2}\} \nonumber\\
        &\sup_{z \in \mathbb{R}} \left|\mathbb{P}(S_h(x)\leqq z)-\Phi(z)-\phi(z)\Biggl[(nh)^{-1/2}p_1(z)+(nh)^{-1}p_2(z)\Biggl]\right|=o\{(nh)^{-1}\}, \label{EE Third with OPT}
    \end{align}
    where $\Phi(z)$ and $\phi(z)$ are the distribution and density functions at $z$ of a standard normal random variable, respectively, and:
    \begin{align}
        p_1(z)&=-\frac{1}{6}\mu_{20}(h)^{-3/2}\mu_{30}(h)(z^2-1), \nonumber\\
        p_2(z)&=-\frac{1}{24}\mu_{20}(h)^{-2}\mu_{40}(h)(z^3-3z)-\frac{1}{72}\mu_{20}(h)^{-3}\mu_{30}^2(z^5-10z^3+15z). \nonumber
    \end{align}
\end{proposition}
\noindent See \cite{Hall91} for the proof.

These results are the Edgeworth expansion of KDE up to the order of $O(\{(nh)^{-1/2}\})$ and $O(\{(nh)^{-1}\})$, respectively. However, bandwidth in his results is still deterministic. In this paper, we study KDE with data-driven bandwidth $\hat{f}_{\hat{h}}$ at a point $x$. The next proposition decomposes the $\hat{f}_{\hat{h}}$ into terms that include the effect of bandwidth selection and ones that do not.
\begin{assumption} \label{A: Kernel ''}
Kernel function $K$ is twice continuously differentiable.
\end{assumption}
\begin{proposition}[Expansion of KDE with Data-Driven Bandwidth]\label{p:Higher-Order Expansion of KDE with Data-Driven Bandwidth} 
Under Assumptions \ref{A:DGP}, \ref{A: interior point}, \ref{A: nh infty}, \ref{A:2L+L_p times differentiable}, \ref{A:pilot kernel}, \ref{A: nh/long infty}, and \ref{A: Kernel ''}, expanding $\hat{f}_{\hat{h}}(x)$ around $\hat{h}=h_0$ yields:
\begin{align}
    \hat{f}_{\hat{h}}(x) &\equiv \frac{1}{n\hat{h}}\sum_{i=1}^nK_{i,\hat{h}}(x) \nonumber\\
    &= \hat{f}_{h_0}(x)-\left(\frac{\hat{h}-h_0}{h_0}\right)\Gamma_{KDE_1}+\frac{1}{2}\left(\frac{\hat{h}-h_0}{h_0}\right)^2\Gamma_{KDE_2}+o_p\left(\left(\frac{\hat{h}-h_0}{h_0}\right)^2\Gamma_{KDE_2}\right), \label{HO expansion of KDE with data driven bandwidth}
\end{align}
where letting $u_{i,h}(x) \equiv \left(\frac{X_i-x}{h}\right) $, $\Gamma_{KDE_1}$ and $\Gamma_{KDE_2}$ are defined as follows.
\begin{align}
    & \Gamma_{KDE_1} \equiv \frac{1}{nh_0}\sum_{i=1}^n\Bigl\{K'_{i,h_0}(x)u_{i,h_0}(x)+K_{i,h_0}(x)\Bigl\}, \nonumber\\
    & \Gamma_{KDE_2} \equiv  \frac{1}{nh_0}\sum_{i=1}^n\Bigl\{2K_{i,h_0}(x)+4K'_{i,h_0}(x)u_{i,h_0}(x)+K''_{i,h_0}(x)u_{i,h_0}(x)^2\Bigl\}. \nonumber
\end{align}
\end{proposition}
 Let $S_{PI}(x)$ be the Standardised version of KDE with global plug-in bandwidth and define $\mu_{kl}=\mu_{kl}(h_0)$. Noting that expanding $\hat{h}^{1/2}$ around $\hat{h}=h_0$ yields $\hat{h}^{1/2}=h_0^{1/2}+\frac{1}{2}h_0^{-1/2}(\hat{h}-h_0)+O_p\{(\hat{h}-h_0)^2h_0^{-3/2}\}$, we have,
\begin{align}
    S_{PI}(x)&\equiv\frac{\sqrt{n\hat{h}}\{\hat{f}_{\hat{h}}(x)-\mathbb{E}\hat{f}_{h_0}(x)\}}{\mu_{20}^{1/2}}\nonumber\\&=S_{h_0}(x)-\frac{\sqrt{nh_0}\left(\frac{\hat{h}-h_0}{h_0}\right)\Gamma_{KDE_1}-\frac{\sqrt{nh_0}}{2}\left(\frac{\hat{h}-h_0}{h_0}\right)^2\Gamma_{KDE_2}}{\mu_{20}^{1/2}} \nonumber\\
    &\qquad + \frac{1}{2}S_{h_0}(x)\left(\frac{\hat{h}-h_0}{h_0}\right)-\frac{\sqrt{nh_0}\left(\frac{\hat{h}-h_0}{h_0}\right)^2\Gamma_{KDE_1}}{2\mu_{20}^{1/2}}+\frac{1}{6}S_{h_0}(x)\left(\frac{\hat{h}-h_0}{h_0}\right)^2+s.o.\label{def SPI}
\end{align}

\begin{assumption} \label{A:Ku}
$\lim_{u\rightarrow\pm \infty}|K(u)u|\rightarrow 0$
\end{assumption}
The following theorem generalises the kernel orders of Theorem 3.2 in \cite{HK01}. Their theorem  specifically sets the order of the kernels to be $L=2$ and $L_p=6$, and we prove that it holds for general kernel orders $L$ and $L_p$.

\begin{assumption} \label{A:L_p order 1}
$L_p > (4L+1)/2$.
\end{assumption}
\begin{remark}
     As stated in Remark \ref{R:degenerate I_L}, this assumption is also interpreted as the generalisation of $(A_{gpi})$ in \cite{HK01}.
\end{remark}

\begin{theorem}[Second Order Equivalence] \label{T: Hall and Kang}
Under Assumptions \ref{A:DGP}, \ref{A: interior point}, \ref{A:2L+L_p times differentiable}, \ref{A:pilot kernel}, \ref{A:Cramer}, \ref{A: Kernel ''}, \ref{A:Ku} and \ref{A:L_p order 1} the following expansion is valid:
\begin{equation}
        \sup_{z \in \mathbb{R}} \left|\mathbb{P}(S_{\text{PI}}(x)\leqq z)-\Phi(z)-\phi(z)\Biggl[(nh_0)^{-1/2}p_1(z)\Biggl]\right|=o\{(nh_0)^{-1/2}+h_0^{L}\}.
\end{equation}
\end{theorem}
\noindent See \ref{Proof:Hall and Kang} for the proof. We note that $(nh_0)^{-1/2}$ and $h_0^{L}$ have the same order of $O(n^{-L/(2L+1)})$, but we write the right hand side in this manner to clarify the effect of the variance and bias. Comparing this result with the first half of Proposition \ref{p:Edgeworth Expansion for KDE with optimal bandwidth}, we see that the bandwidth selection via the global plug-in method has no effect on the asymptotic structure of KDE up to the order of $O\{(nh_0)^{-1/2}+h_0^L\}$ as long as the kernel order of the kernel function for pilot estimation is large enough to satisfy Assumption \ref{A:L_p order 1}.

\begin{remark}
     When Assumption \ref{A:L_p order 1} breaks, Theorem \ref{T: Hall and Kang} does not hold. In other words, the bandwidth selection has an effect on the asymptotic structure of KDE up to the order of $O\{(nh_0)^{-1/2}+h_0^L\}$. Theorem \ref{T:Edgeworth Expansion Including Pilot Bandwidth} deals with this issue.
\end{remark}

\section{Main Results} \label{section:main result}

As stated in Theorem 3.2 in \cite{HK01} or our Theorem \ref{T: Hall and Kang}, bandwidth
selection via the global plug-in method has no effect on the asymptotic properties of KDE up
to the order of $O\{(nh_0)^{-1/2}+h_0^L\}=O(n^{-L/(2L+1)})$ when one uses a sufficiently high order kernel for the estimation of $I_L$. Section \ref{subsec: second order Edgeworth Expansion} provides a valid Edgeworth expansion for
 KDE with plug-in bandwidth up to the order of $O\{(nh_0)^{-1}+h_0^{2L}\}=O(n^{-2L/(2L+1)})$ in Theorem \ref{T:Main Theorem}. This expansion possesses a form comparable with that in \cite{Hall91}. 
In Section \ref{section: comparison}, we rewrite the expansions in Proposition \ref{p:Edgeworth Expansion for KDE with optimal bandwidth} and Theorem \ref{T:Main Theorem} to derive
the expansions only in terms of $n$ and the $n$-independent coefficient functions without $h_0$ in Corollary \ref{Coro: aggregated Hall 1991} 
and \ref{coro:with PI}. Using these results, we scrutinise the higher-order difference between the theoretical and plug-in bandwidths in Section \ref{subsection:Difference}. We realise that the global plug-in bandwidth selection starts to have an 
impact from on the order of $O\{(nh_0)^{-1/2}h_0+h_0^{L+1}\}=O(h_0^{L+1})$, which is stated in Theorem \ref{Coro: Exact Order}. 
Section \ref{section:special case} provides a comprehensive example by considering the special case of $L=2$.

In section \ref{section:Edgeworth Expansion Including Pilot Bandwidth and Studentisation}, we develop the Edgeworth expansion for any $L_p$. This expansion implies that, when $L_p$ is small, bandwidth selection by the global plug-in method has an effect on the asymptotic structure up to the order of \\ $O_p\{(nh_0)^{-1/2}+h_0^L\}$. In addition, we study the intersectional effect of the plug-in method, pilot bandwidth, and the Studentisation.

\subsection{Edgeworth Expansion for KDE with Global Plug-In Bandwidth up to the order of $O\{(nh_0)^{-1}\}$} \label{subsec: second order Edgeworth Expansion}

We introduce the following assumption:
\begin{assumption} \label{A:Ku^2}
For $1\leqq k\leqq L-1$, $\lim_{u\rightarrow\pm\infty}|K(u)u^{k}|\rightarrow0 
\text{ and }
\lim_{u\rightarrow\pm\infty}|K'(u)u^2|\rightarrow0$.
\end{assumption}

\begin{assumption} \label{A:L_p order 2}
$L_p > 8L - \frac{1}{2}$
\end{assumption}
\begin{remark}
     Assumption \ref{A:L_p order 2} guarantees that the pilot bandwidth has no effect on the asymptotic structures of KDE up to the order of $O_p\{(nh_0)^{-1}+h_0^{2L}\}$. We obtain this assumption from Edgeworth expansion including the pilot bandwidth (Theorem \ref{T:Edgeworth Expansion Including Pilot Bandwidth}). Although this assumption may be  unrealistic for empirical analysis (for example, when one uses $L=2$, $L_p\ge 16$ is necessary), it can be considered of theoretical value in the sense that this theorem clarifies the inevitable effect of bandwidth selection via the global plug-in method (i.e. estimating $\int f''(x)^2dx$). For empirical application, Theorem \ref{T:Edgeworth Expansion Including Pilot Bandwidth and the effect of Studentisation}, which considers the simultaneous effect of the bandwidth selection via the global plug-in method, its accompanying pilot bandwidth, and Studentisation,  is more valuable.
\end{remark}

We have the following theorem which is proved in \ref{Proof:Main Theorem}.

\begin{theorem}[The Effect of Estimation of $\int f^{(L)}(x)^2dx$] \label{T:Main Theorem}
Under Assumptions \ref{A:DGP}, \ref{A: interior point}, \ref{A:rate b}, \ref{A:2L+L_p times differentiable}, \ref{A:pilot kernel}, \ref{A:Cramer}, \ref{A: Kernel ''}, \ref{A:Ku^2} and \ref{A:L_p order 2}, the following expansion is valid:
    \begin{align}
        \sup_{z \in \mathbb{R}} &\Biggl|\mathbb{P}(S_{\text{PI}}(x)\leqq z)-\Phi(z) \nonumber\\
        &-\phi(z)\Biggl[(nh_0)^{-1/2}p_1(z)+\sum_{l=0}^{L-1}h_0^{L+l+1}p_{3,l}(z)+n^{-1/2}h_0^{1/2}p_4(z)+(nh_0)^{-1}p_2(z)\Biggl]\Biggl|\nonumber\\
        &=o\{(nh_0)^{-1}+h_0^{2L}\}, \label{EE Third with PI}
    \end{align}
    where
    \begin{align}
        p_{3,l}(z)&=-C_{PI}C_{\Gamma,l}(x)\rho_{11}\mu_{20}^{-1}z, \nonumber\\
        p_4(z)&=-C_{PI}\rho_{11}\xi_{11}\mu_{20}^{-3/2}(z^2-1)+\frac{1}{2}C_{PI}\rho_{11}\mu_{20}^{-1/2}z^2, \nonumber\\
        C_{\Gamma,l}(x)&= -\left(\int u^{L+l}K(u)du\right)\frac{f^{(L+l)}(x)}{(L+l-1)!}, \nonumber
    \end{align}
    \begin{align}
        \xi_{kl}&=h_0^{-\mathbbm{1}(\{k\geqq 1\}~\cup~ \{l\geqq 1\})}\mathbb{E}\Biggl[\left\{K_{i,h_0}(x)-\mathbb{E}[K_{i,h_0}(x)]\right\}^k\nonumber\\
        &\times\left\{K'_{i,h_0}(x)u_{i,h_0}(x)+K_{i,h_0}(x)-\mathbb{E}[K'_{i,h_0}(x)u_{i,h_0}(x)+K_{i,h_0}(x)]\right\}^l\Biggl], \nonumber
    \end{align}
    \begin{equation}
        \rho_{kl}=h_0^{-\mathbbm{1}(k\geqq 1)}\mathbb{E}\left[\left\{K_{i,h_0}(x)-\mathbb{E}[K_{i,h_0}]\right\}^k\Biggl\{f^{(2L)}(X_i)-\mathbb{E}[f^{(2L)}(X_i)]\Biggl\}^l\right]. \nonumber
    \end{equation}
\end{theorem}
\begin{remark}
     The effects of bandwidth selection emerging in this theorem come from $\Gamma_{KDE_1}$ (this effect is independent of the bandwidth selection method) and from the largest component of the projection term of Hoeffding-decomposition of $\hat{I}_L$. Since the head term of the projection term of decomposed $\hat{I}_L^{convo}$, the result of this theorem does not change, even if the estimator is changed from $\hat{I}_L$ to $\hat{I}_L^{convo}$.
\end{remark}

\subsection{Edgeworth Expansions in Powers of $n^{-1/(2L+1)}$} \label{section: comparison}

\cite{Hall91} does not specify the bandwidth order, but we consider the use of plug-in bandwidth with a fixed convergence rate. Note that $h_0$ satisfies Assumption \ref{A: nh/long infty}. Comparing (\ref{EE Third with OPT}) with $h=h_0$ and (\ref{EE Third with PI}), we see that $h_0^{L+l+1}p_{3,l}(z)$ and $n^{-1/2}h_0^{1/2}p_4(z)$ reflect the effect of bandwidth selection via global plug-in methods.  However, the results in Proposition \ref{p:Edgeworth Expansion for KDE with optimal bandwidth}, Theorem \ref{T: Hall and Kang}, and Theorem \ref{T:Main Theorem} are still insufficient for identifying the exact difference because $\mu_{kl},\rho_{kl},\xi_{kl}$. Accordingly $p_1(z),p_2(z),p_{3,l}(z)$, and $p_4(z)$ in the expansions depend on $h_0$ and, consequently, the relationship between the terms in the expansions is unclear.

For $S_{h_0}(x)$, we have to expand $p_1(z)$ and $p_2(z)$ in terms of only $n$, without $h_0$. $p_2(z)$ is easy to handle because only its leading term affects the Edgeworth expansion up to the order of $O\{(nh_0)^{-1}\}$. For $p_1(z)$, recalling that $p_1(z)=-\frac{1}{6}\mu_{20}^{-3/2}\mu_{30}(z^2-1)$, we  expand $\mu_{30}\mu_{20}^{-3/2}$ up to the term whose convergence rate is $O(h_0^L)$. Letting $\kappa_{st}\equiv\int u^sK(u)^tdu$ and, from straightforward computation, we can expand $\mu_{20},\mu_{30}$ as follows:
\begin{align}
    \mu_{20}&=\kappa_{02}f(x)-f(x)^2h_0+\sum_{l=2}^L\kappa_{l2}\frac{f^{(l)}(x)}{l!}h_0^l+o(h_0^L), \label{linearize mu20}\\
    \mu_{30}&=\kappa_{03}f(x)-3\kappa_{02}f(x)^2h_0+\left\{\kappa_{23}\frac{f^{(2)}(x)}{2!}+2f(x)^3\right\}h_0^2 \nonumber\\
    &+\sum_{l=3}^L\left\{\kappa_{l3}\frac{f^{(l)}(x)}{l!}-3\kappa_{l-1,2}\frac{f^{(l-1)}(x)}{(l-1)!}f(x)\right\}h_0^l+o(h_0^L). \label{linearize mu30}
\end{align}
For notational simplicity, we rewrite $\mu_{20},\mu_{30}$ as a series of $h_0$:
\begin{align}
    \mu_{20}\equiv \sum_{l=0}^Lm_{2,l}(x)h_0^l+o(h_0^L),\qquad\mu_{30}\equiv \sum_{l=0}^Lm_{3,l}(x)h_0^l+o(h_0^L). \nonumber
\end{align}
Then, expanding $\mu_{30}\mu_{20}^{-3/2}$ yields:
\begin{align}
    \mu_{30}\mu_{20}^{-3/2}&=\left\{\sum_{l=0}^Lm_{3,l}(x)h_0^l+o(h_0^L)\right\}\left\{\sum_{j=0}^Lm_{2,j}(x)h_0^j+o(h_0^L)\right\}^{-3/2} \nonumber\\
    &=\left\{\sum_{l=0}^Lm_{3,l}(x)h_0^l+o(h_0^L)\right\} \nonumber\\
    &\qquad\times\Biggl[\sum_{k=0}^L\frac{(-1)^k(2k+1)!!}{2^k k!}\{m_{2,0}(x)\}^{\frac{-(2k+3)}{2}}\left(\sum_{j=1}^{L}m_{2,j}(x)h_0^j\right)^k+o(h_0^L)\Biggl] \nonumber\\
    &=\sum_{l=0}^L\sum_{k=0}^l\frac{(-1)^k(2k+1)!!}{2^kk!}\Bigl\{m_{2,0}(x)\Bigl\}^{\frac{-(2k+3)}{2}} \nonumber\\
    &\qquad\times\underset{k\leqq i_1+\cdots+i_k+l\leqq L}{\sum\cdots\sum}m_{3,l}(x)m_{2,i_1}(x)\cdots m_{2,i_k}(x)h_0^{i_1+\cdots i_k+l}+o(h_0^L). \nonumber
\end{align}
We define $\gamma_{1,0}(x),\gamma_{1,1}(x),\gamma_{2,1,0}(x)$ and $\gamma_{2,2,0}(x)$ as follows:
\begin{align}
    \gamma_{1,0}(x)&\equiv\frac{-1}{6}\kappa_{02}^{-3/2}\kappa_{03}f(x), \nonumber\\
    \gamma_{1,1}(x)&\equiv\frac{1}{2}\Bigl(\kappa_{02}^{-1/2}f(x)^{1/2}-\kappa_{02}^{-5/2}\kappa_{03}f(x)^{1/2}/2\Bigl), \nonumber\\
    \gamma_{2,1,0}(x)&\equiv\frac{-1}{24}\kappa_{02}^{-2}\kappa_{04}f(x)^{-1}, \nonumber\\
    \gamma_{2,2,0}(x)&\equiv\frac{-1}{72}\kappa_{02}^{-3}\kappa_{03}^2f(x)^{-1}. \nonumber
\end{align}
From the above results, we obtain the following corollary.
\begin{corollary}[Expansion of \cite{Hall91} in powers of $n^{-1/(2L+1)}$] \label{Coro: aggregated Hall 1991}
Under Assumptions \ref{A:DGP}, \ref{A: Kernel}, \ref{A: interior point}, and \ref{A:Cramer}:
\begin{align}
    \sup_{z\in\mathbb{R}}\Biggl|\mathbb{P}(S_{h_0}(x)\leqq z)&-\Phi(z)-\phi(z)\sum_{j=0}^L a_j(z,x)n^{\frac{-(L+j)}{2L+1}}\Biggl|=o\{(nh_0)^{-1}+h_0^{2L}\}, \nonumber
\end{align}
where the definitions of $a_j(z,x)$ are given as follows for $2\leqq q\leqq L-1$:
\begin{align}
    a_0(z,x)&=\gamma_{1,0}(x)(z^2-1), \nonumber\\
    a_1(z,x)&=\gamma_{1,1}(x)(z^2-1), \nonumber\\
    a_q(z,x)&=\sum_{l=0}^L\sum_{k=0}^l\frac{(-1)^k(2k+1)!!}{2^kk!}\Bigl\{m_{2,0}(x)\Bigl\}^{\frac{-(2k+3)}{2}} \nonumber\\
    &\qquad\qquad\times\underset{i_1+\cdots+i_k+l= q}{\sum\cdots\sum}m_{3,l}(x)m_{2,i_1}(x)\cdots m_{2,i_k}(x)h_0^{i_1+\cdots i_k+l}(z^2-1), \nonumber\\
    a_L(z,x)&=\gamma_{2,0,1}(x)(z^3-3z)+\gamma_{2,0,2}(x)(z^5-10z+15) \nonumber\\
    &\qquad + \sum_{l=0}^L\sum_{k=0}^l\frac{(-1)^k(2k+1)!!}{2^kk!}\Bigl\{m_{2,0}(x)\Bigl\}^{\frac{-(2k+3)}{2}} \nonumber\\
    &\qquad\qquad\times\underset{i_1+\cdots+i_k+l= L}{\sum\cdots\sum}m_{3,l}(x)m_{2,i_1}(x)\cdots m_{2,i_k}(x)h_0^{i_1+\cdots i_k+l}(z^2-1). \nonumber
\end{align}
\end{corollary}
\begin{remark}
     Note that $a_0$ and $a_1$ are special cases of $a_q$, but we explicitly write these terms for comparison of this result with the next corollary. 
\end{remark}

\begin{remark}
     From this corollary, we identify that the Edgeworth expansion of the Standardised KDE with deterministic bandwidth has the term of order $O\{(nh_0)^{-1/2}\}=O(n^{\frac{-L}{2L+1}})$ right after the term $\Phi(z)$, with a gap between them, but the subsequent terms decrease at the rate of $O(h_0)=O(n^{\frac{-1}{2L+1}})$, which is not clear in \cite{Hall91}.
\end{remark}

Next, for (\ref{EE Third with PI}), we also have to expand $p_{3,l}(z)$ and $p_4(z)$. Although we do not provide the details here, one can use a similar process for $p_1(z)$. We define:
\begin{align}
    \tau_{l}&\equiv\int u^l\{K(u)K'(u)u+K(u)^2\}du, \nonumber\\
    \mathcal{L}(x)&\equiv f^{(2L)}(x)-\mathbb{E}[f^{(2L)}(x)], \nonumber\\
    \gamma_{3,1,0}(x)&\equiv -C_{PI}C_{\Gamma,0}(x)\kappa_{02}^{-1}\mathcal{L}(x), \quad
    \gamma_{3,1,1}(x)\equiv C_{PI}C_{\Gamma,0}(x)\kappa_{02}^{-2}\mathcal{L}(x) f(x) \nonumber\\
    \gamma_{4,1,0}(x)&\equiv -C_{PI}\kappa_{02}^{-3/2}\tau_0\mathcal{L}(x) f^{1/2}(x), \quad
    \gamma_{4,1,1}(x)\equiv \frac{3}{2}C_{PI}\kappa_{02}^{-5/2}\tau_0\mathcal{L}(x) f(x)^{3/2}. \nonumber\\
    \gamma_{4,2,0}(x)&\equiv \frac{1}{2}C_{PI}\kappa_{02}^{-1/2}\mathcal{L}(x)f(x)^{1/2}, \quad
    \gamma_{4,2,1}(x)\equiv \frac{-1}{4}C_{PI}\kappa_{02}^{-3/2}\mathcal{L}(x)f(x)^{3/2} \nonumber
\end{align}
Then, we have the next corollary.

\begin{corollary}[Main Theorem in powers of $n^{-1/(2L+1)}$] \label{coro:with PI} Under the same assumptions as in Theorem \ref{T:Main Theorem}:
\begin{equation}
    \sup_{z\in\mathbb{R}}\left|\mathbb{P}(S_{PI}(x)\leqq z)-\Phi(z)-\phi(z)\sum_{j=0}^Lb_j(z,x)n^{\frac{-(L+j)}{2L+1}}\right|=o\{(nh_0)^{-1}+h_0^{2L}\} \nonumber,
\end{equation}
where
\begin{align}
    b_0(z,x)&=a_0(z,x), \nonumber\\
    b_1(z,x)&=a_1(z,x)+\gamma_{3,1,0}(x)z+\gamma_{4,1,0}(x)(z^2-1)+\gamma_{4,2,0}(x)z^2. \nonumber\\
    b_2(z,x)&=a_2(z,x)+\gamma_{3,1,1}(x)z+\gamma_{4,1,1}(x)(z^2-1)+\gamma_{4,2,1}(x)z^2. \nonumber
\end{align}
\end{corollary}
\noindent Here, we do not provide the definitions of $ b_j(z,x),~j=3,...$ because they are too lengthy and tedious, but they can be obtained in a straightforward manner. 

\cite{HK01} and Theorem \ref{T: Hall and Kang} state that the global plug-in method has no effect on the terms up to whose convergence rates are $O\{(nh_0)^{-1/2}\}$; in other words, $b_0(z,x)$ does not include the effect of bandwidth selection in view of Corollary \ref{coro:with PI}. Comparing $a_1(z,x)$ and $b_1(z,x)$, the bandwidth selection via the global plug-in method starts to have an effect on the term with the order of $O\{(nh_0)^{-1/2}h_0\}=O(n^{\frac{-(L+1)}{2L+1}})$. The deviation between $b_0(z,x)$ (the smallest term not affected by bandwidth selection) and $b_1(z,x)$ (the largest term affected by bandwidth selection) is only of order $O(h_0)=O(n^{-1/(2L+1)})$.

\begin{remark}
     Although we omit $b_3(z,x)$ and the subsequent terms, we can show that these terms are also affected by bandwidth selection via the global plug-in method in the same way as the process of deriving Corollary \ref{Coro: aggregated Hall 1991}. However, the most important point is that the influence of the bandwidth selection via the global plug-in method starts to appear at $b_1(z,x)$.
\end{remark}

\subsection{Difference between $S_{h_0}(x)$ and $S_{PI}(x)$} \label{subsection:Difference}
From Corollaries \ref{Coro: aggregated Hall 1991} and \ref{coro:with PI}, we can easily deduce the following theorem, which states the exact order of the difference between $S_{h_0}(x)$ and $S_{PI}(x)$. See Appendix \ref{Proof of Exact Order} for the proof.
\begin{theorem}[Exact Evaluation of the Deviation] \label{Coro: Exact Order} Under the same assumptions as in Theorem \ref{T:Main Theorem}:
\begin{align}
    &\sup_{z\in\mathbb{R}}\Bigl|\mathbb{P}(S_{h_0}(x)\leqq z)-\mathbb{P}(S_{PI}(x)\leqq z)\nonumber\\
    &\qquad-\phi(z)\Bigl[\{\gamma_{3,1,0}(x)z+\gamma_{4,1,0}(x)(z^2-1)+\gamma_{4,2,0}(x)z^2\}n^{\frac{-(L+1)}{2L+1}}\Bigl]\Bigl|=O(n^{\frac{-(L+2)}{2L+1}}), \nonumber
\end{align}
and the order is exact.
\end{theorem}
\noindent This theorem implies that:
\begin{equation}
    \sup_{z\in\mathbb{R}}|\mathbb{P}(S_{h_0}(x)\leqq z)-\mathbb{P}(S_{PI}(x)\leqq z)|=O(n^{\frac{-(L+1)}{2L+1}}). \nonumber
\end{equation}
We can only claim that this deviation is $o\{(nh_0)^{-1/2}\}=o(n^{-L/(2L+1)})$ from Theorem 3.2 in \cite{HK01} and our Theorem \ref{T: Hall and Kang}, whereas Theorem \ref{Coro: Exact Order} gives a stronger result, stating that the convergence rate is exactly $O\{(nh_0)^{-1/2}h_0\}=O(n^{-(L+1)/(2L+1)})$. 

\begin{remark}
     The larger the kernel order $L$ we use, the slower the convergence rate of the approximation in Theorems \ref{T: Hall and Kang}, \ref{T:Main Theorem}, and \ref{Coro: Exact Order} will be. This is because we centralise at $\mathbb{E}\hat{f}_{h_0}(x)$. However, as stated in Section \ref{Section:Discussion}, one of the final goals would be to examine the effect of bandwidth selection and `debias’ simultaneously (we are in the process of working on it), and it is unclear if the second-order kernel $L=2$ is optimal. 
\end{remark}

\subsection{Special Case} \label{section:special case}
Since the previous results are difficult to interpret because of their generality, we consider a special case of $L=2$. Here, we also provide the details of the expansions of $p_{3,l}(z)$ and $p_4(z)$ as well as that of $p_1(z)$.

First, we have to expand $p_1(z)$ and $p_2(z)$. From (\ref{linearize mu20}) and (\ref{linearize mu30}), we can expand $p_1(z)$ as follows (see \ref{Appendix C}):
\begin{align}
    p_1(z)&=\frac{-1}{6}\mu_{30}\mu_{20}^{-3/2}(z^2-1) \nonumber\\
    &=\frac{-1}{6}\Biggl[\kappa_{02}^{-3/2}\kappa_{03}f(x)-3\Biggl\{\frac{f(x)^{1/2}}{\kappa_{02}^{1/2}}-\frac{\kappa_{03}f(x)^{1/2}}{2\kappa_{02}^{5/2}}\Biggl\}h_0 \nonumber\\
    &\qquad+\Biggl\{\frac{-3}{4}\{\kappa_{02}f(x)\}^{-5/2}\kappa_{03}\kappa_{23}f^{(2)}(x)f(x)-3\{\kappa_{02}f(x)\}^{-5/2}\kappa_{03}f(x)^4 \nonumber\\
    &\qquad\qquad+\frac{15}{8}\{\kappa_{02}f(x)\}^{-7/2}\kappa_{03}f(x)^5+\frac{9}{2}\kappa_{02}^{-3/2}f(x)^{3/2}\Biggl\}h_0^2\Biggl](z^2-1)+o(h_0^2), \nonumber\label{Linearized p_1}
\end{align}
and since for $p_2(z)$ we need only the leading term; a straightforward computation yields:
\begin{equation}
    p_2(z)=\frac{-1}{24}\kappa_{02}^{-2}\kappa_{04}f(x)^{-1}(z^3-3z)-
    \frac{1}{72}\kappa_{02}^{-3}\kappa_{03}^2f(x)^{-1}(z^5-10z^3+15)+o(1). \nonumber
\end{equation}
From the above results, in the special case of $L=2$, expansion (\ref{EE Third with OPT}) is as follows:
\begin{align}
    &\sup_{z\in\mathbb{R}}\Biggl|\mathbb{P}(S(x)\leqq z)-\Phi(z)-\phi(z)\Bigl[a_0(z,x)n^{-2/5}+a_1(z,x)n^{-3/5}+a_2(z,x)n^{-4/5}\Bigl]\Biggl| \nonumber\\
    &=o\{(nh_0)^{-1}+h_0^4\}, \nonumber
\end{align}
where
\begin{align}
    a_0(z,x)&=\gamma_{1,0}(x)(z^2-1), \nonumber\\
    a_1(z,x)&=\gamma_{1,1}(x)(z^2-1), \nonumber\\
    a_2(z,x)&=\gamma_{1,2}(x)(z^2-1)+\gamma_{2,1,0}(x)(z^3-3z)+\gamma_{2,2,0}(x)(z^5-10z^3+15). \nonumber
\end{align}
Next, we expand $p_{3,0}(z),p_{3,1}(z)$ and $p_4(z)$. From a straightforward computation, noting $\tau_1=0$ from the properties of the odd function, we can expand $\rho_{11}$ and $\xi_{11}$ as follows: 
\begin{align}
    \rho_{11}&=\mathcal{L}(x) f(x)+O(h_0^L), \nonumber\\
    \xi_{11}&=\tau_0f(x)+\tau_1f^{(1)}(x)h_0+o(h_0)=\tau_0+o(h_0). \nonumber
\end{align}
These imply:
\begin{align}
    p_{3,0}(z)&=-C_{PI}C_{\Gamma,0}(x)\rho_{11}\mu_{20}^{-1}z \nonumber\\
    &=-C_{PI}C_{\Gamma,0}(x)\kappa_{02}^{-1}\mathcal{L}(x) z+C_{PI}C_{\Gamma,0}(x)\kappa_{02}^{-2}\mathcal{L}(x) f(x)zh_0+o(h_0). \nonumber
\end{align}
See \ref{Appendix C} for the second equality. Noting that $C_{\Gamma,1}(x)=0$ from the properties of the odd function:
\begin{align}
    p_{3,1}(z)&=-C_{PI}C_{\Gamma,1}(x)\rho_{11}\mu_{20}^{-1}z=0, \nonumber
\end{align}
and, as shown in \ref{Appendix C}:
\begin{align}
    p_4(z)&=-C_{PI}\rho_{11}\xi_{11}\mu_{20}^{-3/2}(z^2-1)+\frac{1}{2}C_{PI}\rho_{11}\mu_{20}^{-1/2}z^2 \nonumber\\
    &=-C_{PI}\kappa_{02}^{-3/2}\tau_0\mathcal{L}(x) f(x)^{1/2}(z^2-1)+\frac{3}{2}C_{PI}\kappa_{02}^{-5/2}\tau_0\mathcal{L}(x) f(x)^{3/2}(z^2-1)h_0+o(h_0) \nonumber\\
    &\quad+\frac{1}{2}C_{PI}\kappa_{02}^{-1/2}\mathcal{L}(x)f(x)^{1/2}z^2-\frac{1}{4}C_{PI}\kappa_{02}^{-3/2}\mathcal{L}(x)f(x)^{3/2}z^2h_0+o(h_0). \nonumber
\end{align}
From the above results, in the special case of $L=2$, the expansion (\ref{EE Third with PI}) is as follows.
\begin{align}
    &\sup_{z\in\mathbb{R}}\left|\mathbb{P}(S_{PI}(x)\leqq z)-\Phi(z)-\phi(z)\Bigl[b_0(z,x)n^{-2/5}+b_1(z,x)n^{-3/5}+b_2(z,x)n^{-4/5}\Bigl]\right| \nonumber\\
    & \qquad =o\{(nh_0)^{-1}+h_0^{4}\}, \nonumber 
\end{align}
where the definitions of $b_0(z,x)$, $b_1(z,x)$ are given as follow.
\begin{align}
    b_0(z,x)&=a_0(z,x) \nonumber\\
    b_1(z,x)&=a_1(z,x)+\gamma_{3,1,0}(x)z+\gamma_{4,1,0}(x)(z^2-1)+\gamma_{4,2,0}(x)z^2\nonumber\\
    b_2(z,x)&=a_2(z,x)+\gamma_{3,1,1}(x)z+\gamma_{4,1,1}(x)(z^2-1)+\gamma_{4,2,1}(x)z^2 \nonumber
\end{align}

\subsection{Edgeworth Expansion Including Pilot Bandwidth and Studentisation} \label{section:Edgeworth Expansion Including Pilot Bandwidth and Studentisation}
In this section, we provide two more expansions. One is the Edgeworth Expansion of Standardised KDE with global plug-in bandwidth and its accompanying pilot bandwidth. Here, we allow $L_p$ to be small so that $b$ affects the expansion (note that Theorems \ref{T: Hall and Kang} and \ref{T:Main Theorem} set $L_p$ sufficiently large such that $b$ does not appear in the expansion). The other is the Edgeworth Expansion of Studentised KDE with the global plug-in bandwidth and its accompanying pilot bandwidth.

Let $H^{(2L)}\left(\frac{X_i-X_j}{b}\right)=H^{(2L)}_{ij,b}$ and define 
\begin{align*}
    & \omega_{111}\equiv h_0^{-1}b^{-1}\mathbb{E}\Biggl[\left\{K_{i,h_0}(x)-\mathbb{E}[K_{i,h_0}(x)]\right\}\left\{K_{i,h_0}(x)-\mathbb{E}[K_{i,h_0}(x)]\right\}\\
    &\qquad\qquad\times\left\{H^{(2L)}_{ij,b}-\mathbb{E}\left[H^{(2L)}_{ij,b}|X_i\right]-\mathbb{E}\left[H^{(2L)}_{ij,b}|X_j\right]+\mathbb{E}\left[H^{(2L)}_{ij,b}\right]\right\}\Biggl]\\
    & \psi_{111}\equiv h_0^{-1}b^{-1}\mathbb{E}\Biggl[\left\{K_{i,h_0}(x)-\mathbb{E}\left[K_{i,h_0}(x)\right]\right\} \left\{K_{i,h_0}(x)+K'_{i,h_0}(x)u_{i,h_0}(x)-\mathbb{E}[K_{i,h_0}(x)+K'_{i,h_0}(x)u_{i,h_0}(x)]\right\}\\
    &\qquad\qquad\times\left\{H^{(2L)}_{ij,b}-\mathbb{E}\left[H^{(2L)}_{ij,b}|X_i\right]-\mathbb{E}\left[H^{(2L)}_{ij,b}|X_j\right]+\mathbb{E}\left[H^{(2L)}_{ij,b}\right]\right\}\Biggl]
\end{align*}

\begin{assumption} \label{A:L_p order 3}
    $L_p > \frac{2L}{5} + \frac{1}{10}$
\end{assumption}

\begin{assumption} \label{A:L_p order 4}
     $L_p > \frac{4L}{3} + 1$
\end{assumption}

\begin{remark}
     Assumption \ref{A:L_p order 3} is for the following  expansion  (\ref{eq:EE with pilot 1}) and  Assumption \ref{A:L_p order 4} is for the following  expansion  (\ref{eq:EE with pilot 2}). Owing to these assumptions, we can assume $[(\hat{h}-h_0)/h_0]^2$ in (\ref{def SPI}) is negligible. Without this assumption, when $L$ is large enough, Edgeworth expansion (\ref{eq:EE with pilot 1}) and (\ref{eq:EE with pilot 2}) have the term associated with $[(\hat{h}-h_0)/h_0]^3, [(\hat{h}-h_0)/h_0]^4,\dots$. However, when $L=2$, $L_p=4$ is sufficient for Assumption \ref{A:L_p order 4} and any $L_p$ satisfies Assumption \ref{A:L_p order 3}, so these assumptions are not unrealistic unlike Assumption \ref{A:L_p order 2}.
\end{remark}

\begin{remark}
     We do not provide the mathematically rigorous proof for the following expansions. However, one can prove their validity of them in the same way as our Theorem \ref{T:Main Theorem}.
\end{remark}

\begin{theorem}[Edgeworth Expansion Including Pilot Bandwidth] \label{T:Edgeworth Expansion Including Pilot Bandwidth}
Under Assumptions \ref{A:DGP}, \ref{A: interior point}, \ref{A:rate b}, \ref{A:2L+L_p times differentiable}, \ref{A:pilot kernel}, \ref{A: Kernel ''}, \ref{A:Ku^2} and \ref{A:L_p order 3}
\begin{align}
    \mathbb{P}(S_{PI}(x)\leqq z) & = \Phi(z) + \phi(z)\Biggl[p_1(z)(nh_0)^{-1/2} + \sum_{l=0}^{L-1}\mathfrak{p}_{1,l}(z)n^{-1/2}h_0^{(2L+2l+1)/2}b^{-2L} \nonumber\\
    &\qquad+\mathfrak{p}_{2}(z)n^{-1}b^{-2L}\Biggl] + o\{(nh_0)^{-1/2}+h_0^L\} \label{eq:EE with pilot 1}\\
    \mathbb{P}(S_{PI}(x)\leqq z) & = \Phi(z) + \phi(z)\Biggl[p_1(z)(nh_0)^{-1/2}+p_2(z)(nh_0)^{-1}+\sum_{l=0}^{L-1}p_{3,l}(z)h_0^{L+l+1} \nonumber\\
    & \qquad + \sum_{l=0}^{L-1}\mathfrak{p}_{1,l}(z)n^{-1/2}h_0^{(2L+2l+1)/2}b^{-2L}+p_4(z)n^{-1/2}h_0^{1/2}+\mathfrak{p}_{2}(z)n^{-1}b^{-2L}\Biggl] \nonumber\\
    &\qquad + o\{(nh_0)^{-1}+h_0^{2L}\}\label{eq:EE with pilot 2}
\end{align}
where the definitions of $\mathfrak{p}_{1,l}(z)$ and $\mathfrak{p}_{2}(z)$ are
\begin{align*}
    & \mathfrak{p}_{1,l}(z) \equiv -\frac{C_{PI}C_{\Gamma,l}(x)}{2}\mu_{20}^{-3/2}\omega_{111}(z^2-1)\\
    & \mathfrak{p}_{2}(z) \equiv -C_{PI}\left(\frac{1}{2}\mu_{20}^{-2}\xi_{11}\omega_{111}(z^3-3z)+\mu_{20}^{-1}\psi_{111}z-\frac{1}{4}\mu_{20}^{-1}\omega_{111}(z^3-z)\right)
\end{align*}
\end{theorem}
\begin{remark}
     $\mathfrak{p}_{1,l}(z)$ and $\mathfrak{p}_{2}(z)$ reflect the effect of pilot bandwidth.
\end{remark}

\begin{remark}
     When one uses $\hat{I}_L^{convo}$ instead of $\hat{I}_L$, the definitions of $\eta_k, \omega_{111}$ and $\psi_{111}$ are changed as follows.
     \begin{align*}
         & \omega_{111}^{convo}\equiv h_0^{-1}b^{-1}\mathbb{E}\Biggl[\left\{K_{i,h_0}(x)-\mathbb{E}[K_{i,h_0}(x)]\right\}\left\{K_{i,h_0}(x)-\mathbb{E}[K_{i,h_0}(x)]\right\}\\
         &\qquad\qquad\times\left\{\bar{H}^{(L)}_{ij,b}-\mathbb{E}\left[\bar{H}^{(L)}_{ij,b}|X_i\right]-\mathbb{E}\left[\bar{H}^{(L)}_{ij,b}|X_j\right]+\mathbb{E}\left[\bar{H}^{(L)}_{ij,b}\right]\right\}\Biggl]\\
        & \psi_{111}^{convo}\equiv h_0^{-1}b^{-1}\mathbb{E}\Biggl[\left\{K_{i,h_0}(x)-\mathbb{E}\left[K_{i,h_0}(x)\right]\right\} \left\{K_{i,h_0}(x)+K'_{i,h_0}(x)u_{i,h_0}(x)-\mathbb{E}[K_{i,h_0}(x)+K'_{i,h_0}(x)u_{i,h_0}(x)]\right\}\\
        &\qquad\qquad\times\left\{\bar{H}^{(L)}_{ij,b}-\mathbb{E}\left[\bar{H}^{(L)}_{ij,b}|X_i\right]-\mathbb{E}\left[\bar{H}^{(L)}_{ij,b}|X_j\right]+\mathbb{E}\left[\bar{H}^{(L)}_{ij,b}\right]\right\}\Biggl]
     \end{align*}
where $\bar{H}^{(L)}(u) \equiv \int H^{(L)}(v)H^{(L)}(u-v)dv$.
\end{remark}

The following Theorem \ref{T:Edgeworth Expansion Including Pilot Bandwidth and the effect of Studentisation} is a formal expansion of the Studentised KDE with global plug-in bandwidth. Although we have options for variance estimation, we adopt the following natural estimator as \cite{Hall91,Hall92book} and \cite{HK01}.
\begin{align}
    \hat{\mu}_{20}(h)\equiv h^{-1}\left\{\frac{1}{n}\sum_{i=1}^nK_{i,h}(x)^2-\left(\frac{1}{n}\sum_{i=1}^nK_{i,h}(x)\right)^2\right\}. \nonumber
\end{align}
In addition, as \cite{HK01}, let bandwidth used for the estimation of $\mu_{20}$ be $\hat{h}$.
Concequently, Studentised KDE with global plug-in bandwidth is given by
\begin{align}
    T_{PI}(x) = \frac{\sqrt{n\hat{h}}(\hat{f}_{\hat{h}}(x)-\mathbb{E}\hat{f}_{h_0}(x))}{\hat{\mu}_{20}(\hat{h})^{1/2}} \nonumber
\end{align}

Define
\begin{align}
    \delta \equiv h_0^{-1}\left\{\mathbb{E}\left[K_{i,h_0}(x)K'_{i,h_0}(x)u_{i,h_0}(x)\right]-\mathbb{E}\left[K_{i,h_0}(x)\right]\left[K'_{j,h_0}(x)u_{j,h_0}(x)\right]\right\} \nonumber
\end{align}

\begin{theorem}[Edgeworth Expansion Including Pilot Bandwidth and the effect of Studentisation] \label{T:Edgeworth Expansion Including Pilot Bandwidth and the effect of Studentisation}
under Assumptions \ref{A:DGP}, \ref{A: interior point}, \ref{A:rate b}, \ref{A:2L+L_p times differentiable}, \ref{A:pilot kernel}, \ref{A: Kernel ''}, \ref{A:Ku^2} and \ref{A:L_p order 4}
\begin{align*}
    \mathbb{P}(T_{PI}(x)\leqq z) & = \Phi(z) + \phi(z)\Biggl[q_1(z)(nh_0)^{-1/2}+\Bigl\{q_2(z)+p_4(z)+\mathfrak{q}_1(z)\Bigl\}n^{-1/2}h_0^{1/2}+q_3(z)(nh_0)^{-1}\\
    & \quad + \sum_{l=0}^{L-1}p_{3,l}(z)h_0^{L+l+1} + \sum_{l=0}^{L-1}\mathfrak{p}_{1,l}(z)n^{-1/2}h_0^{(2L+2l+1)/2}b^{-2L} + \Bigl\{\mathfrak{p}_2(z)+\mathfrak{q}_2(z)\Bigl\}n^{-1}b^{-2L}\Biggl] \\
    & \quad + o_p\{(nh_0)^{-1}+h_0^{2L}\}
\end{align*}
where the definitions of $q_1(z), q_2(z), q_3(z), \mathfrak{q}_1(z), \mathfrak{q}_{2}(z)$ and $\mathfrak{q}_{3}(z)$ are
\begin{align*}
    & q_1(z) \equiv \frac{1}{2}\mu_{20}^{-3/2}\mu_{11}-\frac{1}{6}\mu_{20}^{-3/2}(\mu_{30}-3\mu_{11})(z^2-1)\\
    & q_2(z) \equiv -f(x)\mu_{20}^{-1}z^2\\
    & q_3(z) \equiv -\mu_{20}^{-3}\mu_{30}^2z - (\frac{2}{3}\mu_{20}^{-3}\mu_{30}^2 - \frac{1}{12}\mu_{20}^{-2}\mu_{40})(z^3-3z) - \frac{1}{18}\mu_{20}^{-3}\mu_{30}^3(z^5-10z^3+15z)\\
    & \mathfrak{q}_1(z) \equiv \frac{C_{PI}}{2}\Bigl\{1+\delta\mu_{20}^{-1}\Bigl\}\mu_{20}^{-1/2}\rho_{11} z^2\\
    & \mathfrak{q}_2(z) \equiv \frac{C_{PI}}{4}\mu_{20}^{-1}\omega_{111}\Bigl\{1+\delta\mu_{20}^{-1}\Bigl\} (z^3-2z)
\end{align*}
\end{theorem}

\begin{remark}
     $q_1(z), q_2(z)$ and  $q_3(z)$ reflect the effect of Studentisation and $\mathfrak{q}_1(z)$ and $\mathfrak{q}_{2}(z)$ reflect the simultaneous effect of Studentisation and bandwidth selection.
\end{remark}

\begin{remark}
     Although, all expansions in our paper are for the KDE centralised at $\mathbb{E}\hat{f}_{h_0}(x)$ as \cite{HK01}, centring at $f(x)$, as \cite{Hall92} and \cite{CCF18}, is more desirable from an empirical point of view. Additionally, one of the final goals of the theoretical analysis for KDE with data-driven bandwidth is to simultaneously clarify the effect of bandwidth selection, Studentisation and debias. However, \cite{HK01} and this paper retain some value in the sense that they extract the pure effect of bandwidth selection and the simultaneous effect of bandwidth selection and Studentisation. 
\end{remark}

\section{Simulation Study} \label{section:simulation study}

\subsection{Simulation Settings and Confidence Interval estimation}

In order to examine the higher order improvements by the Edgeworth expansions, we compare the coverage accuracies of the normal approximation, the Cornish-Fisher expansion with optimal bandwidth (\cite{Hall91}), and the Cornish-Fisher expansions with plug-in bandwidth (Theorem \ref{T:Main Theorem} and \ref{T:Edgeworth Expansion Including Pilot Bandwidth}). Following  \cite{MW92}, the underlying distributions are chosen to be a standard normal distribution $N(0,1)$ and a skewed unimodal density constructed as a mixture of $N(0,1)$, $N(1/2,(2/3)^2)$ and $N(13/12,(5/9)^2)$ in the proportions of $1:1:3$. We use the following kernel functions as \cite{HK01}:
\begin{align}
    K(u) = \frac{1}{\sqrt{2\pi}}e^{\frac{-u^2}{2}}, \quad H(u) = \frac{1}{8\sqrt{2\pi}}(u^4-10u^2+15)e^{\frac{-u^2}{2}}, \nonumber
\end{align}
namely $L=2$ and $L_p=6$.

 Let $z_{\alpha}, w_{\alpha}$, $w_{\alpha}^{PI}$, and $w_{\alpha}^{pilot}$ be the $100\alpha\%$-quantile point of normal distribution, Cornish-Fisher expansion of KDE with optimal bandwidth (\cite{Hall91}), Cornish-Fisher expansion of the KDE with plug-in bandwidth and Cornish-Fisher expansion (Theorem \ref{T:Main Theorem}) of the KDE with plug-in and its accompanying pilot bandwidth, respectively (Theorem \ref{T:Edgeworth Expansion Including Pilot Bandwidth}). In this experiment, we set $\alpha=0.05$. We construct the following confidence intervals and count the number of intervals that include $\mathbb{E}\hat{f}_{h_0}(x)$ out of $2000$ iterations. We divide it by $2000$ to compute the empirical coverage probability, and evaluate the performance of each approximation by its closeness to the nominal coverage probability of $0.9500$:
\begin{align}
    & I_{N} = \left[\hat{f}_{\hat{h}_{PI}}-\frac{z_{\alpha/2}\mu_{20}^{1/2}}{\sqrt{n\hat{h}_{PI}}}, \hat{f}_{\hat{h}_{PI}}-\frac{z_{(1-\alpha/2)}\mu_{20}^{1/2}}{\sqrt{n\hat{h}_{PI}}}\right], \quad
    I_{H} = \left[\hat{f}_{\hat{h}_{PI}}-\frac{w_{\alpha/2}\mu_{20}^{1/2}}{\sqrt{n\hat{h}_{PI}}}, \hat{f}_{\hat{h}_{PI}}-\frac{w_{(1-\alpha/2)}\mu_{20}^{1/2}}{\sqrt{n\hat{h}_{PI}}}\right], \nonumber\\
    & I_{PI} = \left[\hat{f}_{\hat{h}_{PI}}-\frac{w_{\alpha/2}^{PI}\mu_{20}^{1/2}}{\sqrt{n\hat{h}_{PI}}},  \hat{f}_{\hat{h}_{PI}}-\frac{w^{PI}_{(1-\alpha/2)}\mu_{20}^{1/2}}{\sqrt{n\hat{h}_{PI}}}\right], \quad 
    I_{pilot}^L = \left[\hat{f}_{\hat{h}_{PI}}-\frac{w_{\alpha/2}^{pilot}\mu_{20}^{1/2}}{\sqrt{n\hat{h}_{PI}}}, \hat{f}_{\hat{h}_{PI}}-\frac{w^{pilot}_{(1-\alpha/2)}\mu_{20}^{1/2}}{\sqrt{n\hat{h}_{PI}}}\right] \nonumber.
\end{align}

 The experiment is conducted with MSE-optimal pilot bandwidth for sample sizes $n=50, 100, 400$, and $1000$. The MSE-optimal pilot bandwidth is defined as follows (See Lemma 3.1 of \cite{HM87} for the proof):
\begin{align}
    & b_0 \equiv \left(\frac{(4L+1)\left\{\int f(x)^2dx\right\}^2\int (H^{(2L)}*H)(u)^2du}{L_p(L_p!)^{-2}\left\{\int u^{L_p}H(u)du\right\}^2\left\{\int f^{(L)}(x)f^{(L+L_p)}(x)dx\right\}^2}\right)^{1/(4L+2L_p+1)}n^{-2/(4L+2L_p+1)} \nonumber,\\
    & b_0^{convo} \equiv \left(\frac{(4L+1)\left\{\int f(x)^2dx\right\}^2\int H^{(2L)}(u)^2du}{L_p(L_p!)^{-2}\left\{\int u^{L_p}(H*H)(u)du\right\}^2\left\{\int f^{(L)}(x)f^{(L+L_p)}(x)dx\right\}^2}\right)^{1/(4L+2L_p+1)}n^{-2/(4L+2L_p+1)} \nonumber.
\end{align}
where $*$ denotes the convolution.

\subsection{Simulation Results with $\hat{I}_L$}
Tables \ref{tab:SNx0}-\ref{tab:SNx2} in Section \ref{subsection:simulation N(0,1)} report the nominal coverage probabilities for five evaluation points $x=0, 0.5, 1.0, 1.5$ and $2.0$ in the case of $N(0,1)$ observations. In each table, results for the sample sizes of $n=\{50,100,400,1000\}$ are shown when we approximate the distribution of $S_{PI}(x)$ by $N(0,1)$, \cite{Hall91}'s Edgeworth expansion, Theorem \ref{T:Main Theorem}, and \ref{T:Edgeworth Expansion Including Pilot Bandwidth}. In each row, ** and * indicate the closest and second-closest value to the nominal coverage probability of $0.9500$. Similarly, Tables \ref{tab:SUxm2}-\ref{tab:SUx2} in Section \ref{subsection:skewed unimodal} present the results for $x=-2, -1.5, -1.0, -0.5, 0, 0.5, 1.0, 1.5$ and $2.0$ with skewed unimodal normal mixture. 

We also conducted a simulation when we adopted $\hat {I}_L^{convo}$ to estimate $I_L$ but we suppressed the results because they are qualitatively similar. We provide them in the supplemental material (\ref{section:supp-simulation}).

\subsubsection{Standard Normal} \label{subsection:simulation N(0,1)}
We adopt $\# 1: N(0,1)$ in \cite{MW92}. For sample size $n=(50, 100, 400, 1000)$, MSE optimal pilot bandwidths are $b_0=(0.8448, 0.7908, 0.6930, 0.6351)$. We evaluate the accuracy at the point of $x=0, 0.5, 1, 1.5$, and $x=2$ in Tables \ref{tab:SNx0}-\ref{tab:SNx2} respectively.

From Tables \ref{tab:SNx0}-\ref{tab:SNx2}, we observe that approximation by Theorems \ref{T:Main Theorem} and \ref{T:Edgeworth Expansion Including Pilot Bandwidth} outperform the N(0,1) or Hall's approximations with some exceptions with mainly small $n$ (see Tables \ref{tab:SNx1} and \ref{tab:SNx15}). We also see, in Table \ref{tab:SNx1}, that N(0,1) and Hall's expansion provide the closest coverage rate to 0.9500, but the differences in the coverage probability with Theorems \ref{T:Main Theorem} and \ref{T:Edgeworth Expansion Including Pilot Bandwidth} are only marginal. It is not clear which performs better Theorem \ref{T:Main Theorem} or \ref{T:Edgeworth Expansion Including Pilot Bandwidth} depending on the evaluation points and sample size. We conclude that Edgeworth expansions obtained mostly improve the confidence interval estimation in this case.

We point out that coverage ratios in Tables \ref{tab:SNx1} and \ref{tab:SNx15} are satisfactory in the level, that is, are close to the nominal probability of 0.9500, while Tables \ref{tab:SNx0}, \ref{tab:SNx05}, and \ref{tab:SNx2} provide dismal performance independent of the approximation methods. We further find, in Tables \ref{tab:SNx0} , \ref{tab:SNx15}, and \ref{tab:SNx2}, that increase in sample size does not improve the confidence interval estimation. We discuss this issue at the end of this section. 

Naturally, in any situation, the average length of intervals gets shorter as the sample size increases. Moreover, in most cases, the confidence intervals created by Theorem \ref{T:Main Theorem} and \ref{T:Edgeworth Expansion Including Pilot Bandwidth} are longer than those by $N(0,1)$ and \cite{Hall91}'s expansion. Except for the case of $n=50$ in Table \ref{tab:SNx15}, the coverage probabilities by the approximation of $N(0,1)$ are much less than $0.9500$, and Theorem \ref{T:Main Theorem} and \ref{T:Edgeworth Expansion Including Pilot Bandwidth} correct the approximation error by providing relatively long confidence intervals.

\begin{table}[H]
\small
\centering
    \caption{$x=0, b=\text{MSE optimal},~~$ scaled second derivative$=0.4122$}
    {\begin{tabular}{c||llllllll}
      \hline & \multicolumn{2}{c|}{$n=50$}   &  \multicolumn{2}{c|}{$n=100$} & \multicolumn{2}{c|}{$n=400$} & \multicolumn{2}{c}{$n=1000$}\\
     \hline & \multicolumn{1}{|c|}{CP} & \multicolumn{1}{|c|}{Ave.Length} & \multicolumn{1}{|c|}{CP} & \multicolumn{1}{|c|}{Ave.Length} & \multicolumn{1}{|c|}{CP} & \multicolumn{1}{|c|}{Ave.Length} & \multicolumn{1}{|c|}{CP} & \multicolumn{1}{|c}{Ave.Length} \\
     \hline\hline $N(0,1)$ & $0.5200$ & $0.1264$ & $0.5065$ & 0.1007 & 0.5660 & 0.0670 & 0.5730 & 0.0492 \\
     Hall (1991) & 0.5150 & 0.1265 & 0.5005 & 0.1008 & 0.5615 & 0.0670 & 0.5660 & 0.0492 \\
     Theorem 3.1 & 0.7070$^{*}$ & 0.1762 & 0.6005$^{*}$ & 0.1267 & 0.6270$^{**}$ & 0.0744 & 0.6110$^{**}$ & 0.0523 \\
     Theorem 3.5 & 0.7180$^{**}$ & 0.2115 & 0.6020$^{**}$ & 0.1438 & 0.6240$^{*}$ & 0.0785 & 0.6085$^{*}$ & 0.0538
    \end{tabular}}
    \label{tab:SNx0}
\end{table}

\begin{table}[H]
\small
\centering
    \caption{$x=0.5, b=\text{MSE optimal},~~$ scaled second derivative$=0.2728$}
    {\begin{tabular}{c||llllllll}
      \hline & \multicolumn{2}{c|}{$n=50$}   &  \multicolumn{2}{c|}{$n=100$} & \multicolumn{2}{c|}{$n=400$} & \multicolumn{2}{c}{$n=1000$}\\
     \hline & \multicolumn{1}{|c|}{CP} & \multicolumn{1}{|c|}{Ave.Length} & \multicolumn{1}{|c|}{CP} & \multicolumn{1}{|c|}{Ave.Length} & \multicolumn{1}{|c|}{CP} & \multicolumn{1}{|c|}{Ave.Length} & \multicolumn{1}{|c|}{CP} & \multicolumn{1}{|c}{Ave.Length} \\
     \hline\hline $N(0,1)$ & 0.6780 & 0.1266 & 0.6070 & 0.0999 & 0.7500 & 0.0654 & 0.7605 & 0.0477 \\
     Hall (1991) & 0.6700 & 0.1267 & 0.6010 & 0.1000 & 0.7460 & 0.0654 & 0.7500 & 0.0477 \\
     Theorem 3.1 & 0.7195$^{*}$ & 0.1434 & 0.6310$^{*}$
     & 0.1081 & 0.7585$^{*}$ & 0.0675 & 0.7665$^{*}$ & 0.0485 \\
     Theorem 3.5 & 0.7685$^{**}$ & 0.1729 & 0.6550$^{**}$ & 0.1222 & 0.7690$^{**}$ & 0.0708 & 0.7700$^{**}$ & 0.0497
    \end{tabular}}
    \label{tab:SNx05}
\end{table}

\begin{table}[H]
\small
\centering
    \caption{$x=1, b=\text{MSE optimal},~~$ scaled second derivative$=0$}
    {\begin{tabular}{c||llllllll}
      \hline & \multicolumn{2}{c|}{$n=50$}   &  \multicolumn{2}{c|}{$n=100$} & \multicolumn{2}{c|}{$n=400$} & \multicolumn{2}{c}{$n=1000$}\\
     \hline & \multicolumn{1}{|c|}{CP} & \multicolumn{1}{|c|}{Ave.Length} & \multicolumn{1}{|c|}{CP} & \multicolumn{1}{|c|}{Ave.Length} & \multicolumn{1}{|c|}{CP} & \multicolumn{1}{|c|}{Ave.Length} & \multicolumn{1}{|c|}{CP} & \multicolumn{1}{|c}{Ave.Length} \\
     \hline\hline $N(0,1)$ & 0.9265$^{*}$ & 0.1209 & 0.9130$^{*}$ & 0.0934 & 0.9655 & 0.0592 & 0.9650$^{**}$ & 0.0424 \\
     Hall (1991) & 0.9200 & 0.1209 & 0.9060 & 0.0934 & 0.9640$^{**}$ & 0.0592 & 0.9650$^{**}$ & 0.0424 \\
     Theorem 3.1 & 0.9240 & 0.1209 & 0.9065 & 0.0934 & 0.9640$^{**}$ & 0.0592 & 0.9655 & 0.0424 \\
     Theorem 3.5 & 0.9540$^{**}$ & 0.1375 & 0.9340$^{**}$ & 0.1011 & 0.9680 & 0.0609 & 0.9670 & 0.0430
    \end{tabular}}
    \label{tab:SNx1}
\end{table}

\begin{table}[H]
\small
\centering
    \caption{$x=1.5, b=\text{MSE optimal},~~$ scaled second derivative$=0.1673$}
    {\begin{tabular}{c||llllllll}
      \hline & \multicolumn{2}{c|}{$n=50$}   &  \multicolumn{2}{c|}{$n=100$} & \multicolumn{2}{c|}{$n=400$} & \multicolumn{2}{c}{$n=1000$}\\
     \hline & \multicolumn{1}{|c|}{CP} & \multicolumn{1}{|c|}{Ave.Length} & \multicolumn{1}{|c|}{CP} & \multicolumn{1}{|c|}{Ave.Length} & \multicolumn{1}{|c|}{CP} & \multicolumn{1}{|c|}{Ave.Length} & \multicolumn{1}{|c|}{CP} & \multicolumn{1}{|c}{Ave.Length} \\
     \hline\hline $N(0,1)$ & 0.9495$^{**}$ & 0.1022 & 0.9140 & 0.0773 & 0.8645 & 0.0476 & 0.8330 & 0.0333 \\
     Hall (1991) & 0.9550$^{*}$ & 0.1021 & 0.9305 & 0.0772 & 0.8755 & 0.0472 & 0.8405 & 0.0333 \\
     Theorem 3.1 & 0.9760 & 0.1131 & 0.9520$^{**}$ & 0.0830 & 0.8880$^{*}$ & 0.0489 & 0.8520$^{**}$ & 0.0339 \\
     Theorem 3.5 & 0.9790 & 0.1188 & 0.9530$^{*}$ & 0.0856 & 0.8885$^{**}$ & 0.0495 & 0.8520$^{**}$ & 0.0341
    \end{tabular}}
    \label{tab:SNx15}
\end{table}

\begin{table}[H]
\small
\centering
    \caption{$x=2, b=\text{MSE optimal},~~$ scaled second derivative$=0.1673$}
    {\begin{tabular}{c||llllllll}
      \hline & \multicolumn{2}{c|}{$n=50$}   &  \multicolumn{2}{c|}{$n=100$} & \multicolumn{2}{c|}{$n=400$} & \multicolumn{2}{c}{$n=1000$}\\
     \hline & \multicolumn{1}{|c|}{CP} & \multicolumn{1}{|c|}{Ave.Length} & \multicolumn{1}{|c|}{CP} & \multicolumn{1}{|c|}{Ave.Length} & \multicolumn{1}{|c|}{CP} & \multicolumn{1}{|c|}{Ave.Length} & \multicolumn{1}{|c|}{CP} & \multicolumn{1}{|c}{Ave.Length} \\
     \hline\hline $N(0,1)$ & 0.6865 & 0.0746 & 0.5995 & 0.0552 & 0.6410 & 0.0327 & 0.6240 & 0.0227 \\
     Hall (1991) & 0.7295 & 0.0742 & 0.6350 & 0.0550 & 0.6655 & 0.0327 & 0.6445 & 0.0227 \\
     Theorem 3.1 & 0.7865$^{**}$ & 0.0823 & 0.6625$^{**}$ & 0.0588 & 0.6850$^{**}$ & 0.0335 & 0.6570$^{**}$ & 0.0230 \\
     Theorem 3.5 & 0.7780$^{*}$ & 0.0833 & 0.6580$^{*}$ & 0.0593 & 0.6845$^{*}$ & 0.0336 & 0.6550$^{*}$ & 0.0230
    \end{tabular}}
    \label{tab:SNx2}
\end{table}

\subsubsection{Skewed Unimodal} \label{subsection:skewed unimodal}
We adopt $\# 2: \frac{1}{5}N(0,1)+\frac{1}{5}N(\frac{1}{2},(\frac{2}{3})^2)+\frac{3}{5}N(\frac{13}{12},(\frac{5}{9})^2)$ in \cite{MW92}. \\For sample size $n=(50, 100, 400, 1000)$, MSE optimal pilot bandwidths are $b_0=(0.5227, 0.4893, 0.4287, 0.3929)$. We evaluate the accuracy at the point of $x=-2, -1.5, -1, -0.5, 0, 0.5, 1, 1.5$, and $x=2$.

Similar to the case of standard normal observations in the previous section, we find that Theorems \ref{T:Main Theorem} and \ref{T:Edgeworth Expansion Including Pilot Bandwidth} outperform $N(0,1)$ and Hall's approximations in general. 

We observe in Tables \ref{tab:SUx1} and \ref{tab:SUx2} that the general coverage probability level significantly differs from the nominal value of whichever approximation we adopt and further the results look to contradict the asymptotic theory. We discuss this in the next subsection. 

As with the case of standard normal, as the sample size increases, the intervals also get shorter. Moreover, Theorem \ref{T:Main Theorem} and \ref{T:Edgeworth Expansion Including Pilot Bandwidth} provide the longer confidence intervals and thereby achieve  coverage probabilities closer to 0.95 than $N(0,1)$ and \cite{Hall91}'s approximation, at the point with poor coverage. 

\begin{table}[H]
\small
\centering
    \caption{$x=-2, b=\text{MSE optimal},~~$ scaled secand derivative$=0.0173$}
    {\begin{tabular}{c||llllllll}
      \hline & \multicolumn{2}{c|}{$n=50$}   &  \multicolumn{2}{c|}{$n=100$} & \multicolumn{2}{c|}{$n=400$} & \multicolumn{2}{c}{$n=1000$}\\
     \hline & \multicolumn{1}{|c|}{CP} & \multicolumn{1}{|c|}{Ave.Length} & \multicolumn{1}{|c|}{CP} & \multicolumn{1}{|c|}{Ave.Length} & \multicolumn{1}{|c|}{CP} & \multicolumn{1}{|c|}{Ave.Length} & \multicolumn{1}{|c|}{CP} & \multicolumn{1}{|c}{Ave.Length} \\
     \hline\hline $N(0,1)$ & 0.8675 & 0.0424 & 0.8150 & 0.0299 & 0.9190 & 0.0179 & 0.9270 & 0.0125 \\
     Hall (1991) & 0.9155 & 0.0396 & 0.8615$^{*}$ & 0.0288 & 0.9395$^{**}$ & 0.0176 & 0.9440$^{**}$ & 0.0124 \\
     Theorem 3.1 & 0.9160$^{*}$ & 0.0397 & 0.8615$^{*}$ & 0.0288 & 0.9395$^{**}$ & 0.0176 & 0.9440$^{**}$ & 0.0124 \\
     Theorem 3.5 & 0.9165$^{**}$ & 0.0398 & 0.8630$^{**}$ & 0.0288 & 0.9395$^{**}$ & 0.0177 & 0.9440$^{**}$ & 0.0124
    \end{tabular}}
    \label{tab:SUxm2}
\end{table}

\begin{table}[H]
\small
\centering
    \caption{$x=-1.5, b=\text{MSE optimal},~~$ scaled second derivative$=0.0278$}
    {\begin{tabular}{c||llllllll}
      \hline & \multicolumn{2}{c|}{$n=50$}   &  \multicolumn{2}{c|}{$n=100$} & \multicolumn{2}{c|}{$n=400$} & \multicolumn{2}{c}{$n=1000$}\\
     \hline & \multicolumn{1}{|c|}{CP} & \multicolumn{1}{|c|}{Ave.Length} & \multicolumn{1}{|c|}{CP} & \multicolumn{1}{|c|}{Ave.Length} & \multicolumn{1}{|c|}{CP} & \multicolumn{1}{|c|}{Ave.Length} & \multicolumn{1}{|c|}{CP} & \multicolumn{1}{|c}{Ave.Length}\\
     \hline\hline $N(0,1)$ & 0.8635 & 0.0654 & 0.8070 & 0.0464 & 0.9210 & 0.0279 & 0.9430 & 0.0196 \\
     Hall (1991) & 0.8925$^{*}$ & 0.0637 & 0.8335 & 0.0457 & 0.9355$^{**}$ & 0.0277 & 0.9525$^{**}$ & 0.0195 \\
     Theorem 3.1 & 0.8925$^{*}$ & 0.0639 & 0.8340$^{*}$ & 0.0457 & 0.9355$^{**}$ & 0.0278 & 0.9525$^{**}$ & 0.0196 \\
     Theorem 3.5 & 0.8955$^{**}$ & 0.0643 & 0.8350$^{**}$ & 0.0460 & 0.9355$^{**}$ & 0.0278 & 0.9525$^{**}$ & 0.0196
    \end{tabular}}
    \label{tab:SUxm15}
\end{table}

\begin{table}[H]
\small
\centering
    \caption{$x=-1, b=\text{MSE optimal},~~$ scaled second derivative$=0.0503$}
    {\begin{tabular}{c||llllllll}
      \hline & \multicolumn{2}{c|}{$n=50$}   &  \multicolumn{2}{c|}{$n=100$} & \multicolumn{2}{c|}{$n=400$} & \multicolumn{2}{c}{$n=1000$}\\
     \hline & \multicolumn{1}{|c|}{CP} & \multicolumn{1}{|c|}{Ave.Length} & \multicolumn{1}{|c|}{CP} & \multicolumn{1}{|c|}{Ave.Length} & \multicolumn{1}{|c|}{CP} & \multicolumn{1}{|c|}{Ave.Length} & \multicolumn{1}{|c|}{CP} & \multicolumn{1}{|c}{Ave.Length}\\
     \hline\hline $N(0,1)$ & 0.8490 & 0.0933 & 0.7550 & 0.0664 & 0.8870 & 0.0402 & 0.9375 & 0.0283 \\
     Hall (1991) & 0.8780 & 0.0923 & 0.7870 & 0.0660 & 0.9070 & 0.0401 & 0.9450 & 0.0283 \\
     Theorem 3.1 & 0.8795$^{*}$ & 0.0927 & 0.7880$^{*}$ & 0.0662 & 0.9075$^{*}$ & 0.0402 & 0.9455$^{**}$ & 0.0283 \\
     Theorem 3.5 & 0.8855$^{**}$ & 0.0941 & 0.7910$^{**}$ & 0.0669 & 0.9085$^{**}$ & 0.0403 & 0.9455$^{**}$ & 0.0284
    \end{tabular}}
    \label{tab:SUxm1}
\end{table}

\begin{table}[H]
\small
\centering
    \caption{$x=-0.5, b=\text{MSE optimal},~~$ scaled second derivative$=0.1112$}
    {\begin{tabular}{c||llllllll}
      \hline & \multicolumn{2}{c|}{$n=50$}   &  \multicolumn{2}{c|}{$n=100$} & \multicolumn{2}{c|}{$n=400$} & \multicolumn{2}{c}{$n=1000$}\\
     \hline & \multicolumn{1}{|c|}{CP} & \multicolumn{1}{|c|}{Ave.Length} & \multicolumn{1}{|c|}{CP} & \multicolumn{1}{|c|}{Ave.Length} & \multicolumn{1}{|c|}{CP} & \multicolumn{1}{|c|}{Ave.Length} & \multicolumn{1}{|c|}{CP} & \multicolumn{1}{|c}{Ave.Length}\\
     \hline\hline $N(0,1)$ & 0.8685 & 0.1265 & 0.7410 & 0.0907 & 0.8540 & 0.0555 & 0.8950 & 0.0393 \\
     Hall (1991) & 0.8945 & 0.1260 & 0.7675 & 0.0905 & 0.8720 & 0.0554 & 0.9030$^{*}$ & 0.0393 \\
     Theorem 3.1 & 0.9055$^{*}$ & 0.1286 & 0.7765$^{*}$ & 0.0916 & 0.8725$^{*}$ & 0.0556 & 0.9030$^{*}$ & 0.0393 \\
     Theorem 3.5 & 0.9120$^{**}$ & 0.1328 & 0.7830$^{**}$ & 0.0936 & 0.8735$^{**}$ & 0.0562 & 0.9055$^{**}$ & 0.0395
    \end{tabular}}
    \label{tab:SUxm05}
\end{table}

\begin{table}[H]
\small
\centering
    \caption{$x=0, b=\text{MSE optimal},~~$ scaled second derivative$=0.2029$}
    {\begin{tabular}{c||llllllll}
      \hline & \multicolumn{2}{c|}{$n=50$}   &  \multicolumn{2}{c|}{$n=100$} & \multicolumn{2}{c|}{$n=400$} & \multicolumn{2}{c}{$n=1000$}\\
     \hline & \multicolumn{1}{|c|}{CP} & \multicolumn{1}{|c|}{Ave.Length} & \multicolumn{1}{|c|}{CP} & \multicolumn{1}{|c|}{Ave.Length} & \multicolumn{1}{|c|}{CP} & \multicolumn{1}{|c|}{Ave.Length} & \multicolumn{1}{|c|}{CP} & \multicolumn{1}{|c}{Ave.Length}\\
     \hline\hline $N(0,1)$ & 0.9630$^{**}$ & 0.1626 & 0.9510$^{**}$ & 0.1187 & 0.9145 & 0.0744 & 0.9205 & 0.0533 \\
     Hall (1991) & 0.9670$^{*}$ & 0.1625 & 0.9590$^{*}$ & 0.1186 & 0.9225 & 0.0744 & 0.9225 & 0.0533 \\
     Theorem 3.1 & 0.9795 & 0.1777 & 0.9710 & 0.1259 & 0.9350$^{*}$ & 0.0763 & 0.9290$^{*}$ & 0.0540 \\
     Theorem 3.5 & 0.9900 & 0.1929 & 0.9780 & 0.1328 & 0.9375$^{**}$ & 0.0780 & 0.9305$^{**}$ & 0.0547    
     \end{tabular}}
    \label{tab:SUx0}
\end{table}

\begin{table}[H]
\small
\centering
    \caption{$x=0.5, b=\text{MSE optimal},~~$ scaled second derivative$=0.1170$}
    {\begin{tabular}{c||llllllll}
      \hline & \multicolumn{2}{c|}{$n=50$}   &  \multicolumn{2}{c|}{$n=100$} & \multicolumn{2}{c|}{$n=400$} & \multicolumn{2}{c}{$n=1000$}\\
     \hline & \multicolumn{1}{|c|}{CP} & \multicolumn{1}{|c|}{Ave.Length} & \multicolumn{1}{|c|}{CP} & \multicolumn{1}{|c|}{Ave.Length} & \multicolumn{1}{|c|}{CP} & \multicolumn{1}{|c|}{Ave.Length} & \multicolumn{1}{|c|}{CP} & \multicolumn{1}{|c}{Ave.Length}\\
     \hline\hline $N(0,1)$ & 0.8595$^{*}$ & 0.1867 & 0.7375$^{*}$ & 0.1399 & 0.9185$^{*}$ & 0.0915 & 0.9545 & 0.0669 \\
     Hall (1991) & 0.8520 & 0.1867 & 0.7290 & 0.1400 & 0.9105 & 0.0915 & 0.9500$^{**}$ & 0.0670 \\
     Theorem 3.1 & 0.8510 & 0.1841 & 0.7260 & 0.1380 & 0.9065 & 0.0907 & 0.9475$^{*}$ & 0.0665 \\
     Theorem 3.5 & 0.9120 & 0.2312 & 0.7905 & 0.1597 & 0.9265 & 0.0960 & 0.9560 & 0.0685
     \end{tabular}}
    \label{tab:SUx05}
\end{table}

\begin{table}[H]
\small
\centering
    \caption{$x=1, b=\text{MSE optimal},~~$ scaled second derivative$=0.7019$}
    {\begin{tabular}{c||llllllll}
      \hline & \multicolumn{2}{c|}{$n=50$}   &  \multicolumn{2}{c|}{$n=100$} & \multicolumn{2}{c|}{$n=400$} & \multicolumn{2}{c}{$n=1000$}\\
     \hline & \multicolumn{1}{|c|}{CP} & \multicolumn{1}{|c|}{Ave.Length} & \multicolumn{1}{|c|}{CP} & \multicolumn{1}{|c|}{Ave.Length} & \multicolumn{1}{|c|}{CP} & \multicolumn{1}{|c|}{Ave.Length} & \multicolumn{1}{|c|}{CP} & \multicolumn{1}{|c}{Ave.Length}\\
     \hline\hline $N(0,1)$ & 0.6345 & 0.1923 & 0.3620 & 0.1457 & 0.4835 & 0.0973 & 0.5265 & 0.0720 \\
     Hall (1991) & 0.6240 & 0.1925 & 0.3545 & 0.1458 & 0.4810 & 0.0973 & 0.5225 & 0.0720 \\
     Theorem 3.1 & 0.8370$^{*}$ & 0.2930 & 0.5015$^{*}$ & 0.1977 & 0.5630$^{**}$ & 0.1130 & 0.5845$^{**}$ & 0.0788 \\
     Theorem 3.5 & 0.8405$^{**}$ & 0.3556 & 0.5035$^{**}$ & 0.2279 & 0.5595$^{*}$ & 0.1209 & 0.5785$^{*}$ & 0.0819   
     \end{tabular}}
    \label{tab:SUx1}
\end{table}

\begin{table}[H]
\small
\centering
    \caption{$x=1.5, b=\text{MSE optimal},~~$ scaled second derivative$=0.1559$}
    {\begin{tabular}{c||llllllll}
      \hline & \multicolumn{2}{c|}{$n=50$}   &  \multicolumn{2}{c|}{$n=100$} & \multicolumn{2}{c|}{$n=400$} & \multicolumn{2}{c}{$n=1000$}\\
     \hline & \multicolumn{1}{|c|}{CP} & \multicolumn{1}{|c|}{Ave.Length} & \multicolumn{1}{|c|}{CP} & \multicolumn{1}{|c|}{Ave.Length} & \multicolumn{1}{|c|}{CP} & \multicolumn{1}{|c|}{Ave.Length} & \multicolumn{1}{|c|}{CP} & \multicolumn{1}{|c}{Ave.Length}\\
     \hline\hline $N(0,1)$ & 0.8530 & 0.1849 & 0.7390 & 0.1373 & 0.8480 & 0.0886 & 0.8895 & 0.0644 \\
     Hall (1991) & 0.8485 & 0.1850 & 0.7290 & 0.1373 & 0.8430 & 0.0886 & 0.8865 & 0.0644 \\
     Theorem 3.1 & 0.8525$^{*}$ & 0.1893 & 0.7305$^{*}$ & 0.1387 & 0.8430$^{*}$ & 0.0886 & 0.8865$^{*}$ & 0.0643 \\
     Theorem 3.5 & 0.8985$^{**}$ & 0.2259 & 0.7765$^{**}$ & 0.1556 & 0.8625$^{**}$ & 0.0928 & 0.8905$^{**}$ & 0.0659    
     \end{tabular}}
    \label{tab:SUx15}
\end{table}

\begin{table}[H]
\small
\centering
    \caption{$x=2, b=\text{MSE optimal},~~$ scaled second derivative$=0.3591$}
    {\begin{tabular}{c||llllllll}
      \hline & \multicolumn{2}{c|}{$n=50$}   &  \multicolumn{2}{c|}{$n=100$} & \multicolumn{2}{c|}{$n=400$} & \multicolumn{2}{c}{$n=1000$}\\
     \hline & \multicolumn{1}{|c|}{CP} & \multicolumn{1}{|c|}{Ave.Length} & \multicolumn{1}{|c|}{CP} & \multicolumn{1}{|c|}{Ave.Length} & \multicolumn{1}{|c|}{CP} & \multicolumn{1}{|c|}{Ave.Length} & \multicolumn{1}{|c|}{CP} & \multicolumn{1}{|c}{Ave.Length}\\
     \hline\hline $N(0,1)$ & 0.8245 & 0.1349 & 0.6745 & 0.0968 & 0.6275 & 0.0591 & 0.6070 & 0.0418 \\
     Hall (1991) & 0.8450 & 0.1346 & 0.6995 & 0.1346 & 0.6460 & 0.0590 & 0.6215 & 0.0418 \\
     Theorem 3.1 & 0.9040$^{**}$ & 0.1599 & 0.7710$^{**}$ & 0.1098 & 0.6900$^{**}$ & 0.0628 & 0.6450$^{**}$ & 0.0433 \\
     Theorem 3.5 & 0.9000$^{*}$ & 0.1659 & 0.7650$^{*}$ & 0.1124 & 0.6865$^{*}$ & 0.0635 & 0.6435$^{*}$ & 0.0436    
     \end{tabular}}
    \label{tab:SUx2}
\end{table}

\subsection{Difficulties in Confidence Interval Estimation with High Curvature} \label{subsection:difficulties}

\begin{figure}[H]
    \centering
    \includegraphics[scale=0.3]{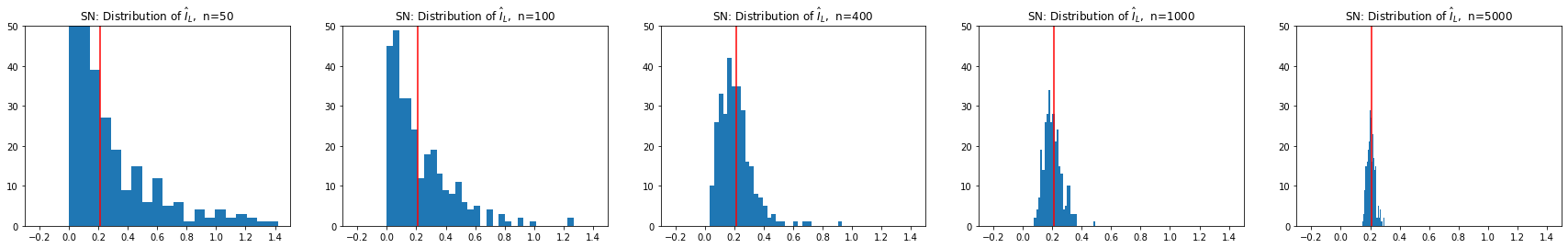}
    \caption{Distributions of $\hat{I}_L$ estimated with $\hat{h}$}
    \label{fig:SN_dist_IL_hat}
\end{figure}

\begin{figure}[H]
    \centering
    \includegraphics[scale=0.3]{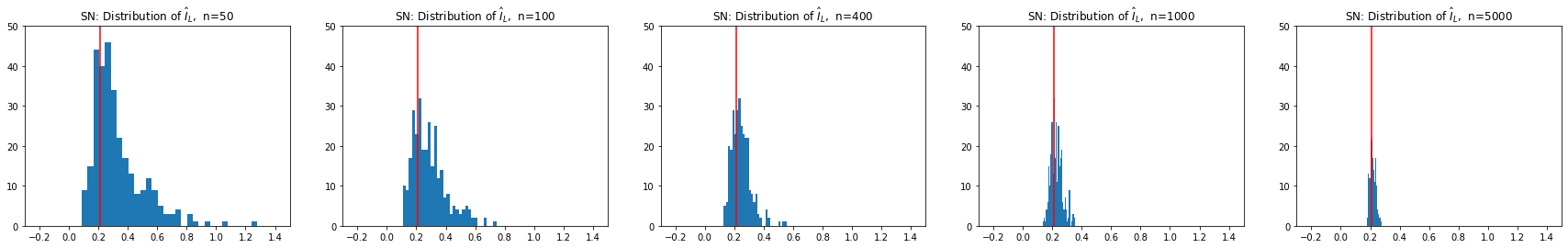}
    \caption{Distributions of $\hat{I}_L$ estimated with $h_{convo}$}
    \label{fig:SN_dist_IL_convo}
\end{figure}

\begin{figure}[H]
    \centering
    \includegraphics[scale=0.3]{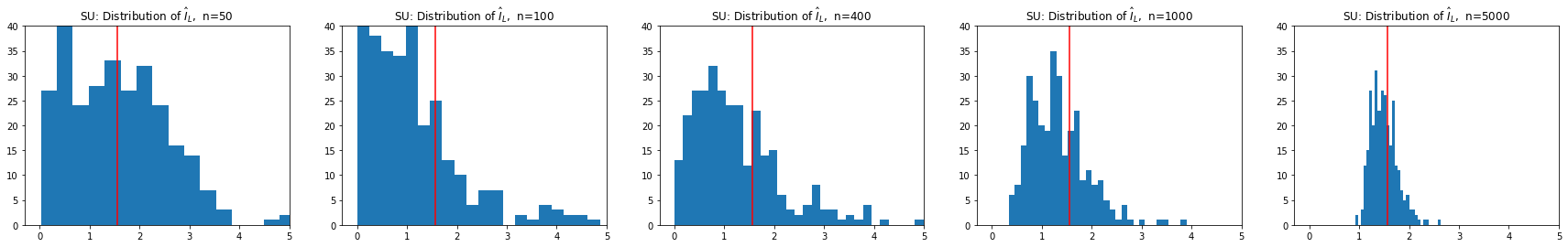}
    \caption{Distributions of $\hat{I}_L$ estimated with $\hat{h}$}
    \label{fig:SU_dist_IL_hat}
\end{figure}

\begin{figure}[H]
    \centering
    \includegraphics[scale=0.3]{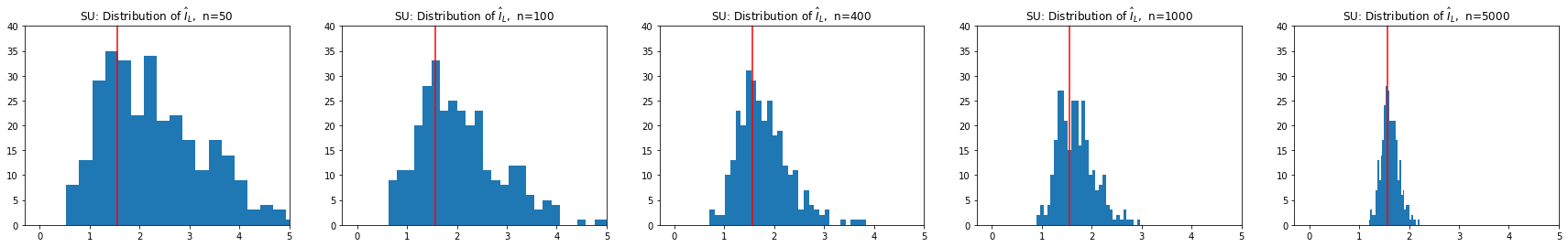}
    \caption{Distributions of $\hat{I}_L$ estimated with $h_{convo}$}
    \label{fig:SU_dist_IL_convo}
\end{figure}

\begin{figure}[H]
    \centering
    \includegraphics[scale=0.3]{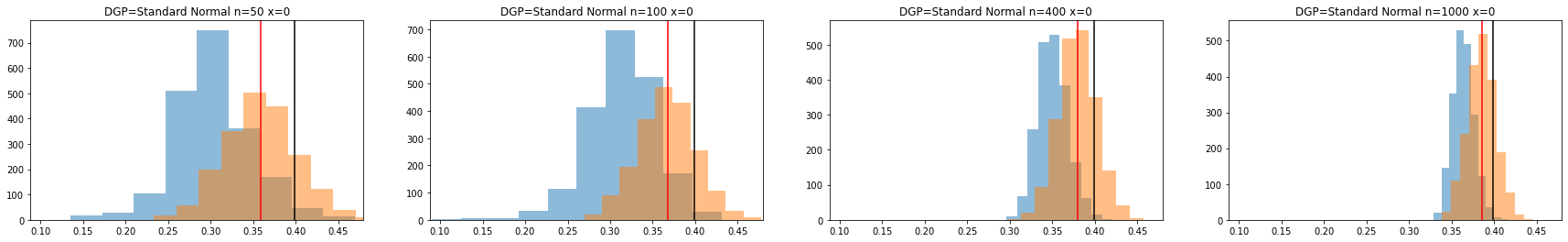}
    \caption{Distributions of $\hat{f}_{\hat{h}}$ and $\hat{f}_{h_0}$}
    \label{fig:SN_dist_KDE_x0_hat}
\end{figure}

\begin{figure}[H]
    \centering
    \includegraphics[scale=0.3]{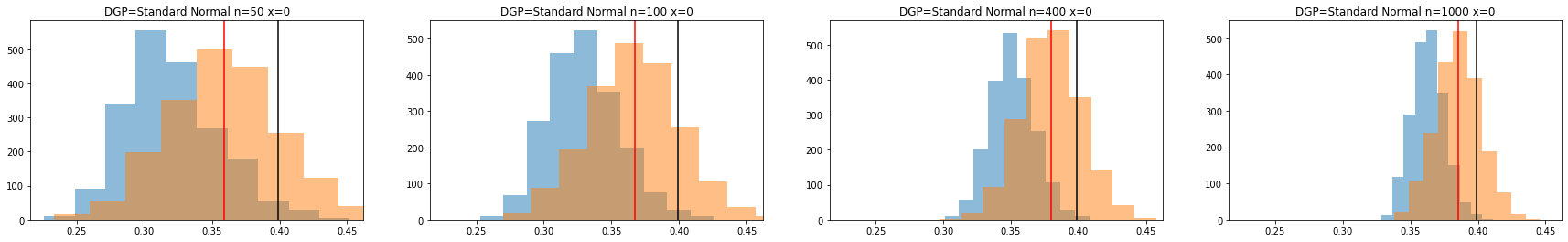}
    \caption{Distributions of $\hat{f}_{h_{convo}}$ and $\hat{f}_{h_0}$}
    \label{fig:SN_dist_KDE_x0_convo}
\end{figure}

\begin{figure}[H]
    \centering
    \includegraphics[scale=0.3]{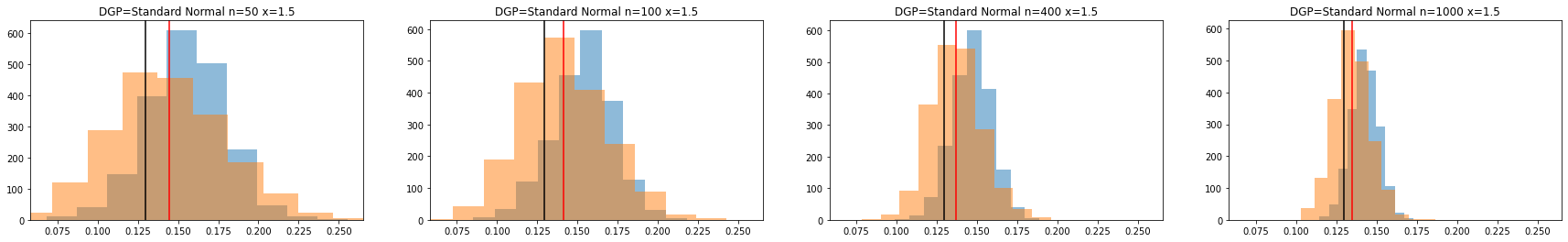}
    \caption{Distributions of $\hat{f}_{\hat{h}}$ and $\hat{f}_{h_0}$}
    \label{fig:SN_dist_KDE_x15_hat}
\end{figure}

\begin{figure}[H]
    \centering
    \includegraphics[scale=0.3]{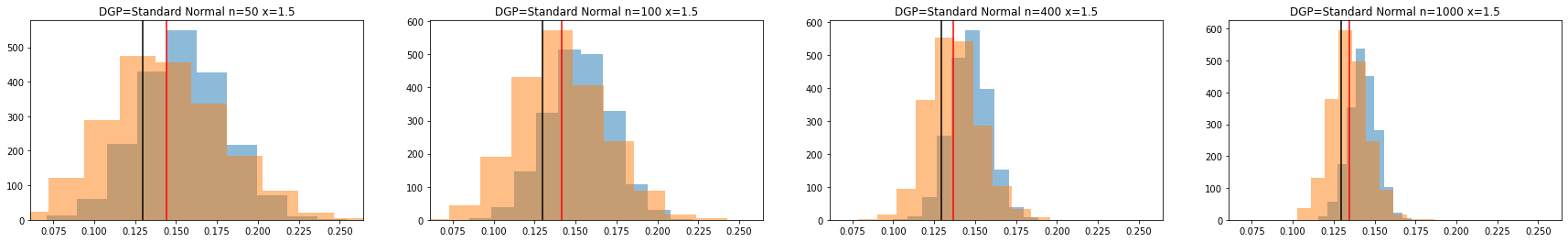}
    \caption{Distributions of $\hat{f}_{h_{convo}}$ and $\hat{f}_{h_0}$}
    \label{fig:SN_dist_KDE_x15_convo}
\end{figure}

\begin{figure}[H]
    \centering
    \includegraphics[scale=0.3]{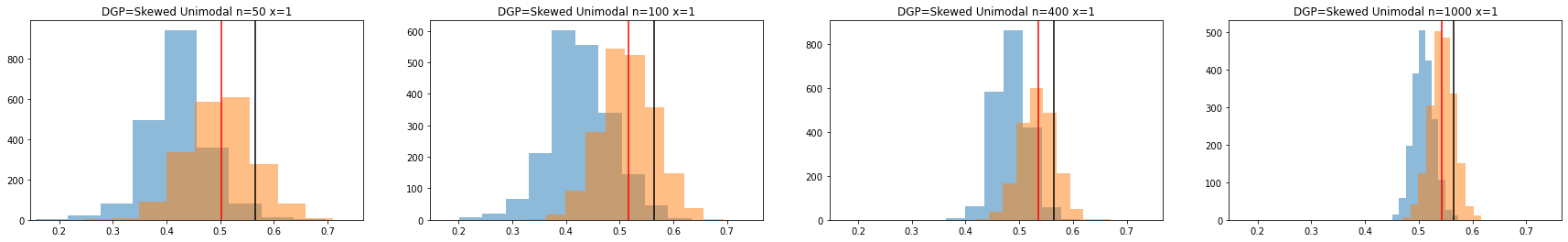}
    \caption{Distributions of $\hat{f}_{\hat{h}}$ and $\hat{f}_{h_0}$}
    \label{fig:SU_dist_KDE_x1_hat}
\end{figure}

\begin{figure}[H]
    \centering
    \includegraphics[scale=0.3]{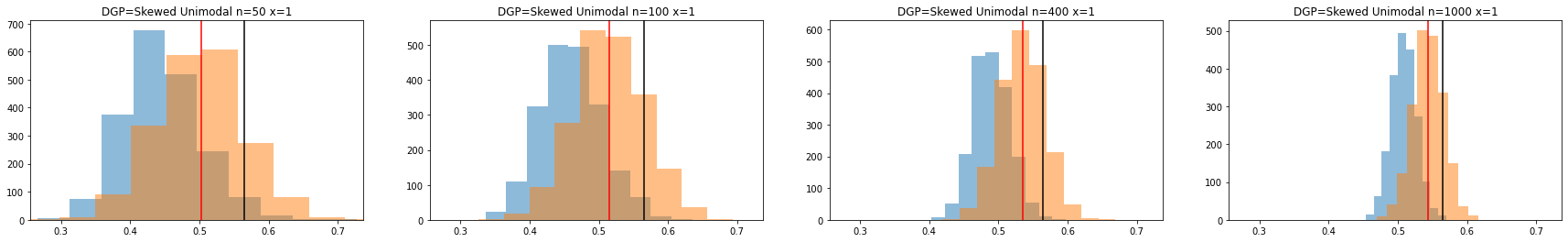}
    \caption{Distributions of $\hat{f}_{h_{convo}}$ and $\hat{f}_{h_0}$}
    \label{fig:SU_dist_KDE_x1_convo}
\end{figure}

We provide the tables of simulation results in the above subsections. Note that at some points, the results contradict the asymptotic theory. Such phenomena seem to occur at points of high curvature of the density function (e.g. Table \ref{tab:SNx0}, \ref{tab:SNx2}, \ref{tab:SUx1}, \ref{tab:SUx2}, \ref{tab:cSNx0}, \ref{tab:cSNx2}, \ref{tab:cSUx1} and \ref{tab:cSUx2}.) Nevertheless, tough Table \ref{tab:SNx05}, \ref{tab:SUxm2}, \ref{tab:SUxm15}, \ref{tab:SUxm1}, \ref{tab:SUxm05}, \ref{tab:SUx05}, \ref{tab:SUx1}, \ref{tab:SUx15}, \ref{tab:cSUxm1}, \ref{tab:cSUxm05}, \ref{tab:cSUx05} and \ref{tab:cSUx15} also contradict the asymptotic theory in cases of a small sample size. However, their accuracy is recovered in cases of a large sample size, indicating that the asymptotic theory would work in situations where the sample size is literally $\infty$. Some studies (i.e. \cite{BGH93} and \cite{FHMP96}) have already found similar phenomenon where curve estimation with global bandwidth tends to be oversmoothing and displays have poor performance at points of large curvature. \cite{HTF07} introduce this phenomenon as 'trimming the hills' and 'filling the valleys' in the literature of local linear regression. However, since our Standardised statistics are centred at $\mathbb{E}[f_{h_0}(x)]$, this phenomenon cannot occur. Additional simulation results show that, using optimal bandwidth $h_0$, KDEs are distributed around $\mathbb{E}[\hat{f}_{h_0}(x)]$ (See Figure \ref{fig:SN_dist_KDE_x0_hat},\ref{fig:SN_dist_KDE_x0_convo},\ref{fig:SN_dist_KDE_x15_hat},\ref{fig:SN_dist_KDE_x15_convo},\ref{fig:SU_dist_KDE_x1_hat} and \ref{fig:SU_dist_KDE_x1_convo}), so oversmoothing at the high curveature point in our simulation studies comes from the bandwidth selection. This is despite the fact that  $\hat{I}_L$'s are distributed around $I_L$ in a good manner for large sample sizes (See Figure \ref{fig:SN_dist_IL_hat},\ref{fig:SN_dist_IL_convo},\ref{fig:SU_dist_IL_hat}, and \ref{fig:SU_dist_IL_convo}). \cite{HTF07} state that one can avoid oversmoothing from 'trimming the hills' and 'filling the valleys' by using local polynomial regressions higher than second-order (for density estimation, one has to use local polynomial density \cite{CJM20} higher than third-order). However, one cannot not avoid oversmoothing from bandwidth selection in the way.  The scope of this paper is to develop higher-order approximation of KDE with global plug-in bandwidth and solving the puzzle on the curvature is out of scope. Strategies for dealing with this difficulty are discussed in Section \ref{Section:Discussion}. However, for almost all points, our expansions provide more precise approximation than the normal approximations and the Edgeworth expansion with deterministic bandwidth.

\section{Discussion and Conclusions} \label{Section:Discussion}

    This study investigated the higher-order asymptotic properties of KDE with global plug-in bandwidth. The first contribution is that we provide the Edgeworth expansion of KDE with global plug-in bandwidth up to the order of $O\{(nh_0)^{-1}+h_0^{2L}\}=O(n^{\frac{-2L}{2L+1}})$ and show that the bandwidth selection by the plug-in method starts to have an effect from on the term whose convergence rate is $O\{(nh_0)^{-1/2}h_0+h_0^{L+1}\}=O(n^{\frac{-(L+1)}{2L+1}})$ under the condition that $L_p$ is large enough. 
    Second, we generalise Theorem 3.2 of \cite{HK01}, which states that bandwidth selection via the global plug-in method has no effect on the asymptotic structure of KDE up to the order of $O\{(nh_0)^{-1/2}+h_0^{L}\}=O(n^{\frac{-L}{2L+1}})$.
    Their results limit the order of kernel functions $K(u)$ and $H(u)$ to $L=2,L_p=6$ respectively, but we show that they are valid for general orders $L$  as well under the condition that $L_p$ is large enough. 
    Third, we explore Edgeworth expansion of KDE with deterministic bandwidth in more detail than \cite{Hall91}. We show that Edgeworth expansion of Standardized KDE with deterministic bandwidth has the term of order $O\{(nh_0)^{-1/2}+h_0^{L}\}=O(n^{\frac{-L}{2L+1}})$ right after the term $\Phi(z)$ with a gap between them. After that however, the terms decrease at the rate of $O(h_0)=O(n^{\frac{-1}{2L+1}})$.
    Fourth, we weaken this condition on $L_p$ assumed by \cite{HK01} and our Theorem \ref{T:Main Theorem} and provide the Edgeworth expansion including the effect of pilot bandwidth up to the order of $O\{(nh_0)^{-1}+h_0^{2L}\}$. In this situation, the bandwidth selection via the global plug-in method possibly has an effect on the asymptotic structure of KDE even up to the order of $O\{(nh_0)^{-1/2}+h_0^L\}$ (for example, when $L=2$ and $L_p=2$).
    Finally, we consider the intersectional effect of the bandwidth selection via the global plug-in method, its accompanying pilot bandwidth, and Studentisation.

    Simulation studies show that our higher-order approximation is more precise at the point where coverage probability of normal approximation is away from $0.9500$ while less precise at the point where  normal approximation is nearly $0.9500$.

    Another implication of simulation studies is that the estimation at the points of large curvature is difficult. One possible method to avoid this problem is to use locally adaptive bandwidth. However, locally adaptive bandwidth also has disadvantages. First, selecting bandwidth at each $x$ is computationally expensive,  especially in multivariate case. Second, \cite{HK01} have shown that nonparametric bootstrap procesures for KDE with locally adaptive bandwidth lack the asymptotic refinement, while those with global bandwidth do not.  Finally, and most importantly, some authors state that locally adaptive procedures are not suited for the construction of confidence intervals; we quote \cite[p. 212]{wasserman06} [...do adaptive methods work or not? If one needs accurate function estimates and the noise level is low, then the answer is that adaptive function estimators are very effective. But if we are facing a standard nonparametric regression problem and we are interested in confidence sets, then adaptive methods do not perform significantly better than other methods such as fixed bandwidth local regression.] Another strategy employs partially adaptive bandwidth as \cite{HMT95} and in domains where global bandwidths are used, our approximations might be useful.

    As stated in Remark \ref{R:Asymptotic Bias}, centring at $\mathbb{E}\hat{f}_h(x)$ leaves asymptotic bias under standard conditions. Two standard methods to deal with asymptotic bias (\textit{debias}) are `undersmoothing' and `explicit bias reduction'. The former refers to choosing the bandwidth satisfying $\sqrt{nh}h^L\rightarrow 0$ and the latter directly estimates and removes the bias term. \cite{Hall92} examined the effect of undersmoothing and explicit bias reduction on the asymptotic structure via the Edgeworth expansion up to the order of $O\{(nh)^{-1}\}$, and stated that undersmoothing provides better coverage than explicit bias correction. After that, \cite{CCF18} have proposed alternative bias correction methods and show that thier method is comparable with undersmoothing by Edgeworth expansion up to the order of $O\{(nh)^{-1}\}$. However, the bandwidth in their expansion is still deterministic. We can interpret that \cite{HK01}, our study, and \cite{Hall92}, \cite{CCF18} studied these effects separately, that is, the pure effect of bandwidth selection and the pure effect of debias respectively. A goal for future research will be  investigating the effect of bandwidth selection and debias simultaneously, on which we are working at the moment.

    Among the recent topics in which the density estimator plays an important role is the manipulation test of regression discontinuity designs (RDD). \cite{CJM20} proposed a local polynomial density estimator for adaptability at or near the boundary points. We expect that the asymptotic structure of their estimator with the corresponding plug-in bandwidth has a similar structure to that of the KDE provided in this paper.

    One of the other possible extensions of this work is, which we are in the process of working on, is investigating the effects of cross-validation methods on the asymptotic structure. 

\section*{Acknowledgement}
    This work was supported by JSPS KAKENHI, Grant Number 19H01473, and the Joint Usage and Research Project of Institute of Economic Research, Kyoto University.
    The authors are grateful to Takahide Yanagi for the useful discussions. We would also like to thank Daisuke Kurisu, Takuya Ishihara, Masahiko Sagae, Yoshihiko Maesono, Kanta Naito, Masamune Iwasawa and the participants of the several meetings and conferences for their useful comments.

    \clearpage
    \newgeometry{margin=20truemm}
    \begin{center}
    \textbf{Supplemental Materials for 'Higher-Order Asymptotic Properties of Kernel Density Estimator with Global Plug-In and Its Accompanying Pilot Bandwidth'} (not for publication)
    \end{center}
    \begin{center}
        Shunsuke Imai$^{*}$ and Yoshihiko Nishiyama$^{\dagger}$
    \end{center}
    \begin{center}
         $^{*}$Graduate School of Economics, Kyoto University\\
        $^{\dagger}$Institute of Economic Research, Kyoto University
    \end{center}
    
    \appendix

    \section{Proofs of Results}

    \subsection{Proof of Proposition \ref{p:Linearization of Plug-In Bandwidth}}

\begin{proof} \label{Proof: Linearization of PI bandwidth}
Recall that the unknown part $I_L$ of the theoretically optimal bandwidth which minimize MISE is estimated by 
\begin{equation}
    \hat{I}_L=\binom{n}{2}^{-1}\sum_{i=1}^{n-1}\sum_{j=i+1}^nb^{-(2L+1)}H^{(2L)}\left(\frac{X_i-X_j}{b}\right)\equiv\binom{n}{2}^{-1}\sum_{i=1}^{n-1}\sum_{j=i+1}^n \hat{I}_{Lij}. \nonumber
\end{equation}
Since $\hat{I}_{L}$ has a U-statistic form, we can use Hoeffding-Decomposition,
\begin{align}
    \hat{I}_L&=\mathbb{E}\hat{I}_{Lij}+\frac{2}{n}\sum_{i=1}^n\Bigl\{\hat{I}_{Li}-\mathbb{E}\hat{I}_{Lij}\Bigl\}+\binom{n}{2}^{-1}\sum_{i=1}^{n-1}\sum_{j=1}^n\Bigl\{\hat{I}_{Lij}-\hat{I}_{Li}-\hat{I}_{Lj}+\mathbb{E}\hat{I}_{Lij}\Bigl\}, \label{I hat}
\end{align}
where $\hat{I}_{Li}=\mathbb{E}[\hat{I}_{Lij}|X_i]$. In order to examine $\hat{I}_L$, we have to compute $\mathbb{E}\hat{I}_{Lij}$ and $\hat{I}_{Li}$.
\begin{align}
    \hat{I}_{Li}=\mathbb{E}\Bigl[\hat{I}_{Lij}|X_i\Bigl]&=\int \frac{1}{b^{2L+1}}H^{(2L)}\left(\frac{X_i-x}{b}\right)f(x)dx \nonumber\\
    &=\int \frac{1}{b^{2L}}H^{(2L)}(u)f(X_i+ub)du \nonumber\\
    &=\int H(u)f^{(2L)}(X_i+ub)du \nonumber\\
    &=\int H(u)\left\{f^{(2L)}(X_i)+\frac{f^{2L+L_p}}{L_p!}(ub)^{L_p} + o(b^{L_p})\right\}du \\
    &=f^{(2L)}(X_i)+\frac{b^{L_p}}{(L_p)!}\left(\int u^{L_p}H(u)du\right)f^{(2L+L_p)}(X_i)+o_p(b^{L_p}), \label{Li}
\end{align}
where the third equality follows from integration by part and the fourth equality follows from the expansion of $f^{(2L)}(X_i+ub)$ around $X_i$. This implies
\begin{equation}
    \mathbb{E}\hat{I}_{Lij}=\mathbb{E}\bigl[f^{(2L)}(X_i)\Bigl] + \frac{\int u^{L_p}H(u)du}{(L_p)!}\mathbb{E}[f^{(2L+L_p)}(X_i)]b^{L_p} + o_p(b^{L_p}) \label{Lij}
\end{equation}
From integration by parts the first term of the right-hand side of (\ref{Lij}) is
\begin{align}
    \mathbb{E}\bigl[f^{(2L)}(X_i)\Bigl] = \int f^{(2L)}(x)f(x)dx = \int f^{(L)}(x)^2dx=I_L \label{Lij2},
\end{align}

Inserting (\ref{Li}), (\ref{Lij}) and (\ref{Lij2}) into (\ref{I hat}), we have
\begin{align}
    \hat{I}_L = I_L&+\frac{2}{n}\sum_{i=1}^n\{f^{(2L)}(X_i)-\mathbb{E}f^{(2L)}(X_i)\} \nonumber\\
    &\qquad+\frac{2}{n}\left(\int u^{L_{p}}H(u)du\right) \frac{b^{L_{p}}}{(L_p)!}\sum_{i=1}^n\{f^{(2L+L_p)}(X_i)-\mathbb{E}f^{(2L+L_p)}(X_i)\} \nonumber\\
    &\qquad+\binom{n}{2}^{-1}\sum_{i=1}^{n-1}\sum_{j=1}^n\Bigl\{\hat{I}_{Lij}-\hat{I}_{Li}-\hat{I}_{Lj}+\mathbb{E}\hat{I}_{Lij}\Bigl\}+o_p(n^{-1/2}b^{L_p}).  \label{linearized hat I}
\end{align}
Recall that Plug-In bandwidth is defined as follows,
\begin{align}
  \hat{h}&=\left(\frac{R(K)}{2LC_L^2\hat{I}_L}\right)^{\frac{1}{2L+1}}n^{-\frac{1}{2L+1}}. \label{plug in}
\end{align}
We evaluate the difference between $\hat h$ and $h_0$ using (\ref{linearized hat I}).
\begin{align}
    \hat{I}_L^{\frac{-1}{2L+1}}&=I_L^{\frac{-1}{2L+1}}-\frac{1}{2L+1}I_L^{\frac{-1}{2L+1}-1}\Biggl[\frac{2}{n}\sum_{i=1}^n\{f^{(2L)}(X_i)-\mathbb{E}f^{(2L)}(X_i)\} \nonumber\\
    &\qquad+\frac{2}{n}\left(\int u^{L_{p}}H(u)du\right) \frac{b^{L_{p}}}{(L_p)!}\sum_{i=1}^n\{f^{(2L+L_p)}(X_i)-\mathbb{E}f^{(2L+L_p)}(X_i)\} \nonumber\\
    &\qquad+\binom{n}{2}^{-1}\sum_{i=1}^{n-1}\sum_{j=1}^n\Bigl\{\hat{I}_{Lij}-\hat{I}_{Li}-\hat{I}_{Lj}+\mathbb{E}\hat{I}_{Lij}\Bigl\}\Biggl]+o_p(n^{-1/2}b^{L_p}). \nonumber
\end{align}
Inserting this expansion into (\ref{plug in}) yields
\begin{align*}
    \hat{h}
    &=h_0-\frac{h_0}{2L+1}I_L^{-1}\\
    &\qquad\times\Biggl[\frac{2}{n}\sum_{i=1}^n\{f^{(2L)}(X_i)-\mathbb{E}f^{(2L)}(X_i)\} \\
    &\quad+\frac{2}{n}\left(\int u^{L_{p}}H(u)du\right) \frac{b^{L_{p}}}{(L_p)!}\sum_{i=1}^n\{f^{(2L+L_p)}(X_i)-\mathbb{E}f^{(2L+L_p)}(X_i)\} \nonumber\\
    & \quad +\binom{n}{2}^{-1}\sum_{i=1}^{n-1}\sum_{j=1}^n\Bigl\{\hat{I}_{Lij}-\hat{I}_{Li}-\hat{I}_{Lj}+\mathbb{E}\hat{I}_{Lij}\Bigl\}\Biggl]+o_p(n^{-1/2}b^{L_p})
\end{align*}
This implies
\begin{align*}
    \frac{\hat{h}-h_0}{h_0}&=-\frac{1}{2L+1}I_L^{-1}\Biggl[\frac{2}{n}\sum_{i=1}^n\{f^{(2L)}(X_i)-\mathbb{E}f^{(2L)}(X_i)\}\\
    &\quad+\frac{2}{n}\left(\int u^{L_{p}}H(u)du\right) \frac{b^{L_{p}}}{(L_p)!}\sum_{i=1}^n\{f^{(2L+L_p)}(X_i)-\mathbb{E}f^{(2L+L_p)}(X_i)\} \nonumber\\
    &\quad+\binom{n}{2}^{-1}\sum_{i=1}^{n-1}\sum_{j=1}^n\Bigl\{\hat{I}_{Lij}-\hat{I}_{Li}-\hat{I}_{Lj}+\mathbb{E}\hat{I}_{Lij}\Bigl\}\Biggl]+o_p(n^{-1/2}b^{L_p})\\
    \implies\frac{\hat{h}-h_0}{h_0}&=\frac{-C_{PI}}{n}\sum_{i=1}^n\left(\Bigl\{f^{(2L)}(X_i)-\mathbb{E}f^{(2L)}(X_i)\Bigl\} + \frac{\int u^{L_p}H(u)du}{(L_p)!}\Bigl\{f^{2L+L_p}(X_i)-\mathbb{E}f^{2L+L_p}(X_i)\Bigl\}\right)\\
    &\qquad-\frac{C_{PI}}{2}\binom{n}{2}^{-1}\sum_{i=1}^{n-1}\sum_{j=i+1}^n\Bigl\{\hat{I}_{Lij}-\hat{I}_{Li}-\hat{I}_{Lj}+\mathbb{E}\hat{I}_{Lij}\Bigl\} +o_p\{(nh_0)^{-1}\}.
\end{align*}
\end{proof}

\subsection{Proof of Theorem \ref{T: Hall and Kang}}
\begin{proof} \label{Proof:Hall and Kang}
 In view of (\ref{def SPI}), if the following evaluation is correct,
\begin{align*}
    \mathbb{E}\left|\sqrt{nh}\left(\frac{\hat{h}-h_0}{h_0}\right)\Gamma_{KDE_1}\right|=o\{(nh_0)^{-1/2}\},\quad \mathbb{E}\left|S_{h_0}(x)\left(\frac{\hat{h}-h_0}{h_0}\right)\right|=o\{(nh_0)^{-1/2}\}
\end{align*}
then bandwidth selection has no effect on the asymptotic structure up to the order of $O\{(nh_0)^{-1/2}\}$.
From Cauchy-Schwarz Inequality
\begin{align*}
    \mathbb{E}\left|\sqrt{nh_0}\left(\frac{\hat{h}-h_0}{h_0}\right)\Gamma_{KDE_1}\right|\le\sqrt{nh_0}\left\{\mathbb{E}\left|\frac{\hat{h}-h_0}{h_0}\right|^2\mathbb{E}\Bigl|\Gamma_{KDE_1}\Bigl|^2\right\}^{1/2}
\end{align*}
Since under the Assumption \ref{A:L_p order 1}, $\hat{h}-h_0/h_0$ has the asymptotic linear form, straightforward calculation gives 
\begin{align*}
    \mathbb{E}\left|\frac{\hat{h}-h_0}{h_0}\right|^2=O(n^{-1}).
\end{align*}
Next, we evaluate $\mathbb{E}|\Gamma_{KDE_1}|^2$.
\begin{align*}
    \mathbb{E}\Bigl|\Gamma_{KDE_1}\Bigl|^2&=\frac{1}{(nh_0)^2}\mathbb{E}\Biggl[\sum_{i=1}^n\sum_{j\neq i}^n\left\{K'_{i,h_0}(x)u_{i,h_0}(x)+K_{i,h_0}(x)\right\}\left\{K'_{i,h_0}(x)u_{i,h_0}(x)+K_{i,h_0}(x)\right\}\Biggl]\\
    &\qquad+\frac{1}{(nh_0)^2}\mathbb{E}\left[\sum_{i=1}^n\left\{K'_{i,h_0}(x)u_{i,h_0}(x)+K_{i,h_0}(x)\right\}^2\right]\\
    &=\frac{1}{h_0^2}\mathbb{E}\Biggl[\left\{K'_{i,h_0}(x)u_{i,h_0}(x)+K_{i,h_0}(x)\right\}\left\{K'_{i,h_0}(x)u_{i,h_0}(x)+K_{i,h_0}(x)\right\}\Biggl]+O\{(nh_0)^{-1}\}\\
    &=\frac{1}{h_0^2}\left(\int \left\{K'\left(\frac{z_1-x}{h_0}\right)\left(\frac{z_1-x}{h_0}\right)+K\left(\frac{z_1-x}{h_0}\right)\right\}f(z_1)dz_1\right)^2+O\{(nh_0)^{-1}\}\\
    &=\left(\int K'(u)uf(x+uh_0)du+\int K(u)f(x+uh_0)du\right)^2+O\{(nh_0)^{-1}\}\\
    &=\Biggl(-\int K(u)f(x+uh_0)du-\int K(u)uf'(x+uh_0)h_0du+\int K(u)f(x+uh_0)du\Biggl)^2+O\{(nh_0)^{-1}\}\\
    &=\left(-\int K(u)uf'(x+uh_0)h_0du\right)^2+O\{(nh_0)^{-1})\\
    &=O(h_0^{2L})+O\{(nh_0)^{-1}\}.
\end{align*}
The fifth equality follows from integration by part of the first term and Assumption \ref{A:Ku}, and the final equality follows from the expansion of $f'(x+uh_0)$ around $h_0=0$ and Assumption \ref{A:2L+L_p times differentiable},\ref{A: Kernel ''}.
Therefore form Cauchy-Schwarz inequality,
\begin{align*}
    \mathbb{E}&\left|\sqrt{nh_0}\left(\frac{\hat{h}-h_0}{h_0}\right) \Gamma_{KDE}\right|\le\sqrt{nh_0}\left\{\mathbb{E}\left|\frac{\hat{h}-h_0}{h_0}\right|^2\mathbb{E}\Bigl|\Gamma_{KDE}\Bigl|^2\right\}^{1/2}\\
    &=O(n^{1/2}h_0^{1/2})\Bigl(O(n^{-1})O(h_0^{2L}+(nh_0)^{-1})\Bigl)^{1/2}=O(h_0^{L+\frac{1}{2}}+n^{-1/2})=o\{(nh_0)^{-1/2}\}.
\end{align*}
Similar to above evaluation, Cauchy-Schwarz inequality gives $\mathbb{E}\left|S(x)\left(\frac{\hat{h}-h_0}{h_0}\right)\right|=O(n^{-1/2})=o\{(nh_0)^{-1/2}\}$.
Therefore bandwidth selection via Plug-In Method has no effect on the asymptotic structure up to the order of $O\{(nh_0)^{-1/2}\}$.
\end{proof}
\subsection{Proof of Theorem \ref{T:Main Theorem}}
\begin{proof} \label{Proof:Main Theorem}
From Proposition \ref{p:Higher-Order Expansion of KDE with Data-Driven Bandwidth} and Lemma \ref{L:CS inequality 2}, we have, 
\begin{align}
    \sqrt{n\hat{h}}\Bigl({\hat{f}_{\hat{h}}(x)}-\mathbb{E}\hat{f}_{h_0}(x)\Bigl)&=\sqrt{nh_0}\Bigl(\hat{f}_{h_0}(x)-\mathbb{E}\hat{f}_{h_0}(x)\Bigl)-\sqrt{nh_0}\left(\frac{\hat{h}-h_0}{h_0}\right)\Gamma_{KDE_1} +\frac{1}{2}S_{h_0}(x)\left(\frac{\hat{h}-h_0}{h_0}\right)+o_p\{(nh_0)^{-1}\}. \nonumber
\end{align}
 Noting that we provide Theorem \ref{T:Main Theorem} under Assumption \ref{A:L_p order 2} and this assumption guarantees that the quadratic term of $(\hat{h}-h_0)/h_0$ is negligible, Proposition \ref{p:Linearization of Plug-In Bandwidth} provides the expansion of plug-in bandwidth as follows.
\begin{align}
    \frac{\hat{h}-h_0}{h_0}&=\frac{-C_{PI}}{n}\sum_{i=1}^n\Bigl\{f^{(2L)}(X_i)-\mathbb{E}f^{(2L)}(X_i)\Bigl\} + o_p\{(nh_0)^{-1}\}. \nonumber
\end{align}

Define
\begin{align}
    S_i&\equiv \mu_{20}^{-1/2}\left(K\left(\frac{X_i-x}{h_0}\right)-\mathbb{E}K\left(\frac{X_i-x}{h_0}\right)\right), \nonumber\\
    \Gamma_i&\equiv K'\left(\frac{X_i-x}{h_0}\right)\left(\frac{X_i-x}{h_0}\right)+K\left(\frac{X_i-x}{h_0}\right)-\mathbb{E}\left[K'\left(\frac{X_i-x}{h_0}\right)\left(\frac{X_i-x}{h_0}\right)+K\left(\frac{X_i-x}{h_0}\right)\right], \nonumber\\
    \mathcal{L}_i&\equiv f^{(2L)}(X_i)-\mathbb{E}f^{(2L)}(X_i). \nonumber
\end{align}
Recalling that $S_{PI}(x)$ is defined as (\ref{def SPI}), we have from Lemma \ref{L:EGamma1}, 
\begin{align}
    S_{PI}(x)&=\frac{\sqrt{n\hat{h}}\Bigl({\hat{f}_{\hat{h}}(x)}-\mathbb{E}\hat{f}_{h_0}(x)\Bigl)}{\mu_{20}^{1/2}} \nonumber\\
    &=\frac{\sqrt{nh_0}\Bigl(\hat{f}_{h_0}(x)-\mathbb{E}\hat{f}_{h_0}(x)\Bigl)}{\mu_{20}^{1/2}}-\frac{\sqrt{nh_0}\left(\frac{\hat{h}-h_0}{h_0}\right)\mathbb{E}\Gamma_{KDE_1}}{\mu_{20}^{1/2}} \nonumber\\
    &\qquad-\frac{\left(\frac{\hat{h}-h_0}{h_0}\right)\sqrt{nh_0}(\Gamma_{KDE_1}-\mathbb{E}\Gamma_{KDE_1})}{\mu_{20}^{1/2}}+\frac{\sqrt{nh_0}\Bigl(\hat{f}_{h_0}(x)-\mathbb{E}\hat{f}_{h_0}(x)\Bigl)\left(\frac{\hat{h}-h_0}{h_0}\right)}{2\mu_{20}^{1/2}}+o_p\{(nh_0)^{-1}\} \nonumber\\
    &=\frac{1}{\sqrt{nh_0}}\sum_{i=1}^nS_i+\sqrt{nh_0}\sum_{l=0}^{L-1}\frac{C_{PI}C_{\Gamma,l}(x)h_0^{L+l}}{n}\sum_{i=1}^n\frac{\mathcal{L}_i}{\mu_{20}^{1/2}} \nonumber\\
    &\qquad+\left(\frac{1}{\sqrt{nh_0}}\sum_{i=1}^n\frac{\Gamma_i}{\mu_{20}^{1/2}}\right)\left(\frac{C_{PI}}{n}\sum_{i=1}^n\mathcal{L}_i\right)-\frac{1}{2}\left(\frac{1}{\sqrt{nh_0}}\sum_{i=1}^nS_i\right)\left(\frac{C_{PI}}{n}\sum_{i=1}^n\mathcal{L}_i\right)+o_p\{(nh_0)^{-1}\} \nonumber\\
    &=\frac{1}{\sqrt{nh_0}}\sum_{i=1}^nS_i+\frac{C_{PI}h_0^{\frac{2L+1}{2}}}{n^{1/2}\mu_{20}^{1/2}}\sum_{i=1}^n\mathcal{L}_i \sum_{l=0}^{L-1} C_{\Gamma,l}(x)h_0^l \nonumber\\
    &\qquad+\frac{C_{PI}}{n^{3/2}h_0^{1/2}\mu_{20}^{1/2}}\sum_{i=1}^n\sum_{j\neq i}^n\Gamma_i\mathcal{L}_j+\frac{C_{PI}}{n^{3/2}h_0^{1/2}\mu_{20}^{1/2}}\sum_{i=1}^n\Gamma_i\mathcal{L}_i\nonumber\\&
    \qquad-\frac{C_{PI}}{2n^{3/2}h_0^{1/2}}\sum_{i=1}^n\sum_{j\neq i}^nS_i\mathcal{L}_j-\frac{C_{PI}}{2n^{3/2}h_0^{1/2}}\sum_{i=1}^nS_i\mathcal{L}_i+o_p\{(nh_0)^{-1}\} \nonumber\\
    &\equiv S(x)+\Lambda_1(x)+\Lambda_2(x)+\Lambda_3(x)+\Lambda_4(x)+\Lambda_5(x)+o_p\{(nh_0)^{-1}\}. \label{SPI}
\end{align}
Define $F_{PI}(z)$ and $\Tilde{F}_{PI}(z)$ as follows,
\begin{align}
    F_{PI}(z)&=\mathbb{P}\Bigl(S_{PI}(x)\le z\Bigl), \nonumber\\
    \Tilde{F}_{PI}(z)&=\Phi(z)+\phi(z)\Biggl[(nh_0)^{-1/2}p_1(z)+(nh_0)^{-1}p_2(z)+\sum_{l=0}^{L-1}h_0^{L+l+1}p_{3,l}(z)+n^{-1/2}h_0^{1/2}p_4(z)\Biggl].\nonumber
\end{align}
To show the Edgeworth expansion is valid, we have to confirm $\sup_{z\in\mathbb{R}}\left|F_{PI}(z)-\Tilde{F}_{PI}(z)\right|=o\{(nh_0)^{-1}\}$. First, we evaluate the remainder term.
\begin{align}
    \sup_{z\in\mathbb{R}}&\left|F_{PI}(z)-\Tilde{F}_{PI}(z)\right|\le \sup_{z\in \mathbb{R}}\left|\mathbb{P}\Bigl(S(x)+\Lambda_1(x)+\Lambda_2(x)+\Lambda_3(x)+\Lambda_4(x)+\Lambda_5(x)\le z\Bigl)-\Tilde{F}_{PI}(z)\right| \nonumber\\
    &+\mathbb{P}\left(\left|S_{PI}(x)-\Bigl(S(x)+\Lambda_1(x)+\Lambda_2(x)+\Lambda_3(x)+\Lambda_4(x)+\Lambda_5(x)\Bigl)\right|\ge a_n\right)+O(a_n^{-1}) \nonumber
\end{align}
where $a_n=nh_0(\log n)$. Since
\[\left|S_{PI}(x)-\Bigl(S(x)+\Lambda_1(x)+\Lambda_2(x)+\Lambda_3(x)+\Lambda_4(x)+\Lambda_5(x)\Bigl)\right|=o_p\{(nh_0)^{-1}\},\]
 we have 
\begin{align}
    &\mathbb{P}\left(\left|S_{PI}(x)-\Bigl(S(x)+\Lambda_1(x)+\Lambda_2(x)+\Lambda_3(x)+\Lambda_4(x)+\Lambda_5(x)\Bigl)\right|\ge a_n\right) =O\{(nh_0)^{-1}a_n^{-1}\}=o\{(nh_0)^{-1}\}. \nonumber
\end{align}
Obviously, $O(a_n^{-1})=o\{(nh_0)^{-1}\}$. Then, we only need to evaluate
\[\sup_{z\in \mathbb{R}}\left|\mathbb{P}\Bigl(S(x)+\Lambda_1(x)+\Lambda_2(x)+\Lambda_3(x)+\Lambda_4(x)+\Lambda_5(x)\le z\Bigl)-\Tilde{F}_{PI}(z)\right|.\]

Define $\chi_{PI}(t)$ and  $\Tilde{\chi}_{PI}(t)$ as follows,
\begin{align}
    \chi_{PI}(t)&\equiv\mathbb{E}\left[exp\left\{it\Bigl(S(x)+\Lambda_1(x)+\Lambda_2(x)+\Lambda_3(x)+\Lambda_4(x)+\Lambda_5(x)\Bigl)\right\}\right], \nonumber \\
    \Tilde{\chi}_{PI}(t)&\equiv\exp\left(\frac{-t^2}{2}\right) \Biggl[\Biggl\{1+\frac{\mu_{30}\mu_{20}^{-3/2}}{6n^{1/2}h_0^{1/2}}(it)^3+\frac{\mu_{40}\mu_{20}^{-2}}{24nh_0}(it)^4+\frac{\mu_{30}^2\mu_{20}^{-3}}{72nh_0}(it)^6\Biggl\} \nonumber\\
    &\qquad\qquad+C_{PI}\rho_{11}\mu_{20}^{-1}\left(\sum_{l=0}^{L-1}C_{\Gamma,l}(x)h_0^{L+l+1}\right)(it)^2+C_{PI}\frac{\rho_{11}\xi_{11}\mu_{20}^{-3/2}h_0^{1/2}}{n^{1/2}}(it)^3-C_{PI}\frac{\mu_{20}^{-1/2}\rho_{11}h_0^{1/2}}{2n^{1/2}}\{(it)^3+(it)\}\Biggl]. \nonumber
\end{align}
From \cite{Esseen45} smoothing lemma,
\begin{align}
    &\sup_{z\in \mathbb{R}}\left|P\Bigl(S(x)+\Lambda_1(x)+\Lambda_2(x)+\Lambda_3(x)+\Lambda_4(x)+\Lambda_5(x)\leqq z\Bigl)-\Tilde{F}_{PI}(z)\right|\nonumber\\
    &\lesssim \int_{-n^{\frac{2L}{2L+1}}\log n}^{n^{\frac{2L}{2L+1}}\log n}\left|\frac{\chi_{PI}(t)-\Tilde{\chi}_{PI}(t)}{t}\right|dt+O\left(\frac{1}{n^\frac{2L}{2L+1}\log n}\right)\nonumber\\
    &\le \int_{-p}^p\left|\frac{\chi_{PI}(t)-\Tilde{\chi}_{PI}(t)}{t}\right|dt +\int_{p\leqq |t|\leqq n^{\frac{2L}{2L+1}}\log n}\left|\frac{\chi_{PI}(t)}{t}\right|dt+\int_{p\leqq |t|\le n^{\frac{2L}{2L+1}}\log n}\left|\frac{\Tilde{\chi}_{PI}(t)}{t}\right|dt+o\{(nh_0)^{-1}\} \nonumber\\
    &\le\int_{-p}^p\left|\frac{\chi_{PI}(t)-\Tilde{\chi}_{PI}(t)}{t}\right|dt +\int_{p\leqq |t|\leqq n^{\frac{2L}{2L+1}}\log n}\left|\frac{\chi_{PI}(t)}{t}\right|dt+\int_{p\leqq |t|}\left|\frac{\Tilde{\chi}_{PI}(t)}{t}\right|dt+o\{(nh_0)^{-1}\} \nonumber\\
    &\equiv (A)+(B)+(C)+o\{(nh_0)^{-1}\} \label{esseen plug-in}
\end{align}
where $p=\min\left\{\frac{n^{1/2}h_0^{1/2}}{\mu_{20}^{-3/2}\mu_{30}},\log n\right\}$.
To prove the validity of the Edgeworth expansion, we show that each term of (\ref{esseen plug-in}) has the convergence rate $o\{(nh_0)^{-1}\}$.

In order to evaluate $(A)$, we represent $\chi_{PI}(t)$ as $\Tilde{\chi}_{PI}(t)$ plus a remainder. From Lemmas \ref{L: AM Lambda 1}, \ref{L:AM Lambda 2}, \ref{L:AM Lambda 3}, \ref{L:AM Lambda 4}, and \ref{L:AM Lambda 5},
\begin{align}
    \chi_{PI}(t)&=\mathbb{E}\left[e^{it\Bigl(S(x)+\Lambda_1(x)+\Lambda_2(x)+\Lambda_3(x)+\Lambda_4(x)+\Lambda_5(x)\Bigl)}\right] \nonumber\\
    &=\mathbb{E}\left[e^{itS(x)}\Bigl\{1+it\Lambda_1(x)\Bigl\}\Bigl\{1+it\Lambda_2(x)\Bigl\}\Bigl\{1+it\Lambda_4(x)\Bigl\}\Bigl\{1+it\Lambda_5(x)\Bigl\}\right] \nonumber\\
    &\qquad+O(t^2\mathbb{E}|\Lambda_1(x)|^2)+O(t^2\mathbb{E}|\Lambda_2(x)|^2)+O(|t|\mathbb{E}|\Lambda_3(x)|)+O(t^2\mathbb{E}|\Lambda_4(x)|^2)+O(t^2\mathbb{E}|\Lambda_5(x)|^2) \nonumber\\
    &=\mathbb{E}\left[e^{itS(x)}\Bigl\{1+it\Lambda_1(x)+it\Lambda_2(x)+it\Lambda_4(x)+it\Lambda_5(x)\Bigl\}\right]\nonumber\\
    &\qquad+O(t^2\mathbb{E}|\Lambda_1(x)|^2)+O(t^2\mathbb{E}|\Lambda_2(x)|^2)+O(|t|\mathbb{E}|\Lambda_3(x)|)+O(t^2\mathbb{E}|\Lambda_4(x)|^2)+O(t^2\mathbb{E}|\Lambda_5(x)|^2) \nonumber\\
    &\qquad +O(t^2\mathbb{E}|\Lambda_1(x)\Lambda_2(x)|)+O(t^2\mathbb{E}|\Lambda_1(x)\Lambda_4(x)|)+O(t^2\mathbb{E}|\Lambda_1(x)\Lambda_5(x)|) \nonumber\\
    &\qquad+O(t^2\mathbb{E}|\Lambda_2(x)\Lambda_4(x)|)+O(t^2\mathbb{E}|\Lambda_2(x)\Lambda_5(x)|)+O(t^2\mathbb{E}|\Lambda_4(x)\Lambda_5(x)|)\nonumber\\
    &\equiv(\text{I})+(\text{II})+(\text{III})+(\text{IV})+(\text{V}) \nonumber\\
    &\qquad+O(t^2\mathbb{E}|\Lambda_1(x)|^2)+O(t^2\mathbb{E}|\Lambda_2(x)|^2)+O(|t|\mathbb{E}|\Lambda_3(x)|)+O(t^2\mathbb{E}|\Lambda_4(x)|^2)+O(t^2\mathbb{E}|\Lambda_5(x)|^2) \nonumber\\
    &\qquad +O(t^2\mathbb{E}|\Lambda_1(x)\Lambda_2(x)|)+O(t^2\mathbb{E}|\Lambda_1(x)\Lambda_4(x)|)+O(t^2\mathbb{E}|\Lambda_1(x)\Lambda_5(x)|) \nonumber\\
    &\qquad+O(t^2\mathbb{E}|\Lambda_2(x)\Lambda_4(x)|)+O(t^2\mathbb{E}|\Lambda_2(x)\Lambda_5(x)|)+O(t^2\mathbb{E}|\Lambda_4(x)\Lambda_5(x)|)\nonumber\\
    &=(\text{I})+(\text{II})+(\text{III})+(\text{IV})+(\text{V}) \nonumber\\
    &\qquad+O(t^2h_0^{2L+1})+O(t^2n^{-1})+O(|t|\{n^{-1/2}h_0^{\frac{2L+1}{2}}+n^{-1}\})+O(t^2n^{-1}h_0)+O(t^2n^{-1}h_0) \nonumber\\
    &\qquad +O(t^2n^{-1/2}h_0^{\frac{2L+1}{2}})+O(t^2n^{-1/2}h_0^{L+1})+O(t^2n^{-1/2}h_0^{L+1})\nonumber\\
    &\qquad +O(t^2n^{-1}h_0^{1/2})+O(t^2n^{-1}h_0^{1/2})+O(t^2n^{-1}h_0) \nonumber\\
    &=(\text{I})+(\text{II})+(\text{III})+(\text{IV})+(\text{V}) + O(t^2n^{-1}) + O(|t|n^{-1}),
\end{align}
where the fourth equality follows from Lemmas \ref{L:AM Lambda 12}, \ref{L:AM Lambda 14},\ref{L:AM Lambda 15}, \ref{L:AM Lambda 24}, \ref{L:AM Lambda 25}, and \ref{L:AM Lambda 45} and the final equality uses $h_0=O(n^{-1/(2L+1)})$.

Define $\gamma(t)=\mathbb{E}\left[e^{\frac{it}{\sqrt{nh}}S_i}\right]$. We have
\begin{align}
    (\text{I})=\mathbb{E}\left[e^{\frac{it}{\sqrt{nh}}\sum_{i=1}^nS_i}\right]=\mathbb{E}\left[e^{\frac{it}{\sqrt{nh}}S_1}\right]^n=\gamma(t)^n, \label{Fist compo}
\end{align}
from, Lemma \ref{S^k L},
\begin{align}
    (\text{II})&=\mathbb{E}\left[e^{itS(x)}(it\Lambda_1(x))\right] \nonumber\\
    &=\gamma(t)^{n-1}\frac{C_{PI}h_0^{\frac{2L+1}{2}}}{n^{1/2}\mu_{20}^{1/2}}n\mathbb{E}\left[e^{\frac{it}{\sqrt{nh_0}}S_1}\mathcal{L}_1\right]\left(\sum_{l=0}^{L-1}C_{\Gamma,l}(x)h_0^l\right)(it) \nonumber\\
    &=\gamma(t)^{n-1}\frac{C_{PI}n^{1/2}h_0^{\frac{2L+1}{2}}}{\mu_{20}^{1/2}}\mathbb{E}\left[\Bigl\{1+\frac{it}{(nh_0)^{1/2}}S_1+\frac{(it)^2}{2nh_0}S_1^2\Bigl\}\mathcal{L}_1\right]\left(\sum_{l=0}^{L-1}C_{\Gamma,l}(x)h_0^l\right)(it)+o\{(nh)^{-1}\}\nonumber\\
    &=\gamma(t)^{n-1}\frac{C_{PI}n^{1/2}h_0^{\frac{2L+1}{2}}}{\mu_{20}^{1/2}}\frac{1}{(nh_0)^{1/2}}\mathbb{E}[S_1\mathcal{L}_1]\left(\sum_{l=0}^{L-1}C_{\Gamma,l}(x)h_0^l\right)(it)^2 \nonumber\\
    &\qquad+O(n^{-1/2}h_0^{(2L+1)/2})O(n)O(|t|^3n^{-1}h_0^{-1})O(h_0) \nonumber\\
    &=\gamma(t)^{n-1}\frac{C_{PI}}{\mu_{20}^{1/2}}\mathbb{E}[S_1\mathcal{L}_1]\left(\sum_{l=0}^{L-1}C_{\Gamma,l}(x)h_0^{L+l}\right)(it)^2+O(|t|^3n^{-1/2}h_0^{\frac{(2L+1)}{2}}) \nonumber\\
    &=\gamma(t)^{n-1}\frac{C_{PI}}{\mu_{20}^{1/2}}\mathbb{E}[S_1\mathcal{L}_1]\left(\sum_{l=0}^{L-1}C_{\Gamma,l}(x)h_0^{L+l}\right)(it)^2+O(|t|^3n^{-1}),\label{Second compo}
\end{align}
from Lemma \ref{S^k L} and \ref{S^k gamma},
\begin{align}
    (\text{III})&=\mathbb{E}\left[e^{itS(x)}(it\Lambda_2(x))\right] \nonumber\\
    &=\gamma(t)^{n-2}\frac{C_{PI}n(n-1)}{n^{3/2}h_0^{1/2}\mu_{20}^{1/2}}\mathbb{E}\left[e^{\frac{it}{\sqrt{nh_0}}(S_1+S_2)}\Gamma_1\mathcal{L}_2\right](it) \nonumber\\
    &=\gamma(t)^{n-2}\frac{C_{PI}n(n-1)}{n^{3/2}h_0^{1/2}\mu_{20}^{1/2}}\mathbb{E}\Biggl[\Bigl\{1+\frac{it}{(nh_0)^{1/2}}(S_1+S_2) \nonumber\\
    &\qquad+\frac{(it)^2}{2nh_0}(S_1+S_2)^2+\frac{(it)^3}{6(nh_0)^{3/2}}(S_1+S_2)^3\Bigl\}\Gamma_1\mathcal{L}_2\Biggl](it)+o\{(nh)^{-1}\} \nonumber\\
    &=\gamma(t)^{n-2}\frac{C_{PI}n(n-1)}{n^{5/2}h_0^{3/2}\mu_{20}^{1/2}}\mathbb{E}[S_1\Gamma_1]\mathbb{E}[S_2\mathcal{L}_2](it)^3+O(n^2)O(t^4n^{-3}h_0^{-2})O(h_0)O(h_0) \nonumber\\
    &=\gamma(t)^{n-2}\frac{C_{PI}n(n-1)}{n^{5/2}h_0^{3/2}\mu_{20}^{1/2}}\mathbb{E}[S_1\Gamma_1]\mathbb{E}[S_2\mathcal{L}_2](it)^3+O(t^4n^{-1}) \nonumber\\
    &=\gamma(t)^{n-2}\frac{C_{PI}}{n^{1/2}h_0^{3/2}\mu_{20}^{1/2}}\mathbb{E}[S_1\Gamma_1]\mathbb{E}[S_2\mathcal{L}_2](it)^3+O(t^4n^{-1}), \label{Third compo}
\end{align}
from Lemma \ref{S^k L},
\begin{align}
    (\text{IV})&=\mathbb{E}\left[e^{itS(x)}(it\Lambda_4(x))\right] \nonumber\\
    &=\frac{-C_{PI}n(n-1)}{2n^{3/2}h_0^{1/2}}\gamma(t)^{n-2}\mathbb{E}\left[e^{\frac{it}{\sqrt{nh_0}}(S_1+S_2)}S_1\mathcal{L}_2\right](it) \nonumber\\
    &=\frac{-C_{PI}n^{1/2}}{2h_0^{1/2}}\gamma(t)^{n-2}\mathbb{E}\Biggl[\Biggl\{1+\frac{it}{(nh_0)^{1/2}}(S_1+S_2)\nonumber\\
    &\qquad+\frac{(it)^2}{2nh_0}(S_1+S_2)^2+\frac{(it)^3}{6(nh_0)^{3/2}}(S_1+S_2)^3\Biggl\}S_1\mathcal{L}_2\Biggl](it)+o\{(nh)^{-1}\}\nonumber\\
    &=\frac{-C_{PI}}{2n^{1/2}h_0^{3/2}}\gamma(t)^{n-2}\mathbb{E}[S_1^2]\mathbb{E}[S_2\mathcal{L}_2](it)^3+O(n^{1/2}h_0^{-1/2})O(t^4n^{-3/2}h_0^{-3/2})O(h_0)O(h_0)\nonumber\\
    &=\frac{-C_{PI}}{2n^{1/2}h_0^{1/2}}\gamma(t)^{n-2}\mathbb{E}[S_2\mathcal{L}_2](it)^3 + O(t^4n^{-1})\label{Fourth compo}
\end{align}
where the final equality uses $\mathbb{E}[S_1^2]=h_0$, and from Lemma \ref{S^k L},
\begin{align}
    (\text{V})&= \mathbb{E}\left[e^{itS(x)}(it\Lambda_5(x))\right] \nonumber\\
    &= \frac{-C_{PI}n}{2n^{3/2}h_0^{1/2}}\gamma(t)^{n-1}\mathbb{E}\left[e^{\frac{it}{\sqrt{nh_0}}S_1}S_1\mathcal{L}_1\right](it) \nonumber\\
    &= \frac{-C_{PI}}{2n^{1/2}h_0^{1/2}}\gamma(t)^{n-1}\mathbb{E}\left[\left\{1+\frac{it}{\sqrt{nh_0}}S_1\right\}S_1\mathcal{L}_1\right](it)+o\{(nh)^{-1}\} \nonumber\\
    &= \frac{-C_{PI}}{2n^{1/2}h_0^{1/2}}\gamma(t)^{n-1}\mathbb{E}\left[S_1\mathcal{L}_1\right](it)+O(n^{-1/2}h_0^{-1/2})O(t^2n^{-1/2}h_0^{-1/2})O(h_0) \nonumber\\
    &= \frac{-C_{PI}}{2n^{1/2}h_0^{1/2}}\gamma(t)^{n-1}\mathbb{E}\left[S_1\mathcal{L}_1\right](it) + O(t^2n^{-1})\label{Fifth compo}
\end{align}
then
\begin{align}
    \chi_{PI}(t)&=(\text{I})+(\text{II})+(\text{III})+(\text{IV})+(\text{V})+O(t^2n^{-1})+O(|t|n^{-1}) \nonumber\\
    &=\gamma(t)^n+\gamma(t)^{n-1}\frac{C_{PI}}{\mu_{20}^{1/2}}\mathbb{E}[S_1\mathcal{L}_1]\left(\sum_{l=0}^{L-1}C_{\Gamma,l}h_0^{L+l}\right)(it)^2 \nonumber\\
    &\qquad+\gamma(t)^{n-2}\frac{C_{PI}}{\mu_{20}^{1/2}}\frac{1}{n^{1/2}h_0^{3/2}}\mathbb{E}[S_1\Gamma_1]\mathbb{E}[S_2\mathcal{L}_2](it)^3 \nonumber\\
    &\qquad-\Bigl\{\gamma(t)^{n-2}(it)^3+\gamma(t)^{n-1}(it)\Bigl\}\frac{C_{PI}}{2n^{1/2}h_0^{1/2}}\mathbb{E}[S_1\mathcal{L}_1] \nonumber\\
    &\qquad+O\left((|t|+t^2+|t|^3+t^4)n^{-1}\right). \nonumber
\end{align}
For $m=0,1,2$, by \cite[p535-536]{Feller71},
\begin{align}
    \gamma(t)^{n-m}&=\exp\left(\frac{-t^2}{2}\right)\Biggl\{1+\frac{\mu_{30}\mu_{20}^{-3/2}}{6n^{1/2}h_0^{1/2}}(it)^3+\frac{\mu_{40}\mu_{20}^{-2}}{24nh_0}(it)^4+\frac{\mu_{30}^2\mu_{20}^{-3}}{72nh_0}(it)^6\Biggl\}+o\left((nh_0)^{-1}(t^4+|t|^9)e^{-t^2/4}\right). \nonumber
\end{align}
By (\ref{Fist compo}), (\ref{Second compo}), (\ref{Third compo}), (\ref{Fourth compo}) and (\ref{Fifth compo}), noting $\mathbb{E}[S_1\mathcal{L}_1]=h_0\mu_{20}^{-1/2}\rho_{11}$, $\mathbb{E}[S_1\Gamma_1]=h_0\mu_{20}^{-1/2}\xi_{11},$ and $\mathbb{E}[S_2\mathcal{L}_2]=h_0\mu_{20}^{-1/2}\rho_{11}$,
\begin{align}
    \chi_{PI}(t)&=\exp\left(\frac{-t^2}{2}\right)\Biggl[\Biggl\{1+\frac{\mu_{30}\mu_{20}^{-3/2}}{6n^{1/2}h_0^{1/2}}(it)^3+\frac{\mu_{40}\mu_{20}^{-2}}{24nh_0}(it)^4+\frac{\mu_{30}^2\mu_{20}^{-3}}{72nh_0}(it)^6\Biggl\} \nonumber\\
    &\qquad+C_{PI}\rho_{11}\mu_{20}^{-1}\left(\sum_{l=0}^{L-1}C_{\Gamma,l}(x)h_0^{L+1+l}\right)(it)^2+C_{PI}\frac{\rho_{11}\xi_{11}\mu_{20}^{-3/2}h_0^{1/2}}{n^{1/2}}(it)^3 -C_{PI}\frac{\mu_{20}^{-1/2}\rho_{11}h_0^{1/2}}{2n^{1/2}}\{(it)^3+(it)\}\Biggl] \nonumber\\
    &\qquad+O\left((|t|+t^2+|t|^3+t^4)n^{-1}\right)+o\left((nh_0)^{-1}(t^4+|t|^9)e^{-t^2/4}\right) \nonumber\\
    &=\Tilde{\chi}_{PI}(t)+O\left((|t|+t^2+|t|^3+t^4)n^{-1}\right)+o\left((nh_0)^{-1}(t^4+|t|^9)e^{-t^2/4}\right). \nonumber
\end{align}
This implies
\begin{align}
    (A)=\int_{-p}^p\left|\frac{\chi_{PI}(t)-\Tilde{\chi}_{PI}(t)}{t}\right|dt=o\{(nh_0)^{-1}\} \nonumber
\end{align}

Next, we confirm $(B)=o\{(nh_0)^{-1}\}$, for $p\leqq|t|\leqq n^{\frac{2L}{2L+1}}\log n$.
Define
\begin{align}
    S(x;m)&\equiv\frac{1}{n^{1/2}h_0^{1/2}}\sum_{i=1}^mS_i, \nonumber\\
    \Lambda_1(x;m)&\equiv\frac{C_{PI}h_0^{\frac{2L+1}{2}}}{n^{1/2}\mu_{20}^{1/2}}\sum_{i=1}^m\mathcal{L}_i\left(\sum_{l=0}^{L-1}C_{\Gamma,l}(x)h_0^l\right), \nonumber\\
    \Lambda_2(x;m)&\equiv\frac{C_{PI}}{n^{3/2}h_0^{1/2}\mu_{20}^{1/2}}\sum_{i=1}^m\sum_{j\neq i}^m\Gamma_i\mathcal{L}_j, \nonumber\\
    \Lambda_4(x;m)&\equiv-\frac{C_{PI}}{2n^{3/2}h_0^{1/2}}\sum_{i=1}^m\sum_{j\neq i}^mS_i\mathcal{L}_j \nonumber\\
    \Lambda_5(x;m)&\equiv-\frac{C_{PI}}{2n^{3/2}h_0^{1/2}}\sum_{i=1}^mS_i\mathcal{L}_i \nonumber
\end{align}
then
\begin{align}
    \left|\chi_{PI}(t)\right|&=|\mathbb{E}e^{it(S(x)+\Lambda_1(x)+\Lambda_2(x)+\Lambda_3(x)+\Lambda_4(x)+\Lambda_5(x))}| \nonumber\\
    &\le|\mathbb{E}e^{it(S(x)+\Lambda_1(x)+\Lambda_2(x)+\Lambda_4(x)+\Lambda_5(x)}|+O(|t||\mathbb{E}\Lambda_3(x)|) \nonumber\\
    &\le\Bigl|\mathbb{E}e^{it(S(x)+(\Lambda_1(x)-\Lambda_1(x;m))+(\Lambda_2(x)-\Lambda_2(x;m))+(\Lambda_4(x)-\Lambda_4(x;m))+(\Lambda_5(x)-\Lambda_5(x;m)))} \nonumber\\
    &\qquad\times\Bigl\{1+it\Lambda_1(x;m)\Bigl\}\Bigl\{1+it\Lambda_2(x;m)\Bigl\}\Bigl\{1+it\Lambda_4(x;m)\Bigl\}\Bigl\{1+it\Lambda_5(x;m)\Bigl\}\Bigl| \nonumber\\
    &\qquad+O(t^2\{\mathbb{E}\Lambda_1(x;m)^2+\mathbb{E}\Lambda_2(x;m)^2+\mathbb{E}\Lambda_4(x;m)^2+\mathbb{E}\Lambda_5(x;m)^2\})+O(|t||\mathbb{E}\Lambda_3(x)|) \nonumber\\
    &\le \left|\mathbb{E}e^{it(S(x)+(\Lambda_1(x)-\Lambda_1(x;m))+(\Lambda_2(x)-\Lambda_2(x;m))+(\Lambda_4(x)-\Lambda_4(x;m))+(\Lambda_5(x)-\Lambda_5(x;m)))}\right|\nonumber\\
    &\qquad+|t|\Bigl|\mathbb{E}e^{it(S(x)+(\Lambda_1(x)-\Lambda_1(x;m))+(\Lambda_2(x)-\Lambda_2(x;m))+(\Lambda_4(x)-\Lambda_4(x;m))+(\Lambda_5(x)-\Lambda_5(x;m)))} \nonumber\\
    &\qquad\qquad\times\{\Lambda_1(x;m)+\Lambda_2(x;m)+\Lambda_4(x;m)+\Lambda_5(x;m)\}\Bigl|\nonumber\\
    &\qquad+O(t^2\{\mathbb{E}\Lambda_1(x;m)^2+\mathbb{E}\Lambda_2(x;m)^2+\mathbb{E}\Lambda_4(x;m)^2+\mathbb{E}\Lambda_5(x;m)^2\nonumber\\
    &\qquad\qquad+\mathbb{E}|\Lambda_1(x;m)\Lambda_2(x;m)|+\mathbb{E}|\Lambda_1(x;m)\Lambda_4(x;m)|+\mathbb{E}|\Lambda_1(x;m)\Lambda_5(x;m)| \nonumber\\
    &\qquad\qquad+\mathbb{E}|\Lambda_2(x;m)\Lambda_4(x;m)|+\mathbb{E}|\Lambda_2(x;m)\Lambda_5(x;m)|+\mathbb{E}|\Lambda_4(x;m)\Lambda_5(x;m)|\}) \nonumber\\
    &\qquad+O(|t||\mathbb{E}\Lambda_3(x)|). \label{Absolute ch.f PI}
\end{align}
The first term of (\ref{Absolute ch.f PI}) is bounded as below.
\begin{align}
    &\left|\mathbb{E}e^{itS(x;m)}\mathbb{E}e^{it((S(x)-S(x;m))+(\Lambda_1(x)-\Lambda_1(x;m))+(\Lambda_2(x)-\Lambda_2(x;m)+(\Lambda_4(x)-\Lambda_4(x;m)+(\Lambda_5(x)-\Lambda_5(x;m)))}\right| \nonumber\\
    &=\left|\mathbb{E}e^{itS(x;m)}\right|\left|\mathbb{E}e^{it((S(x)-S(x;m))+(\Lambda_1(x)-\Lambda_1(x;m))+(\Lambda_2(x)-\Lambda_2(x;m)+(\Lambda_4(x)-\Lambda_4(x;m))+(\Lambda_5(x)-\Lambda_5(x;m)))}\right| \nonumber\\
    &\le \left|\mathbb{E}e^{itS(x;m)}\right|=|\gamma(t)|^m. \label{Bound first}
\end{align}
Similarly, the second term of (\ref{Absolute ch.f PI}) devided by $|t|$ is bounded by
\begin{align}
    &|\mathbb{E}\{e^{itS(x;m)}\Lambda_1(x;m)\}|+|\mathbb{E}\{e^{itS(x;m)}\Lambda_2(x;m)\}| +|\mathbb{E}\{e^{itS(x;m)}\Lambda_4(x;m)\}|+|\mathbb{E}\{e^{itS(x;m)}\Lambda_5(x;m)\}|, \nonumber
\end{align}
where each term is bounded as follows. Let $C(x)$ be some positive and bounded generic function.
\begin{align}
    |\mathbb{E}\{e^{itS(x;m)}\Lambda_1(x;m)\}|&=\left|\gamma(t)^{m-1}\frac{C_{PI}C_{\Gamma,0}(x)h_0^{\frac{2L+1}{2}}}{n^{1/2}\mu_{20}^{1/2}}m\mathbb{E}\left[e^{\frac{it}{\sqrt{nh_0}}S_1(x)}\mathcal{L}_1\right]\right|+s.o. \nonumber\\
    &\le |\gamma(t)|^{m-1}\frac{mh_0^{\frac{2L+1}{2}}}{n^{1/2}\mu_{20}^{1/2}}|C_{PI}C_{\Gamma,0}(x)|\left|\mathbb{E}\left[e^{\frac{it}{\sqrt{nh_0}}S_1(x)}\mathcal{L}_1\right]\right|+s.o. \nonumber\\
    &\le |\gamma(t)|^{m-1}\frac{mh_0^{\frac{2L+1}{2}}}{n^{1/2}\mu_{20}^{1/2}}|C_{PI}C_{\Gamma,0}(x)|\mathbb{E}\left|e^{\frac{it}{\sqrt{nh_0}}S_1(x)}\mathcal{L}_1\right|+s.o. \nonumber\\
    &\le |\gamma(t)|^{m-1}\frac{mh_0^{\frac{2L+1}{2}}}{n^{1/2}\mu_{20}^{1/2}}|C_{PI}C_{\Gamma,0}(x)|\mathbb{E}\left|e^{\frac{it}{\sqrt{nh_0}}S_1(x)}\right|E\left|\mathcal{L}_1\right|+s.o. \nonumber\\
    &\le |\gamma(t)|^{m-1}\frac{mh_0^{\frac{2L+1}{2}}}{n^{1/2}\mu_{20}^{1/2}}|C_{PI}C_{\Gamma,0}(x)|\mathbb{E}\left|\mathcal{L}_1\right|+s.o. \nonumber\\
    &\le C(x)|\gamma(t)|^{m-1}\frac{mh_0^{\frac{2L+1}{2}}}{n^{1/2}}, \label{Absolute Moment Lambda 1 m}
\end{align}
where the final inequality uses Lemma \ref{L:Gamma and L exponentiation}.
\begin{align}
    |\mathbb{E}\{e^{itS(x;m)}\Lambda_2(x;m)\}|&=\left|\gamma(t)^{m-2}\frac{C_{PI}}{n^{3/2}h_0^{1/2}\mu_{20}^{1/2}}m(m-1)\mathbb{E}\left[e^{\frac{it}{\sqrt{nh_0}}(S_1+S_2)}\Gamma_1\mathcal{L}_2\right]\right| \nonumber\\
    &\le|\gamma(t)|^{m-2}\frac{m(m-1)}{n^{3/2}h_0^{1/2}\mu_{20}^{1/2}}|C_{PI}|\left|\mathbb{E}\left[e^{\frac{it}{\sqrt{nh_0}}(S_1+S_2)}\Gamma_1\mathcal{L}_2\right]\right| \nonumber\\
    &\le|\gamma(t)|^{m-2}\frac{m(m-1)}{n^{3/2}h_0^{1/2}\mu_{20}^{1/2}}|C_{PI}|\mathbb{E}\left|e^{\frac{it}{\sqrt{nh_0}}(S_1+S_2)}\right|\mathbb{E}\left|\Gamma_1\mathcal{L}_2\right| \nonumber\\
    &\le|\gamma(t)|^{m-2}\frac{m(m-1)}{n^{3/2}h_0^{1/2}\mu_{20}^{1/2}}|C_{PI}|\mathbb{E}\left|\Gamma_1\mathcal{L}_2\right| \nonumber\\
    &\le C(x)|\gamma(t)|^{m-2}\frac{m(m-1)h_0^{1/2}}{n^{3/2}}, \label{Absolute Moment Lambda 2 m}
\end{align}
where the final inequality uses Lemma \ref{L:Gamma and L exponentiation}.
\begin{align}
    |\mathbb{E}\{e^{itS(x;m)}\Lambda_4(x;m)\}|&= \left|\frac{C_{PI}m(m-1)}{2n^{3/2}h_0^{1/2}}\gamma(t)^{m-2}\mathbb{E}\left[e^{\frac{it}{\sqrt{nh_0}}(S_1+S_2)}S_1\mathcal{L}_2\right]\right| \nonumber\\
    &\le |\gamma(t)|^{m-2}\frac{m(m-1)}{2n^{3/2}h_0^{1/2}}|C_{PI}|\left|\mathbb{E}\left[e^{\frac{it}{\sqrt{nh_0}}(S_1+S_2)}S_1\mathcal{L}_2\right]\right|\nonumber\\
    &\le |\gamma(t)|^{m-2}\frac{m(m-1)}{2n^{3/2}h_0^{1/2}}|C_{PI}|\mathbb{E}\left|e^{\frac{it}{\sqrt{nh_0}}(S_1+S_2)}\right|\mathbb{E}\left|S_1\mathcal{L}_2\right|\nonumber\\
    &\le |\gamma(t)|^{m-2}\frac{m(m-1)}{2n^{3/2}h_0^{1/2}}|C_{PI}|\mathbb{E}\left|S_1\mathcal{L}_2\right|\nonumber\\
    &\le C(x)|\gamma(t)|^{m-2}\frac{m(m-1)h_0^{1/2}}{2n^{3/2}}, \label{Absolute Moment Lambda 4 m}
\end{align}
where the final inequality uses Lemma \ref{S^k L} and \ref{L:Gamma and L exponentiation}.
\begin{align}
    |\mathbb{E}\{e^{itS(x;m)}\Lambda_5(x;m)\}|&=\left|\frac{C_{PI}m}{2n^{3/2}h_0^{1/2}}\gamma(t)^{m-1}\mathbb{E}\left[e^{\frac{it}{\sqrt{nh_0}}S_1}S_1\mathcal{L}_1\right]\right| \nonumber\\
    &\le |\gamma(t)|^{m-1}\frac{C_{PI}m}{2n^{3/2}h_0^{1/2}}|C_{PI}|\left|\mathbb{E}\left[e^{\frac{it}{\sqrt{nh_0}}S_1}S_1\mathcal{L}_1\right]\right| \nonumber\\
    &\le |\gamma(t)|^{m-1}\frac{C_{PI}m}{2n^{3/2}h_0^{1/2}}|C_{PI}|\mathbb{E}\left|e^{\frac{it}{\sqrt{nh_0}}S_1}\right|\mathbb{E}\left|S_1\mathcal{L}_1\right|\nonumber\\
    &\le |\gamma(t)|^{m-1}\frac{C_{PI}m}{2n^{3/2}h_0^{1/2}}|C_{PI}|\mathbb{E}\left|S_1\mathcal{L}_1\right|\nonumber\\
    &\le C(x)|\gamma(t)|^{m-1}\frac{mh_0^{1/2}}{2n^{3/2}}\label{Absolute Moment Lambda 5 m}
\end{align}
where the final inequality uses Lemma \ref{S^k L}.
 (\ref{Absolute Moment Lambda 1 m}), (\ref{Absolute Moment Lambda 2 m}), (\ref{Absolute Moment Lambda 4 m}) and (\ref{Absolute Moment Lambda 5 m}) imply
\begin{align}
    &|t|\Bigl|\mathbb{E}e^{it(S(x)+(\Lambda_1(x)-\Lambda_1(x;m))+(\Lambda_2(x)-\Lambda_2(x;m))+(\Lambda_4(x)-\Lambda_4(x;m))+(\Lambda_5(x)-\Lambda_5(x;m))}\nonumber\\
    &\qquad\times\{\Lambda_1(x;m)+\Lambda_2(x;m)+\Lambda_4(x;m)+\Lambda_5(x;m)\}\Bigl| \nonumber\\
    &\le C(x)\Biggl\{|\gamma(t)|^{m-1}\frac{mh_0^{\frac{2L+1}{2}}}{n^{1/2}}+|\gamma(t)|^{m-2}\frac{m^2h_0^{1/2}}{n^{3/2}}+|\gamma(t)|^{m-1}\frac{mh_0^{1/2}}{n^{3/2}}\Biggl\}|t|. \label{Bound second}
\end{align}
Then, (\ref{Absolute ch.f PI}), (\ref{Bound first}), and (\ref{Bound second}) yield 
\begin{align}
    |\chi_{PI}(t)|&\le |\gamma(t)|^m+C(x)\Biggl\{|\gamma(t)|^{m-1}\frac{mh_0^{\frac{2L+1}{2}}}{n^{1/2}}+|\gamma(t)|^{m-2}\frac{m^2h_0^{1/2}}{n^{3/2}}+|\gamma(t)|^{m-1}\frac{mh_0^{1/2}}{n^{3/2}}\Biggl\}|t| \nonumber\\
    &\qquad+O(t^2\{\mathbb{E}\Lambda_1(x;m)^2+\mathbb{E}\Lambda_2(x;m)^2+\mathbb{E}\Lambda_4(x;m)^2+\mathbb{E}\Lambda_5(x;m)^2\nonumber\\
    &\qquad\qquad+\mathbb{E}|\Lambda_1(x;m)\Lambda_2(x;m)|+\mathbb{E}|\Lambda_1(x;m)\Lambda_4(x;m)|+\mathbb{E}|\Lambda_1(x;m)\Lambda_5(x;m)| \nonumber\\
    &\qquad\qquad+\mathbb{E}|\Lambda_2(x;m)\Lambda_4(x;m)|+\mathbb{E}|\Lambda_2(x;m)\Lambda_5(x;m)|+\mathbb{E}|\Lambda_4(x;m)\Lambda_5(x;m)|\}) \nonumber\\
    &\qquad+O(|t||\mathbb{E}\Lambda_3(x)|) \nonumber\\
    &\le C(x)|\gamma(t)|^{m-2}\Biggl[1+\Biggl\{\frac{mh_0^{\frac{2L+1}{2}}}{n^{1/2}}+\frac{m^2h_0^{1/2}}{n^{3/2}}+\frac{mh_0^{1/2}}{n^{3/2}}\Biggl\}|t|\Biggl] \nonumber\\
    &\qquad+O\left(t^2\Biggl\{\frac{mh_0^{2L+1}}{n}+\frac{m^2}{n^3}+\frac{m^{3/2}h_0^{(2L+1)/2}}{n^2}\Biggl\}\right)+O\left(|t|\left\{n^{-1/2}h_0^{(2L+1)/2}+n^{-1}\right\}\right) \nonumber
\end{align}
where second inequality uses $|\gamma(t)|\le1$ and Lemma \ref{L:AM Lambda 3},\ref{L:lambda AM 1m},\ref{L:lambda AM 2m},\ref{L:lambda AM 4m},\ref{L:lambda AM 5m} and \ref{L:lambda AM m product}.

We evaluate $(B)$, partioning its range of integration into two parts, $p\le|t|\le \frac{n^{1/2}h_0^{1/2}}{\mu_{20}^{-3/2}\mu_{30}}$ and $\frac{n^{1/2}h_0^{1/2}}{\mu_{20}^{-3/2}\mu_{30}}\le|t|\le n^{\frac{2L}{2L+1}}\log n$.

$(i)$ For $p\le|t|\le \frac{n^{1/2}h_0^{1/2}}{\mu_{20}^{-3/2}\mu_{30}}$\\
Applying Taylor expansion to $e^{\frac{it}{\sqrt{nh_0}}S_1(x)}$ with respect to $t$, we have
\begin{align}
    |\gamma(t)-1-\frac{t^2}{2n}|\le \frac{|t|^3\mu_{20}^{-3/2}\mu_{30}}{6n^{3/2}h_0^{1/2}}, \nonumber
\end{align}
then for $|t|\le \frac{n^{1/2}h_0^{1/2}}{\mu_{20}^{-3/2}\mu_{30}}$,
\begin{align}
    |\gamma(t)|&\le 1-\frac{t^2}{2n}+\frac{|t|^3\mu_{20}^{-3/2}\mu_{30}}{6n^{3/2}h_0^{1/2}}\le 1-\frac{t^2}{2n}+\frac{t^2}{6n}=1-\frac{t^2}{3n}\le \exp\left(-\frac{t^2}{3n}\right), \nonumber
\end{align}
then
\begin{align}
    |\chi_{PI}(t)|&\le C(x)|\gamma(t)|^{m-2}\Biggl[1+\Biggl\{\frac{mh_0^{\frac{2L+1}{2}}}{n^{1/2}}+\frac{m^2h_0^{1/2}}{n^{3/2}}+\frac{mh_0^{1/2}}{n^{3/2}}\Biggl\}|t|\Biggl] \nonumber\\
    &\qquad+O\left(t^2\Biggl\{\frac{mh_0^{2L+1}}{n}+\frac{m^2}{n^3}+\frac{m^{3/2}h_0^{(2L+1)/2}}{n^2}\Biggl\}\right)+O\left(|t|\left\{n^{-1/2}h_0^{(2L+1)/2}+n^{-1}\right\}\right) \nonumber\\
    &\le C(x)\exp\left(-\frac{(m-2)t^2}{3n}\right)\left[1+\Biggl\{\frac{mh_0^{\frac{2L+1}{2}}}{n^{1/2}}+\frac{m^2h_0^{1/2}}{n^{3/2}}+\frac{mh_0^{1/2}}{n^{3/2}}\Biggl\}|t|\right] \nonumber\\
    &\qquad+O\left(t^2\Biggl\{\frac{m}{n^2}+\frac{m^2}{n^3}+\frac{m^{3/2}}{n^{5/2}}\Biggl\}\right)+O\left(|t|n^{-1}\}\right). \nonumber
\end{align}
Using (A.21) in \cite{NR00}, we can take $m=[9n\log n/t^2]$ since $1\le m \le n-1$ holds for $p\le |t| \le \frac{n^{1/2}h_0^{1/2}}{\mu_{20}^{-3/2}\mu_{30}}$ and sufficiently large $n$.\\
Because $m\ge (9n\log n)/t^2-1$, for $|t| \le \frac{n^{1/2}h_0^{1/2}}{\mu_{20}^{-3/2}\mu_{30}}$
\begin{align}
    \exp\left(-\frac{(m-2)t^2}{3n}\right)=\exp\left(-\frac{(m+1)t^2}{3n}\right)\exp\left(\frac{3t^2}{3n}\right)\le C\exp(-3\log n)\le \frac{C}{n^3}, \nonumber
\end{align}
and this implies, using $m\le (9n\log n)/t^2$,
\begin{align}
    |\chi_{PI}(t)|
    &\le\frac{C(x)}{n^3}\Biggl[1+n^{1/2}(\log n)h_0^{\frac{2L+1}{2}}\frac{1}{|t|}+n^{1/2}(\log n)^2h_0^{1/2}\frac{1}{|t|^3}+n^{-1/2}(\log n)h_0^{1/2}\frac{1}{|t|}\Biggl] \nonumber\\
    &\qquad+ O\left(n^{-1}(\log n)+n^{-1}(\log n)^2\frac{1}{t^2}+n^{-1}(\log n)^{3/2}\frac{1}{|t|}\right)+O(|t|n^{-1})\nonumber
\end{align}
Therefore, dropping the integral range $p\le |t|\le \frac{n^{1/2}h_0^{1/2}}{\mu_{20}^{-3/2}\mu_{30}}$ on the right-hand side,
\begin{align}
    &\int_{p\le |t|\le \frac{n^{1/2}h_0^{1/2}}{\mu_{20}^{-3/2}\mu_{30}} }\left|\frac{\chi_{PI}(t)}{t}\right|dt \nonumber\\
    &\qquad\le C(x)\Biggl[\Bigl\{n^{-3}+n^{-1}(\log n)\Bigl\}\int \frac{dt}{|t|}+\Bigl\{n^{-5/2}(\log n)h_0^{\frac{2L+1}{2}}+n^{-7/2}(\log n)h_0^{1/2} \nonumber\\
    &\qquad\qquad+n^{-1}(\log n)^{3/2}\Bigl\}\int \frac{dt}{t^2}+n^{-1}(\log n)^2\int \frac{dt}{|t|^3} +n^{-5/2}(\log n)^2h_0^{1/2}\int \frac{dt}{t^4}\Biggl]+O(n^{-1}) \nonumber\\
    &\qquad=o\{(nh_0)^{-1}\} \nonumber
\end{align}

$(ii)$ For $\frac{n^{1/2}h_0^{1/2}}{\mu_{20}^{-3/2}\mu_{30}}\le |t|\le n^{\frac{2L}{2L+1}}\log n$, there exist $\eta\in(0,1)$, such that $|\gamma(t)|\le 1-\eta$ from Assumption \ref{A:Cramer}. We can take $m=[-3\log n/\log (1-\eta)]$ since $1\le m\le n-1$ for sufficiently large $n$. Then $\chi_{PI}(t)$ is bounded as follow.
\begin{align}
    &|\chi_{PI}(t)| \nonumber\\
    &\le C(1-\eta)^{-3\log n/\log(1-\eta)}\nonumber\\
    &\qquad\times\left[1+\left\{\frac{h_0^{\frac{2L+1}{2}}}{n^{1/2}}+\frac{h_0^{1/2}}{n^{3/2}}\right\}|t|\left(\frac{-3\log n}{\log(1-\eta)}\right)+\frac{h_0^{1/2}}{n^{3/2}}|t|\left(\frac{-3\log n}{\log(1-\eta)}\right)^2\right] \nonumber\\
    &\qquad+O\left(t^2\left\{n^{-2}\left(\frac{-3\log n}{\log(1-\eta)}\right)+n^{-3}\left(\frac{-3\log n}{\log(1-\eta)}\right)^2+n^{-5/2}\left(\frac{-3\log n}{\log(1-\eta)}\right)^{3/2}\right\}\right) \nonumber
\end{align}
Noting that
\begin{align}
    (1-\eta)^{-3\log n/\log(1-\eta)}=(1-\eta)^{\log n^{-3}/\log(1-\eta)}=(1-\eta)^{\log_{(1-\eta)}n^{-3}}=n^{-3}, \nonumber
\end{align}
\begin{align}
    &\int_{\frac{n^{1/2}h_0^{1/2}}{\mu_{20}^{-3/2}\mu_{30}}\le |t|\le n^{\frac{2L}{2L+1}}\log n}\left|\frac{\chi_{PI}(t)}{t}\right|dt \nonumber\\
    &=O\Biggl(\frac{\log(n^{\frac{2L}{2L+1}}\log n)}{n^3}+\frac{\log n}{n^3}\left\{\frac{h_0^{\frac{2L+1}{2}}}{n^{1/2}}+\frac{h_0^{1/2}}{n^{3/2}}\right\}(n^{\frac{2L}{2L+1}}\log n)+\frac{(\log n)^2}{n^3}\frac{h_0^{1/2}}{n^{3/2}}(n^{\frac{2L}{2L+1}}\log n)\Biggl) \nonumber\\
    &\qquad+O\left(n^{\frac{2L}{2L+1}}\log n\left\{n^{-2}(\log n)+n^{-2}(\log n)^2 + n^{-3/2}(\log n)^{3/2}\right\}\right) \nonumber\\
    &=o\{(nh_0)^{-1}\} \nonumber
\end{align}

Finally, we evaluate $(C)$. For some constant $C$,
\begin{align}
    (C)&=\int_{p\le|t|}\frac{1}{|t|}e^{\frac{-t^2}{2}}\Biggl|1+\frac{\mu_{30}\mu_{20}^{-3/2}}{6n^{1/2}h_0^{1/2}}(it)^3+\frac{\mu_{40}\mu_{20}^{-2}}{24nh_0}(it)^4+\frac{\mu_{30}^2\mu_{20}^{-3}}{72nh_0}(it)^6 \nonumber\\
    &\qquad+C_{PI}\rho_{11}\mu_{20}^{-1}\left(\sum_{l=0}^{L-1}C_{\Gamma,l}(x)h_0^{L+l+1}\right)(it)^2+C_{PI}\frac{\rho_{11}\xi_{11}\mu_{20}^{-3/2}h_0^{1/2}}{n^{1/2}}(it)^3 -C_{PI}\frac{\mu_{20}^{-1/2}\rho_{11}h_0^{1/2}}{2n^{1/2}}\{(it)^3+(it)\}\Biggl|dt \nonumber\\
    &\le C\Biggl[\int_{p}^{\infty}\frac{1}{t}e^{\frac{-t^2}{2}}dt+\frac{h_0^{1/2}}{n^{1/2}}\int_{p}^{\infty}t^2e^{\frac{-t^2}{2}}dt+\frac{1}{nh_0}\int_{p}^{\infty}(t^3+t^5)e^{\frac{-t^2}{2}}dt \nonumber\\
    &\qquad+h_0^{L+1}\int_{p}^{\infty}te^{\frac{-t^2}{2}}dt+\frac{h_0^{1/2}}{n^{1/2}}\int_{p}^{\infty}t^2e^{\frac{-t^2}{2}}dt+\frac{h_0^{1/2}}{n^{1/2}}\int^{\infty}_p(t^2+1)e^{\frac{-t^2}{2}}dt\Biggl] \nonumber
\end{align}
Since $p=\min\left\{\frac{n^{1/2}h_0^{1/2}}{\mu_{20}^{-3/2}\mu_{30}},\log n\right\}$,the first integral is smaller than $p^{-2}\int_{p}^{\infty}te^{-t^2/2}dt=p^{-2}e^{-p^2/2}=o(n^{-1})$, 
the second and fifth integrals are smaller than
$p^{-1}\int_{p}^{\infty}t^3e^{-t^2/2}dt=p^{-1}e^{-p^2/2}(p^2+2)=o(n^{-1})$, 
the third integral is $\int_{p}^{\infty}(t^3+t^5)e^{-t^2/2}dt=e^{-p^2/2}(p^4+5p^2+10)=o(n^{-1})$, 
the fourth integral is $\int_{p}^{\infty}te^{-t^2/2}dt=e^{-p^2/2}=o(n^{-1})$, and
the final integral is $p^{-1}e^{-t^2/2}(p^2+3)=o(n^{-1})$. 
It follows that $(C)=o\{(nh_0)^{-1}\}$. Thus the expansion is valid.
\end{proof}

\subsection{Proof of Theorem \ref{Coro: Exact Order}} \label{Proof of Exact Order}
\begin{proof}
Define $\epsilon$ and $\epsilon_{PI}$ as follows,
\begin{align}
    & \epsilon \equiv \mathbb{P}(S(x)\le z)-\Phi(z)-\phi(z)\Biggl[(nh_0)^{-1/2}p_1(z)+(nh_0)^{-1}p_2(z)\Biggl], \nonumber\\
    &\epsilon_{PI} \equiv \mathbb{P}(S_{\text{PI}}(x)\le z)-\Phi(z) \nonumber\\&\qquad-\phi(z)\Biggl[(nh_0)^{-1/2}p_1(z)+h_0^{L+1}p_{3,0}(z)+n^{-1/2}h_0^{1/2}p_4(z)+(nh_0)^{-1}p_2(z)\Biggl]. \nonumber
\end{align}
Then we have,
\begin{align}
    &\sup_{z\in\mathbb{R}}\left|\mathbb{P}(S(x)\le z)-\mathbb{P}(S_{PI}(x)\le z)-\phi(z)\Bigl[h_0^{L+1}p_{3,0}(z)+n^{-1/2}h_0^{1/2}p_4(z)\Bigl]\right| \nonumber\\
    &=\sup_{z\in\mathbb{R}}\left|\epsilon-\epsilon_{PI}\right|=o(h_0^{L+1}+n^{-1/2}h_0^{1/2}). \nonumber
\end{align}
\end{proof}

\section{Lemmas}

\begin{lemma}\label{L:EGamma1}
Under Assumptions \ref{A:DGP}, \ref{A: interior point}, \ref{A: nh infty}, \ref{A:2L+L_p times differentiable}, \ref{A: Kernel ''} and \ref{A:Ku}
\begin{equation}
    \mathbb{E}\Gamma_{KDE_1}=\sum_{l=0}^{L-1}C_{\Gamma,l}(x)h_0^{L+l}+o(h_0^{2L-1}),\quad \text{where}\quad C_{\Gamma,l}(x)\equiv -\left(\int u^{L+l}K(u)du\right)\frac{f^{(L+l)}(x)}{(L+l-1)!} \nonumber
\end{equation}
\begin{proof}
\begin{align*}
    \mathbb{E}\Gamma_{KDE_1}&=\mathbb{E}\left[\frac{1}{nh_0}\sum_{i=1}^nK'\left(\frac{X_i-x}{h_0}\right)\left(\frac{X_i-x}{h_0}\right)\right]+\mathbb{E}\left[\frac{1}{nh_0}\sum_{i=1}^nK\left(\frac{X_i-x}{h_0}\right)\right]\\&=\frac{1}{h_0}\int K'\left(\frac{z-x}{h_0}\right)\left(\frac{z-x}{h_0}\right)f(z)dz+\frac{1}{h_0}\int K\left(\frac{z-x}{h_0}\right)f(z)dz\\
    &=\int K'(u)uf(x+uh_0)du+\int K(u)f(x+uh_0)du\\
    &=-\int K(u)f(x+uh_0)du-\int K(u)uf'(x+uh_0)h_0du + \int K(u)f(x+uh_0)du\\
    &=-\int K(u)uf'(x+uh_0)h_0du\\
    &=-\int K(u)u\Bigl\{f^{(1)}(x)+\cdots+\frac{f^{(L)}(x)}{(L-1)!}(uh_0)^{L-1}+\cdots+\frac{f^{(2L)}}{(2L-1)!}(uh_0)^{2L-1}\Bigl\}h_0du+o(h_0^{2L-1})\\
    &=-\sum_{l=0}^{L}\left(\int u^{L+l}K(u)du\right)\frac{f^{(L+l)}(x)}{(L+l-1)!}h_0^{L+l}+o(h_0^{2L-1})\\
    &\equiv \sum_{l=0}^{L-1}C_{\Gamma,l}(x)h_0^{L+l}+o(h_0^{2L-1})
\end{align*}
The fourth equality follows from integration by part of the first term and Assumption \ref{A: L times differentiable},\ref{A:Ku}, the seventh equality follows from the expansion of $f'(x+uh_0)$ around $h_0=0$ and Assumption \ref{A: L times differentiable} and the eighth equality follows from Assumption \ref{A: Kernel ''}. 
\end{proof}
\end{lemma}

\begin{lemma}\label{L:EGamma2}
Under Assumptions \ref{A:DGP}, \ref{A: L times differentiable}, \ref{A: interior point},\ref{A: nh infty}, \ref{A: Kernel ''}, and \ref{A:Ku^2},
\begin{equation}
    \mathbb{E}\Gamma_{KDE_2}=O(h_0^L) \nonumber
\end{equation}
\begin{proof}
\begin{align}
    \mathbb{E}\Gamma_{KDE_2}&=\mathbb{E}\left[\frac{1}{nh_0}\sum_{i=1}^nK''\left(\frac{X_i-x}{h_0}\right)\left(\frac{X_i-x}{h_0}\right)^2\right] +\mathbb{E}\left[\frac{4}{nh_0}\sum_{i=1}^nK'\left(\frac{X_i-x}{h_0}\right)\left(\frac{X_i-x}{h_0}\right)]+\mathbb{E}[\frac{2}{nh_0}\sum_{i=1}^nK\left(\frac{X_i-x}{h_0}\right)\right] \nonumber\\
    &=\frac{1}{h_0}\int K''\left(\frac{z-x}{h_0}\right)\left(\frac{z-x}{h_0}\right)^2f(z)dz +\frac{4}{h_0}\int K'\left(\frac{z-x}{h_0}\right)\left(\frac{z-x}{h_0}\right)f(z)dz+\frac{2}{h_0}\int K\left(\frac{z-x}{h_0}\right)f(z)dz \nonumber\\
    &=2\int K(u)f(x+uh_0)du \nonumber\\
    &\qquad+4\left\{-\int K(u)f(x+uh_0)du-\int K(u)uf'(x+uh_0)h_0du\right\} \nonumber\\
    &\qquad+\left\{-2\int K'(u)uf(x+uh_0)-\int K'(u)u^2f'(x+uh_0)h_0du\right\} \nonumber\\
    &=2\int K(u)f(x+uh_0)du \nonumber\\
    &\qquad+4\left\{-\int K(u)f(x+uh_0)du-\int K(u)uf'(x+uh_0)h_0du\right\} \nonumber\\
    &\qquad+\left\{-2\left[-\int K(u)f(x+uh_0)du-\int K(u)uf'(x+uh_0)h_0du\right]\right\} \nonumber\\
    &\qquad+\left\{-\left[-2\int K(u)uf'(x+uh_0)h_0du-\int K(u)u^2f''(x+uh_0)h_0^2du\right]\right\} \nonumber\\
    &=\int K(u)u^2f''(x+uh_0)h_0^2du=O(h_0^L) \nonumber
\end{align}
The third equality follows from integration by parts of the first and second terms and Assumption \ref{A:Ku^2}, the fourth equality follows from integration by parts of the second term and Assumption \ref{A:Ku^2} and the final equality follows from the expansion of $f''(x+uh_0)$ around $h_0=0$ and Assumptions \ref{A: L times differentiable}, \ref{A: Kernel ''}.
\end{proof}
\end{lemma}

\begin{lemma}\label{L:CS inequality 2}
Under Assumptions \ref{A:DGP}, \ref{A: L times differentiable}, \ref{A: interior point},\ref{A: nh infty}, \ref{A: Kernel ''}, \ref{A:L_p order 1} and \ref{A:Ku^2},
\begin{align}
    \mathbb{E}\left|\sqrt{nh_0}\left(\frac{\hat{h}-h_0}{h_0}\right)^2\Gamma_{KDE_2}\right|=o\{(nh_0)^{-1}\} \nonumber
\end{align}
\begin{proof}
Similar to the proof of Theorem 1.
\end{proof}
\end{lemma}

\begin{lemma} \label{gamma L}
Under Assumptions \ref{A:DGP}, \ref{A: L times differentiable}, \ref{A: interior point}, \ref{A:2L+L_p times differentiable},\ref{A: Kernel ''} and \ref{A:Ku^2},
\begin{equation}
    \mathbb{E}[\Gamma_1\mathcal{L}_1]=O(h_0^{L+1}) \nonumber
\end{equation}
\begin{proof}
 Letting $g(x)=f^{(L)}(x)f(x)$
\begin{align}
    \mathbb{E}[\Gamma_1\mathcal{L}_1]
    &=\mathbb{E}\Biggl[\Biggl\{K'\left(\frac{X_1-x}{h_0}\right)\left(\frac{X_1-x}{h_0}\right)+K\left(\frac{X_1-x}{h_0}\right) -\mathbb{E}\left[K'\left(\frac{X_1-x}{h_0}\right)\left(\frac{X_1-x}{h_0}\right)+K\left(\frac{X_1-x}{h_0}\right)\right]\Biggl\}\mathcal{L}_1\Biggl] \nonumber\\
    &=\mathbb{E}\left[\left\{K'\left(\frac{X_1-x}{h_0}\right)\left(\frac{X_1-x}{h_0}\right)+K\left(\frac{X_1-x}{h_0}\right)\right\}f^{(L)}(X_1)\right] -\mathbb{E}\left[\left\{K'\left(\frac{X_1-x}{h_0}\right)\left(\frac{X_1-x}{h_0}\right)+K\left(\frac{X_1-x}{h_0}\right)\right\}\right]\mathbb{E}[f^{(L)}(X_1)] \nonumber
\end{align}
We can compute the first term as follow. 
\begin{align}
    &\mathbb{E}\left[\left\{K'\left(\frac{X_1-x}{h_0}\right)\left(\frac{X_1-x}{h_0}\right)+K\left(\frac{X_1-x}{h_0}\right)\right\}f^{(L)}(X_1)\right]  \nonumber \\&= \int \left\{K'\left(\frac{z-x}{h_0}\right)\left(\frac{z-x}{h_0}\right)+K\left(\frac{z-x}{h_0}\right)\right\}f^{(L)}(z)f(z)dz \nonumber\\
    &= h_0\int \Bigl\{K'(u)u+K(u)\Bigl\}g(x+uh_0)du \nonumber\\
    &= h_0\int \Bigl\{K'(u)u+K(u)\Bigl\}\left\{g(x)+\cdots+\frac{g^{(L)}(x)}{L!}(uh_0)^L+o(h_0^L)\right\}du \nonumber\\
    &= -h_0\int K(u)l(x)du + h_0\int K(u)l(x)du+h_0\int K'(u)u\frac{g^{(L)}(x)}{L!}(uh_0)^Ldu + h_0\int K(u)\frac{g^{(L)}(x)}{L!}(uh_0)^ldu+o(h_0^{L+1})\nonumber\\
    &= \frac{g^{(L)}(x)}{L!}h_0^{L+1}\int K'(u)u^{L+1}du +\frac{g^{(L)}(x)}{L!}h_0^{L+1}\int K(u)u^Ldu+o(h_0^{L+1}) \nonumber\\
    &= -(L+1)\frac{g^{(L)}(x)}{L!}h_0^{L+1}\int K(u)u^Ldu + \frac{g^{(L)}(x)}{L!}h_0^{L+1}\int K(u)u^Ldu+o(h_0^{L+1}) \nonumber\\
    &=O(h_0^{L+1}) \nonumber
\end{align}
The fourth equality follows from the expansion of $l(x+uh_0)$ around $h_0=0$ and Assumption \ref{A:2L+L_p times differentiable}, and the fifth equality follows from integration by parts of the products of $K'(u)u$ and $l^{(k)}(x)u^k,(0\le k\le L-1)$ and Assumption \ref{A:2L+L_p times differentiable} and \ref{A:Ku^2}. Next, we can compute the second term similarly to the first term.
\begin{align}
    &\mathbb{E}[\left\{K'\left(\frac{X_1-x}{h_0}\right)\left(\frac{X_1-x}{h_0}\right)+K\left(\frac{X_1-x}{h_0}\right)\right\}]\mathbb{E}[f^{(L)}(X_1)] \nonumber\\
    &=\left(-(L+1)\frac{f^{(L)}(x)}{L!}h_0^{L+1}\int K(u)u^Ldu + \frac{f^{(L)}(x)}{L!}h_0^{L+1}\int K(u)u^Ldu\right)\mathbb{E}[f^{(L)}(X_1)] \nonumber\\
    &=O(h_0^{L+1}) \nonumber
\end{align}
These imply the lemma holds.
\end{proof}
\end{lemma}

\begin{lemma} \label{S^k L}
For any positive integer $k$ and any non-negative integer $l$, 
\begin{equation}
    \mathbb{E}|S_1^k\mathcal{L}_1^l|=O(h_0) \nonumber
\end{equation}
\begin{proof}
Straightforward.
\end{proof}
\end{lemma}

\begin{lemma} \label{S^k gamma}
For any positive integer $k,l$,
\begin{equation}
    \mathbb{E}|S_1^k\Gamma_1^l|=O(h_0) \nonumber
\end{equation}
\begin{proof}
Straightforward.
\end{proof}
\end{lemma}

\begin{lemma} \label{L:Gamma and L exponentiation}
For any positive integer $k,l\geqq 2$,
\begin{equation}
    \mathbb{E}|\Gamma_1^l|^k=O(h_0^k),~~~\mathbb{E}|\mathcal{L}_1^l|^k=O(1) \nonumber
\end{equation}
\begin{proof}
Straightforward.
\end{proof}
\end{lemma}

\begin{lemma}\label{L: AM Lambda 1}
For any positive integer r,
\begin{equation}
    \mathbb{E}|\Lambda_1(x)|^r=O(h_0^{\frac{r(2L+1)}{2}}) \nonumber
\end{equation}
\begin{proof}
From Lemma \ref{L:Gamma and L exponentiation}, for any positive integer $k$, and some positive bounded function $C(x)$,
\begin{align}
    \mathbb{E}|\Lambda_1(x)|^{2k}=\mathbb{E}\Lambda_1(x)^{2k}&\lesssim\frac{h_0^{k(2L+1)}}{n^k}\mathbb{E}\left[\left(\sum_{i=1}^n\mathcal{L}_i\right)^{2k}\right]+s.o. \nonumber\\
    &=\frac{h_0^{k(2L+1)}}{n^k}n^k\mathbb{E}\left[\mathcal{L}_i^2\right]^k+s.o.=O(h_0^{k(2L+1)}) \nonumber
\end{align}
From Holder's inequality, for $0<r<s$, $\mathbb{E}|X|^r\le \{\mathbb{E}|X|^s\}^{r/s}$, thus for any positive integer $k$,
\begin{equation}
    \mathbb{E}|\Lambda_1(x)|^{2k-1}\le \{\mathbb{E}|\Lambda_1(x)|^{2k}\}^{\frac{2k-1}{2k}}=O(h_0^{\frac{(2k-1)(2L+1)}{2}}) \nonumber
\end{equation}
This implies the lemma holds.
\end{proof}
\end{lemma}
\begin{lemma}\label{L:AM Lambda 2}
For any positive integer $r$,
\begin{equation}
    \mathbb{E}|\Lambda_2(x)|^r=O(n^{-r/2}) \nonumber
\end{equation}
\begin{proof}
From Lemma \ref{L:Gamma and L exponentiation}, for any positive integer $k$,
\begin{align}
    \mathbb{E}|\Lambda_2(x)|^{2k}&\lesssim\frac{1}{n^{3k}h_0^k\mu_{20}^k}\mathbb{E}\left[\left(\sum_{i=1}^n\sum_{j\neq i}^n\Gamma_i\mathcal{L}_j\right)^{2k}\right] \nonumber\\
    &=\frac{n^k(n-1)^k}{n^{3k}h_0^k\mu_{20}^k}\mathbb{E}\left[\Gamma_1^2\right]^{k}\mathbb{E}\left[\mathcal{L}_2^2\right]^{k}+s.o.=O(n^{-k})\nonumber
\end{align}
Then, similarly to the evaluation of $\mathbb{E}|\Lambda_1(x)|^r$, the lemma holds.
\end{proof}
\end{lemma}
\begin{lemma}\label{L:AM Lambda 12}
\begin{equation}
    \mathbb{E}|\Lambda_1(x)\Lambda_2(x)|=O(n^{-1/2}h_0^{\frac{2L+1}{2}}) \nonumber
\end{equation}
\begin{proof}
Lemma \ref{L: AM Lambda 1}, \ref{L:AM Lambda 2} and Holder inequality implies
\begin{equation}
    \mathbb{E}|\Lambda_1(x)\Lambda_2(x)|\le \mathbb{E}|\Lambda_1(x)|\mathbb{E}|\Lambda_2(x)|=O(h_0^{\frac{2L+1}{2}})O(n^{-1/2}) \nonumber
\end{equation}
\end{proof}
\end{lemma}

\begin{lemma}\label{L:AM Lambda 3}
\begin{equation}
    \mathbb{E}|\Lambda_3(x)|=O(n^{-1/2}h_0^{\frac{2L+1}{2}}+n^{-1}) \nonumber
\end{equation}
\begin{proof}
From Lemma \ref{gamma L},
\begin{align}
    \mathbb{E}\Lambda_3(x)^2&\lesssim\frac{1}{n^3h_0\mu_{20}}\mathbb{E}\left[\sum_{i=1}^n\sum_{j=1}^n\Gamma_i\mathcal{L}_i\Gamma_j\mathcal{L}_j\right] \nonumber\\
    &=\frac{1}{n^3h_0\mu_{20}}\mathbb{E}\left[\sum_{i=1}^n\sum_{j\neq i}^n\Gamma_i\mathcal{L}_i\Gamma_j\mathcal{L}_j+\sum_{i=1}^n\Gamma_i^2\mathcal{L}_i^2\right] \nonumber\\
    &=\frac{n(n-1)}{n^3h_0\mu_{20}}\mathbb{E}[\Gamma_1\mathcal{L}_1]^2+\frac{1}{n^2h_0\mu_{20}}\mathbb{E}[\Gamma_1^2\mathcal{L}_1^2] \nonumber\\
    &=O(n^{-1}h_0^{-1})O(h_0^{2(L+1)})+O(n^{-2}h_0^{-1})O(h_0) \nonumber
\end{align}
Similarly to the evaluation of $\mathbb{E}|\Lambda_1(x)|^r$, the lemma holds.
\end{proof}
\end{lemma}

\begin{lemma} \label{L:AM Lambda 4}
\begin{align}
    \mathbb{E}|\Lambda_4(x)|^2=O(n^{-1}h_0) \nonumber
\end{align}
\begin{proof}
The proof is similar to Lemma \ref{L:AM Lambda 2}.
\end{proof}
\end{lemma}

\begin{lemma} \label{L:AM Lambda 5}
\begin{align}
    \mathbb{E}|\Lambda_5(x)|^2=O(n^{-1}h_0) \nonumber
\end{align}
\begin{proof}
The proof is similar to Lemma \ref{L:AM Lambda 3}.
\end{proof}
\end{lemma}

\begin{lemma} \label{L:AM Lambda 14}
\begin{align}
    \mathbb{E}|\Lambda_1(x)\Lambda_4(x)|=O(n^{-1/2}h_0^{L+1}) \nonumber
\end{align}
\begin{proof}
From Cauchy-Schwarz inequality and  Lemma \ref{L: AM Lambda 1} and \ref{L:AM Lambda 4}, this lemma holds.
\end{proof}
\end{lemma}

\begin{lemma} \label{L:AM Lambda 15}
\begin{align}
    \mathbb{E}|\Lambda_1(x)\Lambda_5(x)|=O(n^{-1/2}h_0^{L+1}) \nonumber
\end{align}
\begin{proof}
From Cauchy-Schwarz inequality and  Lemma \ref{L: AM Lambda 1} and \ref{L:AM Lambda 5}, this lemma holds.
\end{proof}
\end{lemma}

\begin{lemma} \label{L:AM Lambda 24}
\begin{align}
    \mathbb{E}|\Lambda_2(x)\Lambda_4(x)|=O(n^{-1}h_0^{1/2}) \nonumber
\end{align}
\begin{proof}
From Cauchy-Schwarz inequality and  Lemma \ref{L:AM Lambda 2} and \ref{L:AM Lambda 4}, this lemma holds.
\end{proof}
\end{lemma}

\begin{lemma} \label{L:AM Lambda 25}
\begin{align}
    \mathbb{E}|\Lambda_2(x)\Lambda_5(x)|=O(n^{-1}h_0^{1/2}) \nonumber
\end{align}
\begin{proof}
From Cauchy-Schwarz inequality and  Lemma \ref{L:AM Lambda 2} and \ref{L:AM Lambda 5}, this lemma holds.
\end{proof}
\end{lemma}

\begin{lemma} \label{L:AM Lambda 45}
\begin{align}
    \mathbb{E}|\Lambda_4(x)\Lambda_5(x)|=O(n^{-1}h_0) \nonumber
\end{align}
\begin{proof}
From Cauchy-Schwarz inequality and  Lemma \ref{L:AM Lambda 4} and \ref{L:AM Lambda 5}, this lemma holds.
\end{proof}
\end{lemma}

\begin{lemma} \label{L:lambda AM 1m}
\begin{align}
    \mathbb{E}|\Lambda_1(x;m)|^2=O\left(\frac{mh_0^{2L+1}}{n}\right) \nonumber
\end{align}
\begin{proof}
From Lemma \ref{L:Gamma and L exponentiation},
\begin{align}
    &\mathbb{E}\left|\Lambda_1(x;m)\right|^2 \nonumber\\
    &\quad\lesssim \frac{h_0^{2L+1}}{n}\mathbb{E}\left[\sum_{i=1}^m\sum_{j=1}^m\mathcal{L}_i\mathcal{L}_j\right]= \frac{h_0^{2L+1}}{n}\mathbb{E}\left[\sum_{i=1}^m\sum_{j\neq i}^m\mathcal{L}_i\mathcal{L}_j+\sum_{i=1}^m\mathcal{L}_i^2\right]=O\left(\frac{mh_0^{2L+1}}{n}\right) \nonumber
\end{align}
\end{proof}
\end{lemma}

\begin{lemma} \label{L:lambda AM 2m}
\begin{align}
    \mathbb{E}|\Lambda_2(x;m)|^2=O\left(\frac{m^2}{n^3}\right) \nonumber
\end{align}
\begin{proof}
\begin{align}
    &\mathbb{E}|\Lambda_2(x;m)|^2 \nonumber\\
    &\quad \lesssim \frac{1}{n^3h_0}\mathbb{E}\left[\sum_{i=1}^m\sum_{j\neq i}^m\sum_{k=1}^m\sum_{l\neq k}^m\Gamma_i\mathcal{L}_j\Gamma_k\mathcal{L}_l\right]=\frac{m(m-1)}{n^3h_0}\mathbb{E}[\Gamma_1^2]\mathbb{E}[\mathcal{L}_2^2]=O\left(\frac{m^2}{n^3}\right) \nonumber
\end{align}
\end{proof}
\end{lemma}

\begin{lemma} \label{L:lambda AM 4m}
\begin{align}
    \mathbb{E}|\Lambda_4(x;m)|^2=O\left(\frac{m^2}{n^3}\right) \nonumber
\end{align}
\begin{proof}
Proof is similar to Lemma \ref{L:lambda AM 2m}
\end{proof}
\end{lemma}

\begin{lemma} \label{L:lambda AM 5m}
\begin{align}
    \mathbb{E}|\Lambda_5(x;m)|^2=O\left(\frac{m}{n^3}\right) \nonumber
\end{align}
\begin{proof}
\begin{align}
    \mathbb{E}|\Lambda_5(x;m)|^2 & \lesssim \frac{1}{n^3h_0}\mathbb{E}\left[\sum_{i=1}^m\sum_{j=1}^mS_i\mathcal{L}_iS_j\mathcal{L}_j\right] \nonumber\\
    & = \frac{1}{n^3h_0}\mathbb{E}\left[\sum_{i=1}^m\sum_{j\neq i}^mS_i\mathcal{L}_iS_j\mathcal{L}_j+\sum_{i=1}^mS_i^2\mathcal{L}_i^2\right] \nonumber\\
    & = \frac{m}{n^3h_0}\mathbb{E}[S_1^2\mathcal{L}_1^2]=O\left(\frac{m}{n^3}\right) \nonumber
\end{align}
\end{proof}
\end{lemma}

\begin{lemma} \label{L:lambda AM m product}
\begin{align}
    & \mathbb{E}|\Lambda_1(x;m)\Lambda_2(x;m)|=O\left(\frac{m^3h_0^{2L+1}}{n^4}\right)^{1/2}=O\left(\frac{m^{3/2}h_0^{(2L+1)/2}}{n^2}\right) \nonumber\\
    & \mathbb{E}|\Lambda_1(x;m)\Lambda_4(x;m)|=O\left(\frac{m^3h_0^{2L+1}}{n^4}\right)^{1/2}=O\left(\frac{m^{3/2}h_0^{(2L+1)/2}}{n^2}\right) \nonumber\\
    & \mathbb{E}|\Lambda_1(x;m)\Lambda_5(x;m)|=O\left(\frac{m^2h_0^{2L+1}}{n^4}\right)^{1/2}=O\left(\frac{mh_0^{(2L+1)/2}}{n^2}\right) \nonumber\\
    & \mathbb{E}|\Lambda_2(x;m)\Lambda_4(x;m)|=O\left(\frac{m^4}{n^6}\right)^{1/2}=O\left(\frac{m^2}{n^3}\right) \nonumber\\
    & \mathbb{E}|\Lambda_2(x;m)\Lambda_5(x;m)|=O\left(\frac{m^3}{n^6}\right)^{1/2}=O\left(\frac{m^{3/2}}{n^3}\right) \nonumber\\
    & \mathbb{E}|\Lambda_4(x;m)\Lambda_5(x;m)|=O\left(\frac{m^3}{n^6}\right)^{1/2}=O\left(\frac{m^{3/2}}{n^3}\right) \nonumber
\end{align}
\begin{proof}
From Cauchy-Schwarz inequality and Lemma \ref{L:lambda AM 1m}, \ref{L:lambda AM 2m}, \ref{L:lambda AM 4m} and \ref{L:lambda AM 5m}, this lemma holds.
\end{proof}
\end{lemma}

\section{Derivation of Expression for $p_1(z),p_3(z)$ and $p_4(z)$}\label{Appendix C}

For $p_1(z)$, we have,
\begin{align}
    p_1(z)&=\frac{-1}{6}\mu_{30}\mu_{20}^{-3/2}(z^2-1) \nonumber\\
    &=\frac{-1}{6}\frac{\kappa_{03}f(x)-3\kappa_{02}f(x)^2h_0+\kappa_{22}f^{(2)}(x)h_0^2/2+o(h_0^2)}{[\kappa_{02}f(x)-f(x)^2h_0+\{\kappa_{23}f^{(2)}(x)/2+2f(x)^3\}h^2+o(h_0^2)]^{3/2}}(z^2-1) \nonumber\\
    &=\frac{-1}{6}[\kappa_{03}f(x)-3\kappa_{02}f(x)^2h_0+\kappa_{22}f^{(2)}(x)h_0^2/2+o(h_0^2)] \nonumber\\
    &\qquad\times[\{\kappa_{02}f(x)\}^{-3/2} \nonumber\\
    &\qquad\qquad-\frac{3}{2}\{\kappa_{02}f(x)\}^{-5/2}(f(x)^2h_0-\{\frac{\kappa_{23}f^{(2)}(x)}{2}+2f(x)^3\}h_0^2) \nonumber\\
    &\qquad\qquad+\frac{15}{8}\{\kappa_{02}f(x)\}^{-7/2}f(x)^4h_0^2](z^2-1)+o(h_0^2) \nonumber\\
    &=\frac{-1}{6}\Biggl[\kappa_{02}^{-3/2}\kappa_{03}f(x)-3\Biggl\{\frac{f(x)^{1/2}}{\kappa_{02}^{1/2}}-\frac{\kappa_{03}f(x)^{1/2}}{2\kappa_{02}^{5/2}}\Biggl\}h_0 \nonumber\\
    &\qquad+\Biggl\{\frac{-3}{4}\{\kappa_{02}f(x)\}^{-5/2}\kappa_{03}\kappa_{23}f^{(2)}(x)f(x)-3\{\kappa_{02}f(x)\}^{-5/2}\kappa_{03}f(x)^4 \nonumber\\
    &\qquad\qquad+\frac{15}{8}\{\kappa_{02}f(x)\}^{-7/2}\kappa_{03}f(x)^5+\frac{9}{2}\kappa_{02}^{-3/2}f(x)^{3/2}\Biggl\}h_0^2\Biggl](z^2-1)+o(h_0^2) \nonumber\\
    &\equiv\gamma_{1,0}(x)(z^2-1)+\gamma_{1,1}(x)(z^2-1)h_0+\gamma_{1,2}(x)(z^2-1)h_0^2+o(h_0^2). \nonumber\label{Linearized p_1}
\end{align}

For $p_{3,0}(z)$, we have

\begin{align}
    p_{3,0}(z)&=-C_{PI}C_{\Gamma,0}(x)\rho_{11}\mu_{20}^{-1}z \nonumber\\
    &=-C_{PI}C_{\Gamma,0}(x)\frac{\mathcal{L}(x) f(x)+O(h_0^L)}{[\kappa_{02}f(x)-f(x)^2h_0+o(h_0)]^{-1}}z \nonumber\\
    &=-C_{PI}C_{\Gamma,0}(x)[\mathcal{L}(x) f(x)+O(h_0^L)] \nonumber\\
    &\qquad\times[\{\kappa_{02}f(x)\}^{-1}-\{\kappa_{02}f(x)\}^{-2}\left(f(x)^2h_0\right)]z+o(h_0) \nonumber\\
    &=-C_{PI}C_{\Gamma,0}(x)\kappa_{02}^{-1}\mathcal{L}(x) z+C_{PI}C_{\Gamma,0}(x)\kappa_{02}^{-2}\mathcal{L}(x) f(x)zh_0+o(h_0) \nonumber\\
    &\equiv \gamma_{3,1,0}(x)z+\gamma_{3,1,1}(x)zh_0+o(h_0), \nonumber
\end{align}

while for $p_4(z)$,
\begin{align}
    p_4(z)&=-C_{PI}\rho_{11}\xi_{11}\mu_{20}^{-3/2}(z^2-1)+\frac{1}{2}C_{PI}\rho_{11}\mu_{20}^{-1/2}z^2 \nonumber\\
    &=-C_{PI}\frac{\mathcal{L}(x) f(x)\left\{\tau_0f(x)+o(h_0)\right\}}{[\kappa_{02}f(x)-f(x)^2h_0+o(h_0)]^{3/2}}(z^2-1)+\frac{1}{2}C_{PI}\frac{\mathcal{L}(x) f(x)+O(h_0^L)}{[\kappa_{02}f(x)-f(x)^2h_0+o(h_0)]^{1/2}}z^2 \nonumber\\
    &=-C_{PI}\mathcal{L}(x) f(x)\left\{\tau_0f(x)+o(h_0)\right\} \nonumber\\
    &\qquad\times [\{\kappa_{02}f(x)\}^{-3/2}-\frac{3}{2}\{\kappa_{02}f(x)\}^{-5/2}f(x)^2h_0+o(h_0)\}](z^2-1) \nonumber\\
    &+\frac{1}{2}C_{PI}\mathcal{L}(x)f(x) \nonumber\\
    &\qquad\times[\{\kappa_{02}f(x)\}^{-1/2}-\frac{1}{2}\{\kappa_{02}f(x)\}^{-3/2}f(x)^2h_0+o(h_0)\}]z^2\nonumber\\
    &=-C_{PI}\kappa_{02}^{-3/2}\tau_0\mathcal{L}(x) f(x)^{1/2}(z^2-1)+\frac{3}{2}C_{PI}\kappa_{02}^{-5/2}\tau_0\mathcal{L}(x) f(x)^{3/2}(z^2-1)h_0 \nonumber \\
    &\qquad+\frac{1}{2}C_{PI}\kappa_{02}^{-1/2}\mathcal{L}(x)f(x)^{1/2}z^2-\frac{1}{4}C_{PI}\kappa_{02}^{-3/2}\mathcal{L}(x)f(x)^{3/2}z^2h_0+o(h_0)\nonumber\\
    &\equiv \{\gamma_{4,1,0}(x)(z^2-1)+\gamma_{4,2,0}(x)z^2\}+\{\gamma_{4,1,1}(x)(z^2-1)+\gamma_{4,2,1}(x)z^2\}h_0+o(h_0). \nonumber
\end{align}
    \hfill$\square$\\

\section{Formal Derivation of Theorem \ref{T:Edgeworth Expansion Including Pilot Bandwidth}}
In this section, we derive Theorem \ref{T:Edgeworth Expansion Including Pilot Bandwidth} formally. There is no guarantee for the mathematical rigor. However, one can validate Theorem \ref{T:Edgeworth Expansion Including Pilot Bandwidth} in the same way as the proof \ref{T:Main Theorem}.

\begin{align*}
    S_{pilot}(x) &= \frac{\sqrt{n\hat{h}}\Bigl({\hat{f}_{\hat{h}}(x)}-\mathbb{E}\hat{f}_{h_0}(x)\Bigl)}{\mu_{20}^{1/2}} \nonumber\\
    & =\frac{1}{\sqrt{nh_0}}\sum_{i=1}^nS_i\\
    & \qquad+\frac{C_{PI}h_0^{\frac{2L+1}{2}}}{n^{1/2}\mu_{20}^{1/2}}\sum_{i=1}^nV_i\left(\sum_{l=0}^{L-1}C_{\Gamma,l}(x)h^l\right)+\frac{C_{PI}h_0^{\frac{2L+1}{2}}}{2n^{1/2}(n-1)\mu_{20}^{1/2}}\sum_{i=1}^n\sum_{j\neq i}^nW_{ij}\left(\sum_{l=0}^{L-1}C_{\Gamma,l}(x)h^l\right)\\
    & \qquad+\frac{C_{PI}}{n^{3/2}h_0^{1/2}\mu_{20}^{1/2}}\left[\sum_{i=1}^n\sum_{j\neq i}^nV_i\Gamma_j+\sum_{i=1}^nV_i\Gamma_i\right]+\frac{C_{PI}}{2n^{3/2}(n-1)h_0^{1/2}\mu_{20}^{1/2}}\left[\sum_{i=1}^n\sum_{j\neq i}^n\sum_{k\neq i,j}^nW_{ij}\Gamma_k+\sum_{i=1}^n\sum_{j\neq i}^n\{W_{ij}\Gamma_i+W_{ij}\Gamma_j\}\right]\\
    & \qquad-\frac{C_{PI}}{2n^{3/2}h_0^{1/2}}\left[\sum_{i=1}^n\sum_{j\neq i}^nS_iV_j+\sum_{i=1}^nS_iV_i\right]-\frac{C_{PI}}{4n^{3/2}(n-1)h_0^{1/2}}\left[\sum_{i=1}^n\sum_{j\neq i}^n\sum_{j\neq i,j}^nW_{ij}S_k+\sum_{i=1}^n\sum_{j\neq i}^n\{W_{ij}S_i+W_{ij}S_j\}\right]\\
    & \qquad+o_p\{(nh_0)^{-1}\}\\
    & \equiv S(x) + \sum_{k=1}^{10}\Lambda_k(x)+o_p\{(nh_0)^{-1}\}
\end{align*}
In the following expansion of characteristic function of $S_{pilot}(x)$, we use
\begin{align*}
    & \mathbb{E}[S_1\Gamma_1]=h_0\mu_{20}^{-1/2}\xi_{11},\\
    & \mathbb{E}[S_1V_1] =\mathbb{E}[S_1\mathcal{L}_1]=h_0\mu_{20}^{-1/2}\rho_{11},\\
    & \mathbb{E}[S_1S_2W_{12}] = b^{-2L}h_0\mu_{20}^{-1}\omega_{111},\\
    & \mathbb{E}[S_1\Gamma_2W_{12}] = b^{-2L}h_0\mu_{20}^{-1/2}\psi_{111}.
\end{align*}
where the reason why $\mathbb{E}[S_1V_1]=\mathbb{E}[S_1\mathcal{L}_1]$ holds is provided as Remark \ref{R:pilot}.
Define the characteristic function of $S_{pilot}(x)$ as $\chi_{pilot}(x)$.
\begin{align*}
    \chi_{pilot}(t) & \equiv \mathbb{E}\left[\exp\left(it\left\{S(x) + \sum_{k=1}^{10}\Lambda_k(x)+o_p\{(nh_0)^{-1}\}\right\}\right)\right] \\
    & = \mathbb{E}\left[e^{itS(x)}\left\{1+\sum_{k=1}^{10}it\Lambda_k(x)\right\}\right]+ o_p\{(nh_0)^{-1}\}\\
    & \equiv (I)+(II)+(III)+(IV)+(V)+(VI)+(VII)+(VIII)+(IX)+(X)+(XI)+o_p\{(nh_0)^{-1}\}
\end{align*}
In the following subsections, we expand each component of ch.f.

\subsection{Expansion of The Second Component}
\begin{align*}
    (II) &= \mathbb{E}\left[e^{itS(x)}it\Lambda_1(x)\right]\\
    &=\frac{C_{PI}h_0^{(2L+1)/2}}{n^{1/2}\mu_{20}^{1/2}}\mathbb{E}\left[e^{itS(x)}\sum_{i=1}^nV_i\left(\sum_{l=0}^{L-1}C_{\Gamma,l}(x)h^l\right)\right](it)\\
    &= \frac{C_{PI}h_0^{\frac{2L+1}{2}}}{n^{1/2}\mu_{20}^{1/2}}n\gamma(t)^{n-1}\mathbb{E}\left[\left\{1+\frac{it}{\sqrt{nh_0}}S_1\right\}V_1\right]\left(\sum_{l=0}^{L-1}C_{\Gamma,l}(x)h_0^l\right)(it)+o\{(nh_0)^{-1}\}\\
    &= C_{PI}h_0^{L}\gamma(t)^{n-1}\mu_{20}^{-1/2}\mathbb{E}[S_1V_1]\left(\sum_{l=0}^{L-1}C_{\Gamma,l}(x)h_0^l\right)(it)^2+o\{(nh_0)^{-1}\}\\
    &= C_{PI}h_0^{L+1}\gamma(t)^{n-1}\mu_{20}^{-1}\sum_{l=0}^{L-1}C_{\Gamma,l}(x)h_0^l\rho_{11}(it)^2+o\{(nh_0)^{-1}\}
\end{align*}

\subsection{The Third Component}
\begin{align*}
    (III) &= \mathbb{E}\left[e^{itS(x)}it\Lambda_2(x)\right]\\
    &=\frac{C_{PI}h_0^{\frac{2L+1}{2}}}{2n^{1/2}(n-1)\mu_{20}^{1/2}}\mathbb{E}\left[e^{itS(x)}\sum_{i=1}^n\sum_{j\neq i}^nW_{ij}\left(\sum_{l=0}^{L-1}C_{\Gamma,l}(x)h^l\right)\right](it)\\
    &= \frac{C_{PI}}{2n^{1/2}(n-1)}\mu_{20}^{-1/2}n(n-1)\gamma(t)^{n-2}\mathbb{E}\left[e^{itS(x)}W_{12}\right]\left(\sum_{l=0}^{L-1}C_{\Gamma,l}(x)h_0^l\right)(it)+o\{(nh_0)^{-1}\}\\
    &=\frac{C_{PI}}{2}n^{1/2}h_0^{\frac{2L+1}{2}}\gamma(t)^{n-2}\mu_{20}^{-1/2}\mathbb{E}\left[\Biggl\{1+\frac{it}{\sqrt{nh_0}}(S_1+S_2)+\frac{(it)^2}{nh_0}(S_1+S_2)^2\Biggl\}W_{12}\right]\left(\sum_{l=0}^{L-1}C_{\Gamma,l}(x)h_0^l\right)(it)+o\{(nh_0)^{-1}\}\\
    &= \frac{C_{PI}h_0^{\frac{2L-1}{2}}}{2n^{1/2}}\gamma(t)^{n-2}\mu_{20}^{-1/2}\mathbb{E}[S_1S_2W_{12}]\left(\sum_{l=0}^{L-1}C_{\Gamma,l}(x)h_0^l\right)(it)^3+o\{(nh_0)^{-1}\}\\
    &=\frac{C_{PI}h_0^{\frac{2L+1}{2}}}{2n^{1/2}b^{2L}}\gamma(t)^{n-2}\mu_{20}^{-3/2}\omega_{111}\left(\sum_{l=0}^{L-1}C_{\Gamma,l}(x)h_0^l\right)(it)^3+o\{(nh_0)^{-1}\}
\end{align*}

\subsection{The Fourth Component}
\begin{align*}
    (IV) &= \mathbb{E}\left[e^{itS(x)}it\Lambda_3(x)\right]\\
    &=\frac{C_{PI}}{n^{3/2}h_0^{1/2}\mu_{20}^{1/2}}\mathbb{E}\left[e^{itS(x)}\sum_{i=1}^n\sum_{j\neq i}^n V_i\Gamma_j\right](it)\\
    &= \frac{C_{PI}}{n^{3/2}h_0^{1/2}\mu_{20}^{1/2}}n(n-1)\gamma(t)^{n-2}\mathbb{E}\left[\left\{1+\frac{it}{\sqrt{nh_0}}(S_1+S_2)+\frac{(it)^2}{nh_0}(S_1+S_2)^2\right\}V_1\Gamma_2\right](it)+o\{(nh_0)^{-1}\}\\
    & = \frac{C_{PI}}{n^{1/2}h_0^{3/2}}\gamma(t)^{n-2}\mu_{20}^{-1/2}\mathbb{E}[S_1V_1]\mathbb{E}[S_2\Gamma_2](it)^3+o\{(nh)^{-1}\}\\
    &= \frac{C_{PI}h_0^{1/2}}{n^{1/2}}\gamma(t)^{n-2}\mu_{20}^{-3/2}\rho_{11}\xi_{11}(it)^3+o\{(nh_0)^{-1}\}
\end{align*}

\subsection{The Fifth Component}

\begin{align*}
    (V)&=\frac{C_{PI}}{n^{3/2}h_0^{1/2}\mu_{20}^{1/2}}\mathbb{E}\left[e^{itS(x)}\sum_{i=1}^n V_i\Gamma_i\right](it)\\
    & = \frac{C_{PI}}{n^{1/2}h_0^{1/2}\mu_{20}^{1/2}}\mathbb{E}\left[V_1\Gamma_1\right](it) + o\{(nh_0)^{-1}\}.
\end{align*}
Since, from Lemma \ref{gamma L}, $\mathbb{E}[V_1\Gamma_1]=O(h_0^{L+1})$, $(V)=o\{(nh_0)^{-1}\}$

\subsection{The Sixth Component}
\begin{align*}
    (VI)&=\frac{C_{PI}}{2n^{3/2}(n-1)h_0^{1/2}\mu_{20}^{1/2}}\mathbb{E}\left[e^{itS(x)}\sum_{i=1}^n\sum_{j\neq i}^n\sum_{k\neq i,j}^nW_{ij}\Gamma_k\right](it)\\ 
    &=\frac{C_{PI}}{2n^{3/2}(n-1)h_0^{1/2}\mu_{20}^{1/2}}n(n-1)(n-2)\gamma(t)^{n-3}\mathbb{E}\left[e^{\frac{it}{\sqrt{nh_0}}(S_1+S_2+S_3)}W_{12}\Gamma_3\right](it)\\
    &=\frac{C_{PI}}{2nh_0^2\mu_{20}^{1/2}}\gamma(t)^{n-3}\mathbb{E}[S_1S_2W_{12}]\mathbb{E}[S_3\Gamma_3](it)^4+o\{(nh_0)^{-1}\}\\
    &= \frac{C_{PI}}{2nb^{2L}}\gamma(t)^{n-3}\mu_{20}^{-2}\xi_{11}\omega_{111}(it)^4+o\{(nh_0)^{-1}\}
\end{align*}

\subsection{The Seventh Component}
\begin{align*}
    (VII)&=\frac{C_{PI}}{2n^{3/2}(n-1)h_0^{1/2}\mu_{20}^{1/2}}\mathbb{E}\left[e^{itS(x)}\sum_{i=1}^n\sum_{j\neq i}^n\{W_{ij}\Gamma_i+W_{ij}\Gamma_j\}\right](it)\\  &= \frac{C_{PI}}{2n^{3/2}(n-1)h_0^{1/2}\mu_{20}^{1/2}}n(n-1)\gamma(t)^{n-2}\mathbb{E}\left[e^{\frac{it}{\sqrt{nh_0}}(S_1+S_2)}\{W_{12}\Gamma_1+W_{12}\Gamma_2\}\right](it)\\
    &= \frac{C_{PI}}{nh_0}\gamma(t)^{n-2}\mu_{20}^{-1/2}\mathbb{E}[W_{12}\Gamma_1S_2](it)^2+o\{(nh_0)^{-1}\}\\
    &= \frac{C_{PI}}{nb^{2L}}\gamma(t)^{n-2}\mu_{20}^{-1}\psi_{111}(it)^2+o\{(nh_0)^{-1}\}
\end{align*}

\subsection{The Eighth Component}
\begin{align*}
    (VIII)&= \frac{-C_{PI}}{2n^{3/2}h_0^{1/2}}\mathbb{E}\left[e^{itS(x)}\sum_{i=1}^n\sum_{j\neq i}^nS_iV_j\right](it)\\
    &=\frac{-C_{PI}}{2n^{3/2}h_0^{1/2}}n(n-1)\gamma(t)^{n-2}\mathbb{E}\left[e^{\frac{it}{\sqrt{nh_0}}(S_1+S_2)}S_1V_2\right](it)\\
    &=\frac{-C_{PI}}{2n^{1/2}h_0^{3/2}}\gamma(t)^{n-2}\mathbb{E}[S_1^2]\mathbb{E}[S_2V_2](it)^3+o\{(nh_0)^{-1}\}\\
    &=\frac{-C_{PI}h_0^{1/2}}{2n^{1/2}}\gamma(t)^{n-2}\mu_{20}^{-1/2}\rho_{11}(it)^3+o\{(nh_0)^{-1}\}
\end{align*}

\subsection{The Ninth Component}
\begin{align*}
    (IX) &=\frac{-C_{PI}}{2n^{3/2}h_0^{1/2}}\mathbb{E}\left[e^{itS(x)}\sum_{i=1}^nS_iV_i\right](it)\\
    &= \frac{-C_{PI}}{2n^{3/2}h_0^{1/2}}n\gamma(t)^{n-1}\mathbb{E}\left[e^{\frac{it}{\sqrt{nh_0}}S_1}S_1V_1\right](it)\\
    &=\frac{-C_{PI}}{2n^{1/2}h_0^{1/2}}\gamma(t)^{n-1}\mathbb{E}[S_1V_1](it)+o\{(nh)^{-1}\}\\
    &=\frac{-C_{PI}h_0^{1/2}}{2n^{1/2}}\gamma(t)^{n-1}\mu_{20}^{-1/2}\rho_{11}(it)+o\{(nh_0)^{-1}\}
\end{align*}

\subsection{The Tenth Component}
\begin{align*}
    (X) &=\frac{-C_{PI}}{4n^{3/2}(n-1)h_0^{1/2}}\mathbb{E}\left[e^{itS(x)}\sum_{i=1}^n\sum_{j\neq i}^n\sum_{k\neq i,j}^nW_{ij}S_k\right](it)\\
    &= \frac{-C_{PI}}{4n^{3/2}(n-1)h_0^{1/2}}n(n-1)(n-2)\gamma(t)^{n-3}\mathbb{E}\left[e^{\frac{it}{\sqrt{nh_0}}(S_1+S_2+S_3)}W_{12}S_3\right](it)\\
    &= \frac{C_{PI}(n-2)}{4n^{1/2}h_0^{1/2}}\gamma(t)^{n-3}\mathbb{E}\left[\left\{1+\cdots+\frac{(it)^3}{(nh_0)^{3/2}}W_{12}S_3\right\}\right](it)\\
    &= \frac{-C_{PI}}{4nh_0^2}\gamma(t)^{n-3}\mathbb{E}[S_1S_2W_{12}]\mathbb{E}[S_3^2](it)^4 +o\{(nh)^{-1}\}\\
    &= \frac{-C_{PI}}{4nb^{2L}}\gamma(t)^{n-3}\mu_{20}^{-1}\omega_{111}(it)^4+o\{(nh_0)^{-1}\}
\end{align*}

\subsection{The Eleventh Comopnent}
\begin{align*}
    (XI) &=\frac{-C_{PI}}{4n^{3/2}(n-1)h_0^{1/2}}\mathbb{E}\left[e^{itS(x)}\sum_{i=1}^n\sum_{j\neq i}^n\{W_{ij}S_i+W_{ij}S_j\}\right](it)\\
    &= \frac{-C_{PI}}{4n^{3/2}(n-1)h_0^{1/2}}n(n-1)\gamma(t)^{n-2}\mathbb{E}\left[e^{\frac{it}{\sqrt{nh_0}}(S_1+S_2)}\left\{W_{12}S_1+W_{12}S_2\right\}\right]\\
    &= \frac{-C_{PI}}{2nh_0}\gamma(t)^{n-2}\mathbb{E}[S_1S_2W_{12}](it)^2 +o\{(nh_0)^{-1}\}\\
    &= \frac{-C_{PI}}{2nb^{2L}}\gamma(t)^{n-2}\mu_{20}^{-1}\omega_{111}(it)^2+o\{(nh_0)^{-1}\}
\end{align*}

Recalling
\begin{align}
    \chi_{pilot}(t) 
    &= (I)+(II)+(III)+(IV)+(V)+(VI)+(VII)+(VIII)+(IX)+(X)+(XI)+o_p\{(nh_0)^{-1}\}, \nonumber
\end{align}
we have
\begin{align*}
    \chi_{pilot}(t)&=\exp\left(\frac{-t^2}{2}\right)\Biggl[\Biggl\{1+\frac{\mu_{30}\mu_{20}^{-3/2}}{6n^{1/2}h_0^{1/2}}(it)^3+\frac{\mu_{40}\mu_{20}^{-2}}{24nh_0}(it)^4+\frac{\mu_{30}^2\mu_{20}^{-3}}{72nh_0}(it)^6\Biggl\}\\
    &\qquad+C_{PI}\mu_{20}^{-1}\rho_{11}\left(\sum_{l=0}^{L-1}C_{\Gamma,l}(x)h_0^{L+l+1}\right)(it)^2\\
    &\qquad+\frac{C_{PI}}{2}n^{-1/2}h_0^{\frac{2L+1}{2}}b^{-2L}\mu_{20}^{-3/2}\omega_{111}\left(\sum_{l=0}^{L-1}C_{\Gamma,l}(x)h_0^{l}\right)(it)^3\\
    &\qquad+C_{PI}n^{-1/2}h_0^{1/2}\Biggl(\mu_{20}^{-3/2}\rho_{11}\xi_{11}(it)^3 -\frac{\mu_{20}^{-1/2}}{2}\rho_{11}\{(it)^3+(it)\}\Biggl) \nonumber\\
    &\qquad+C_{PI}n^{-1}b^{-2L}\left(\frac{1}{2}\mu_{20}^{-2}\xi_{11}\omega_{111}(it)^4+\mu_{20}^{-1}\psi_{111}(it)^2-\frac{1}{4}\mu_{20}^{-1}\omega_{111}\{(it)^4+2(it)^2\}\right)\Biggl]+o\{(nh_0)^{-1}\}\\
    & = \exp\left(\frac{-t^2}{2}\right)\Biggl[\Biggl\{1+\frac{\mu_{30}\mu_{20}^{-3/2}}{6n^{1/2}h_0^{1/2}}(it)^3+\frac{\mu_{40}\mu_{20}^{-2}}{24nh_0}(it)^4+\frac{\mu_{30}^2\mu_{20}^{-3}}{72nh_0}(it)^6\Biggl\}\\
    &\qquad+C_{PI}\mu_{20}^{-1}\rho_{11}\left(\sum_{l=0}^{L-1}C_{\Gamma,l}(x)h_0^{L+l+1}\right)(it)^2 +\frac{C_{PI}}{2}n^{-1/2}h_0^{\frac{2L+1}{2}}b^{-2L}\mu_{20}^{-3/2}\omega_{111}\left(\sum_{l=0}^{L-1}C_{\Gamma,l}(x)h_0^{l}\right)(it)^3\\
    &\qquad+C_{PI}n^{-1/2}h_0^{1/2}\rho_{11}\left(\mu_{20}^{-3/2}\xi_{11}(it)^3-\frac{\mu_{20}^{-1/2}}{2}\{(it)^3+(it)\}\}\right) \nonumber\\
    &\qquad+C_{PI}n^{-1}b^{-2L}\left(\frac{1}{2}\mu_{20}^{-2}\xi_{11}\omega_{111}(it)^4+\mu_{20}^{-1}\psi_{111}(it)^2-\frac{1}{4}\mu_{20}^{-1}\omega_{111}\{(it)^4+2(it)^2\}\right)\Biggl] +o\{(nh_0)^{-1}\}
\end{align*}
Inverting $\chi_{pilot}(t)$, we have Theorem \ref{T:Edgeworth Expansion Including Pilot Bandwidth}.

\section{Formal Derivation of Theorem \ref{T:Edgeworth Expansion Including Pilot Bandwidth and the effect of Studentisation}}
In this section, we derive Theorem \ref{T:Edgeworth Expansion Including Pilot Bandwidth and the effect of Studentisation} formally. There is no guarantee for the mathematical rigor. However, one can validate Theorem \ref{T:Edgeworth Expansion Including Pilot Bandwidth and the effect of Studentisation} in the same way as the proof \ref{T:Main Theorem}.

\subsection{Expansion of $\hat{\mu}_{20}(\hat{h})$}
Let $\hat{\mu}_{20}$ be the natural estimator for $\mu_{20}$,
\begin{align}
    \hat{\mu}_{20}(h)\equiv h^{-1}\left\{\frac{1}{n}\sum_{i=1}^nK_{i,h}(x)^2-\left(\frac{1}{n}\sum_{i=1}^nK_{i,h}(x)\right)^2\right\}, \nonumber
\end{align}
then studentized KDE with data-driven bandwidth $\hat{h}$, the standard deviation is $\hat{\mu}_{20}(\hat{h})$. Under Assumption \ref{A:L_p order 4}, expanding $\hat{\mu}_{20}(\hat{h})$ around $\hat{h}=h_0$ yields
\begin{align}
    \hat{\mu}_{20}(\hat{h})=\hat{\mu}_{20}(h_0)+\hat{\mu}_{20,\partial h}(h_0)(\hat{h}-h_0)+o_p\{(nh_0)^{-1}\}.
\end{align}
where the definition of $\hat{\mu}_{20,\partial h}(h_0)$ is
\begin{align}
    & \hat{\mu}_{20,\partial h}(h_0) = -h_0^{-2}\left\{\frac{1}{n}\sum_{i=1}^nK_{i,h_0}(x)^2-\left(\frac{1}{n}\sum_{i=1}^nK_{i,h_0}(x)\right)^2\right\}\\
    &\quad -h_0^{-2}\left\{\frac{2}{n}\sum_{i=1}^nK_{i,h_0}(x)K'_{i,h_0}(x)u_{i,h_0}(x)-\frac{2}{n^2}\sum_{i=1}^n\sum_{j=1}^nK_{i,h_0}(x)K'_{j,h_0}u_{j,h_0}\right\}.
\end{align}
\subsubsection{Transformation of $\hat{\mu}_{20}(h_0)$}
We first transform $\hat{\mu}_{20}(h_0)$. (\cite{Hall91} and \cite{Hall92book} has already done this transformation, see \cite[p.212-213]{Hall92book}.)
Define
\begin{align}
    \Delta_j(h) = \frac{1}{\sqrt{nh}}\sum_{i=1}^n\left\{K_{i,h}(x)^j-\mathbb{E}K_{i,h}(x)^j\right\}. \nonumber
\end{align}
Then, $\hat{\mu}_{20}(h_0)$ is 
\begin{align}
    \hat{\mu}_{20}(h_0) &= h_0^{-1}\left\{\frac{1}{n}\sum_{i=1}^nK_{i,h_0}(x)^2-\left(\frac{1}{n}\sum_{i=1}^nK_{i,h_0}(x)\right)^2\right\} \nonumber\\
    & = h_0^{-1}\Biggl\{\mathbb{E}\left[K_{i,h_0}(x)^2\right]+\frac{1}{n}\sum_{i=1}^nK_{i,h_0}(x)^2-\mathbb{E}\left[K_{i,h_0}(x)^2\right]+\mathbb{E}\left[K_{i,h_0}(x)\right]^2-\left(\frac{1}{n}\sum_{i=1}^nK_{i,h_0}(x)\right)^2-\mathbb{E}\left[K_{i,h_0}(x)\right]^2\Biggl\} \nonumber\\
    & = h_0^{-1}\Biggl\{\mathbb{E}\left[K_{i,h_0}(x)^2\right]-\mathbb{E}\left[K_{i,h_0}(x)\right]^2\Biggl\} \nonumber\\
    & \quad + h_0^{-1}\left\{\frac{1}{n}\sum_{i=1}^nK_{i,h_0}(x)^2-\mathbb{E}\left[K_{i,h_0}(x)^2\right]\right\}-h_0^{-1}\left\{\left(\frac{1}{n}\sum_{i=1}^nK_{i,h_0}(x)\right)^2-\mathbb{E}\left[K_{i,h_0}(x)\right]^2\right\} \nonumber\\
    &= \mu_{20}(h_0) + (nh_0)^{-1/2}\Delta_2(h_0) - 2\frac{h_0f(x)}{\sqrt{nh_0}}\Delta_1(h_0)+O_p(n^{-1}). \label{expanded mu 20}
\end{align}
where the third term in the final quality follows from the following equation.
\begin{align*}
    &\Delta_1(h_0)^2 = nh_0^{-1}\left\{\left(\frac{1}{n}\sum_{i=1}^nK_{i,h_0}(x)\right)^2-\mathbb{E}\left[K_{i,h_0}(x)\right]^2\right\} -2h_0^{-1}\mathbb{E}\left[K_{i,h_0}(x)\right]\Delta_1(h_0)\\
    & \implies h_0^{-1}\left\{\left(\frac{1}{n}\sum_{i=1}^nK_{i,h_0}(x)\right)^2-\mathbb{E}\left[K_{i,h_0}(x)\right]^2\right\}= n^{-1}\Delta_1(h_0)^2+2\frac{\mathbb{E}\left[K_{i,h_0}(x)\right]}{\sqrt{nh_0}}\Delta_1(h_0)
\end{align*}

\subsubsection{$\hat{\mu}_{20,\partial h}(h_0)$}
Next, we have to transform $\hat{\mu}_{20,\partial h}(h_0)$. Define
\begin{align*}
    & \Psi^{(1)}(h)\equiv \frac{1}{\sqrt{nh}}\sum_{i=1}^n\left\{K_{i,h}(x)K'_{i,h}(x)u_{i,h}(x)-\mathbb{E}\left[K_{i,h}(x)K'_{i,h}(x)u_{i,h}(x)\right]\right\}\\
    & \Delta_1^{(1)}(h) = \frac{1}{\sqrt{nh}}\sum_{i=1}^n\left\{K'_{i,h}(x)u_{i,h}(x)-\mathbb{E}\left[K'_{i,h}(x)u_{i,h}(x)\right]\right\}\\
    & \delta(h) = h^{-1}\left\{\mathbb{E}\left[K_{i,h}(x)K'_{i,h}(x)u_{i,h}(x)\right]-\mathbb{E}\left[K_{i,h}(x)K'_{j,h}(x)u_{j,h}(x)\right]\right\} 
\end{align*}
Since $\hat{h}-h_0=O_p\{n^{-1/2}h_0 \lor n^{-1}b^{-(4L+1)/2}h_0\}=O_p\{n^{-1/2}h_0 \lor n^{\frac{-2L_p}{4L+2L_p+1}}h_0\}$, we can ignore the terms whose convergence rates are faster than $O_p\{n^{-1/2}h_0^{-2} \land b^{(4L+1)/2}h_0^{-2}\}$ for Edgeworth expansion up to the order of $O\{(nh_0)^{-1}\}$.
\begin{align}
    & \hat{\mu}_{20,\partial h}(h_0) = -h_0^{-2}\left\{\frac{1}{n}\sum_{i=1}^nK_{i,h_0}(x)^2-\left(\frac{1}{n}\sum_{i=1}^nK_{i,h_0}(x)\right)^2\right\} \nonumber\\
    &\quad -h_0^{-2}\left\{\frac{2}{n}\sum_{i=1}^nK_{i,h_0}(x)K'_{i,h_0}(x)u_{i,h_0}(x)-\frac{2}{n^2}\sum_{i=1}^n\sum_{j=1}^nK_{i,h_0}(x)K'_{j,h_0}(x)u_{j,h_0}(x)\right\} \nonumber\\
    & = -h_0^{-1}\hat{\mu}_{20}(h_0) \nonumber\\
    &\quad -2h_0^{-2}\left\{\mathbb{E}\left[\frac{1}{n}\sum_{i=1}^nK_{i,h_0}(x)K'_{i,h_0}(x)u_{i,h_0}(x)\right]-\mathbb{E}\left[\frac{1}{n^2}\sum_{i=1}^n\sum_{j=1}^nK_{i,h_0}(x)K'_{j,h_0}(x)u_{j,h_0}(x)\right]\right\} \nonumber\\
    &\quad -2h_0^{-2}\left\{\frac{1}{n}\sum_{i=1}^nK_{i,h_0}(x)K'_{i,h_0}(x)u_{i,h_0}(x)-\mathbb{E}\left[\frac{1}{n}\sum_{i=1}^nK_{i,h_0}(x)K'_{i,h_0}(x)u_{i,h_0}(x)\right]\right\} \nonumber\\
    &\quad +2h_0^{-2}\left\{\frac{1}{n^2}\sum_{i=1}^n\sum_{j=1}^nK_{i,h_0}(x)K'_{j,h_0}(x)u_{j,h_0}(x)-\mathbb{E}\left[\frac{1}{n^2}\sum_{i=1}^n\sum_{j=1}^nK_{i,h_0}(x)K'_{j,h_0}(x)u_{j,h_0}(x)\right]\right\} \nonumber\\
    & = -h_0^{-1}\hat{\mu}_{20}(h_0) \nonumber\\
    &\quad -2h_0^{-2}\left\{\left(1-\frac{1}{n}\right)\mathbb{E}\left[K_{i,h_0}(x)K'_{i,h_0}(x)u_{i,h_0}(x)\right]-\left(1-\frac{1}{n}\right)\mathbb{E}\left[K_{i,h_0}(x)K'_{j,h_0}(x)u_{j,h_0}(x)\right]\right\} \label{mu20 first derivative 1}\\
    &\quad -2h_0^{-2}\left\{\frac{1}{n}\sum_{i=1}^nK_{i,h_0}(x)K'_{i,h_0}(x)u_{i,h_0}(x)-\mathbb{E}\left[\frac{1}{n}\sum_{i=1}^nK_{i,h_0}(x)K'_{i,h_0}(x)u_{i,h_0}(x)\right]\right\} \label{mu20 first derivative 2}\\
    &\quad +2h_0^{-2}\left\{\frac{1}{n^2}\sum_{i=1}^n\sum_{j\neq i}^nK_{i,h_0}(x)K'_{j,h_0}(x)u_{j,h_0}(x)-\mathbb{E}\left[\frac{1}{n^2}\sum_{i=1}^n\sum_{j\neq i}^nK_{i,h_0}(x)K'_{j,h_0}(x)u_{j,h_0}(x)\right]\right\} \label{mu20 first derivative 3}\\
    &\quad - 2h_0^{-2}\left\{\frac{1}{n^2}\sum_{i=1}^nK_{i,h_0}(x)K'_{i,h_0}(x)u_{i,h_0}(x)-\mathbb{E}\left[\frac{1}{n^2}\sum_{i=1}^nK_{i,h_0}(x)K'_{i,h_0}(x)u_{i,h_0}(x)\right]\right\} \label{mu20 first derivative 4}\\
    &= -h_0^{-1}\left\{\hat{\mu}_{20}(h_0)+2\delta(h_0)+2(nh_0)^{-1/2}\Psi^{(1)}(h_0)\right\}-2(nh_0)^{-1/2}f(x)\Delta_1^{(1)}(h_0)+2(nh_0)^{-1/2}f(x)\Delta_1(h_0) \label{mu20 first derivative} \\
    & \quad + o_p\{n^{-1/2}h_0^{-2} \land b^{(4L+1)/2}h_0^{-2}\} \nonumber \\
    &= -h_0^{-1}\hat{\mu}_{20}(h_0)-2h_0^{-1}\delta(h_0) + o_p\{n^{-1/2}h_0^{-2} \land b^{(4L+1)/2}h_0^{-2}\} \nonumber\\
    &= -h_0^{-1}\mu_{20}(h_0)-2h_0^{-1}\delta(h_0)+o_p\{n^{-1/2}h_0^{-2} \land b^{(4L+1)/2}h_0^{-2}\} \nonumber
\end{align}
Note that (\ref{mu20 first derivative 4}) is $o_p\{n^{-1/2}h_0^{-2} \land b^{(4L+1)/2}h_0^{-2}\}$. The second and third terms in (\ref{mu20 first derivative}) are another expression of (\ref{mu20 first derivative 1}) and (\ref{mu20 first derivative 2}) respectively. The last two terms of (\ref{mu20 first derivative}) follows from the following transformation of (\ref{mu20 first derivative 3}).
\begin{align*}
    (nh_0)^{-1}\Delta_1(h_0)\Delta_1^{(1)}(h_0) &= \frac{1}{n^2h_0^2}\sum_{i=1}^n\sum_{j=1}^n\left\{K_{i,h_0}(x)-\mathbb{E}\left[K_{i,h_0}(x)\right]\right\}\left\{K'_{j,h_0}(x)u_{j,h_0}(x)-\mathbb{E}\left[K'_{j,h_0}(x)u_{j,h_0}(x)\right]\right\}\\
    & = \frac{1}{n^2h_0^2}\sum_{i=1}^n\sum_{j=1}^n\left\{K_{i,h_0}(x)K'_{j,h_0}(x)u_{j,h_0}(x)-\mathbb{E}\left[K_{i,h_0}(x)\right]\mathbb{E}\left[K'_{j,h_0}(x)u_{j,h_0}(x)\right]\right\}\\
    &\quad -\frac{1}{nh_0^2}\mathbb{E}\left[K_{i,h_0}(x)\right]\sum_{i=1}^n\left\{K'_{i,h_0}(x)u_{i,h_0}(x)-\mathbb{E}\left[K'_{i,h_0}(x)u_{i,h_0}(x)\right]\right\}\\
    &\quad - \frac{1}{nh_0^2}\mathbb{E}\left[K'_{i,h_0}(x)u_{i,h_0}(x)\right]\sum_{i=1}^n\left\{K_{i,h_0}(x)-\mathbb{E}\left[K_{i,h_0}(x)\right]\right\}\\
    &= \frac{1}{n^2h_0^2}\sum_{i=1}^n\sum_{j\neq i}^n\left\{K_{i,h_0}(x)K'_{j,h_0}(x)u_{j,h_0}(x)-\mathbb{E}\left[K_{i,h_0}(x)K'_{j,h_0}(x)u_{j,h_0}(x)\right]\right\}\\
    &\quad +\frac{1}{n^2h_0^2}\sum_{i=1}^n\left\{K_{i,h_0}(x)K'_{i,h_0}(x)u_{i,h_0}(x)-\mathbb{E}\left[K_{i,h_0}(x)\right]\mathbb{E}\left[K'_{i,h_0}(x)u_{i,h_0}(x)\right]\right\}\\
    &\quad -\frac{1}{nh_0^2}\mathbb{E}\left[K_{i,h_0}(x)\right]\sum_{i=1}^n\left\{K'_{i,h_0}(x)u_{i,h_0}(x)-\mathbb{E}\left[K'_{i,h_0}(x)u_{i,h_0}(x)\right]\right\}\\
    &\quad - \frac{1}{nh_0^2}\mathbb{E}\left[K'_{i,h_0}(x)u_{i,h_0}(x)\right]\sum_{i=1}^n\left\{K_{i,h_0}(x)-\mathbb{E}\left[K_{i,h_0}(x)\right]\right\}\\
    &= \frac{1}{n^2h_0^2}\sum_{i=1}^n\sum_{j\neq i}^n\left\{K_{i,h_0}(x)K'_{j,h_0}(x)u_{j,h_0}(x)-\mathbb{E}\left[K_{i,h_0}(x)K'_{j,h_0}(x)u_{j,h_0}(x)\right]\right\}\\
    &\quad + O_p\{(nh_0)^{-3/2}\} - n^{-1/2}h_0^{-3/2}\mathbb{E}\left[K_{i,h_0}(x)\right]\Delta_1^{(1)}(h_0)-n^{-1/2}h_0^{-3/2}\mathbb{E}\left[K'_{i,h_0}(x)u_{i,h_0}(x)\right]\Delta_1(h_0)\\
    \implies &\frac{1}{n^2h_0^2}\sum_{i=1}^n\sum_{j\neq i}^n\left\{K_{i,h_0}(x)K'_{j,h_0}(x)u_{j,h_0}(x)-\mathbb{E}\left[K_{i,h_0}(x)K'_{j,h_0}(x)u_{j,h_0}(x)\right]\right\}\\
    &= (nh_0)^{-1}\Delta_1(h_0)\Delta_1^{(1)}(h_0) + (nh_0)^{-1/2}f(x)\Delta_1^{(1)}(h_0)\\
    &\quad -(nh_0)^{-1/2}f(x)\Delta_1(h_0)+O_p(n^{-1/2}h_0^{L-3/2})+O_p\{(nh_0)^{-3/2}\}\\
    &= (nh_0)^{-1/2}f(x)\Delta_1^{(1)}(h_0)-(nh_0)^{-1/2}f(x)\Delta_1(h_0) + s.o.
\end{align*}

\subsection{$T_{PI}(x)$}
We write the studentized KDE with plug-in bandwidth as $T_{PI}$. Then, from (\ref{expanded mu 20}) and (\ref{mu20 first derivative}), 
\begin{align*}
    T_{PI}(x)&\equiv\frac{\sqrt{n\hat{h}}\{\hat{f}_{\hat{h}}(x)-\mathbb{E}\hat{f}_{h}(x)\}}{\Bigl[\hat{\mu}_{20}(h_0)+\hat{\mu}_{20,\partial h}(h_0)(\hat{h}-h_0)+o_p\{(nh_0)^{-1}\}\Bigl]^{1/2}}\\
    &= \frac{\sqrt{n\hat{h}}\{\hat{f}_{\hat{h}}(x)-\mathbb{E}\hat{f}_{h}(x)\}}{\hat{\mu}_{20}(h_0)^{1/2}}-\frac{\sqrt{n\hat{h}}\{\hat{f}_{\hat{h}}(x)-\mathbb{E}\hat{f}_{h}(x)\}}{2\hat{\mu}_{20}(h_0)^{3/2}}\left\{\hat{\mu}_{20,\partial h}(h_0)(\hat{h}-h_0)+o_p\{(nh_0)^{-1}\}\right\}+o_p\{(nh_0)^{-1}\}\\
    & = \frac{\sqrt{n\hat{h}}\{\hat{f}_{\hat{h}}(x)-\mathbb{E}\hat{f}_{h}(x)\}}{\left[\mu_{20}(h_0)+(nh_0)^{-1/2}\Delta_2(h_0)-\frac{2h_0f(x)}{\sqrt{nh_0}}\Delta_1(h_0)+o_p\{(nh_0)^{-1}\}\right]^{1/2}}\\
    &\quad +\frac{\sqrt{n\hat{h}}\{\hat{f}_{\hat{h}}(x)-\mathbb{E}\hat{f}_{h}(x)\}}{2\left[\mu_{20}(h_0)+(nh_0)^{-1/2}\Delta_2(h_0)-\frac{2h_0f(x)}{\sqrt{nh_0}}\Delta_1(h_0)+o_p\{(nh_0)^{-1}\}\right]^{3/2}}\left\{\mu_{20}(h_0)+2\delta(h_0)+o_p(n^{-1/2}h_0^{-1})\right\}\left(\frac{\hat{h}-h_0}{h_0}\right)\\
    &\quad +o_p\{(nh_0)^{-1}\}\\
    &= \frac{\sqrt{n\hat{h}}\{\hat{f}_{\hat{h}}(x)-\mathbb{E}\hat{f}_{h}(x)\}}{\mu_{20}(h_0)^{1/2}}-\frac{\sqrt{n\hat{h}}\{\hat{f}_{\hat{h}}(x)-\mathbb{E}\hat{f}_{h}(x)\}}{2\mu_{20}(h_0)^{3/2}}\left\{(nh_0)^{-1/2}\Delta_2(h_0)-\frac{2h_0f(x)}{\sqrt{nh_0}}\Delta_1(h_0)\right\}\\
    &\quad +\frac{3\sqrt{n\hat{h}}\{\hat{f}_{\hat{h}}(x)-\mathbb{E}\hat{f}_{h}(x)\}}{8\mu_{20}^{5/2}}(nh_0)^{-1}\Delta_2(h_0)^2 +\frac{\sqrt{n\hat{h}}\{\hat{f}_{\hat{h}}(x)-\mathbb{E}\hat{f}_{h}(x)\}}{2\mu_{20}(h_0)^{3/2}}\left\{\mu_{20}(h_0)+2\delta(h_0)\right\}\left(\frac{\hat{h}-h_0}{h_0}\right)+o_p\{(nh_0)^{-1}\}.
\end{align*}
Since 
\begin{align*}
    &\sqrt{n\hat{h}}\{\hat{f}_{\hat{h}}(x)-\mathbb{E}\hat{f}_{h_0}(x)\}\\
    & = \sqrt{nh_0}\{\hat{f}_{h_0}(x)-\mathbb{E}\hat{f}_{h_0}(x)\}-\sqrt{nh_0}\left(\frac{\hat{h}-h_0}{h_0}\right)\Gamma_{KDE_1}+\frac{\sqrt{nh_0}\{\hat{f}_{h_0}(x)-\mathbb{E}\hat{f}_{h_0}(x)\}}{2}\left(\frac{\hat{h}-h_0}{h_0}\right)+o_p\{(nh_0)^{-1}\},
\end{align*}
$T_{PI}(x)$ is 
\begin{align*}
    T_{PI}(x) &= \frac{\sqrt{nh_0}\{\hat{f}_{h_0}(x)-\mathbb{E}\hat{f}_{h_0}(x)\}}{\mu_{20}(h_0)^{1/2}}-\frac{\sqrt{nh_0}}{\mu_{20}(h_0)^{1/2}}\left(\frac{\hat{h}-h_0}{h_0}\right)\Gamma_{KDE_1}+\frac{\sqrt{nh_0}\{\hat{f}_{h_0}(x)-\mathbb{E}\hat{f}_{h_0}(x)\}}{2\mu_{20}(h_0)^{1/2}}\left(\frac{\hat{h}-h_0}{h_0}\right)\\
    &\qquad-\frac{\sqrt{nh_0}\{\hat{f}_{h_0}(x)-\mathbb{E}\hat{f}_{h_0}(x)\}}{2\mu_{20}(h_0)^{3/2}}\left\{(nh_0)^{-1/2}\Delta_2(h_0)-\frac{2h_0f(x)}{\sqrt{nh_0}}\Delta_1(h_0)\right\}+\frac{3\sqrt{nh_0}\{\hat{f}_{h_0}(x)-\mathbb{E}\hat{f}_{h_0}(x)\}}{8\mu_{20}(h_0)^{5/2}}(nh_0)^{-1}\Delta_2(h_0)^2\\
    &\qquad+\frac{\sqrt{nh_0}\{\hat{f}_{h_0}(x)-\mathbb{E}\hat{f}_{h_0}(x)\}}{2\mu_{20}(h_0)^{3/2}}\left\{\mu_{20}(h_0)+2\delta(h_0)\right\}\left(\frac{\hat{h}-h_0}{h_0}\right)+o_p\{(nh)^{-1}\}\\
    &= \frac{\Delta_1(h_0)}{\mu_{20}(h_0)^{1/2}}-\frac{\sqrt{nh_0}}{\mu_{20}(h_0)^{1/2}}\left(\frac{\hat{h}-h_0}{h_0}\right)\{\mathbb{E}\Gamma_{KDE_1}+(\Gamma_{KDE_1}-\mathbb{E}\Gamma_{KDE_1})\}+\frac{\Delta_1(h_0)}{2\mu_{20}(h_0)^{1/2}}\left(\frac{\hat{h}-h_0}{h_0}\right)\\
    &\qquad-\frac{\Delta_1(h_0)}{2\mu_{20}(h_0)^{3/2}}\left\{(nh_0)^{-1/2}\Delta_2(h_0)-\frac{2h_0f(x)}{\sqrt{nh_0}}\Delta_1(h_0)\right\}+\frac{3(nh_0)^{-1}\Delta_1(h_0)\Delta_2(h_0)^2}{8\mu_{20}(h_0)^{5/2}}\\
    &\qquad+\frac{\Delta_1(h_0)}{2\mu_{20}(h_0)^{3/2}}\left\{\mu_{20}(h_0)+2\delta(h_0)\right\}\left(\frac{\hat{h}-h_0}{h_0}\right)+o_p\{(nh)^{-1}\}\\
    &\equiv S(x)+\Lambda_1(x)+\Lambda_2(x)+\Lambda_3(x)+\Lambda_4(x)+\Lambda_5(x)+\Lambda_6(x)+\Lambda_7(x)+\Lambda_8(x)+o_p\{(nh_0)^{-1}\}
\end{align*}
\begin{itemize}
    \item $S(x)$ is standardized KDE.
    \item $\Lambda_1(x),\Lambda_2(x)$ and $\Lambda_3(x)$ include the effect of global plug-in bandwidth.
    \item $\Lambda_4(x),\Lambda_5(x)$ and $\Lambda_6(x)$ include the effect of studentization. 
    \item $\Lambda_7(x)$ and $\Lambda_8(x)$ interaction of studentization and global plug-in bandwidth.
\end{itemize}

\subsection{Review of the Definitions of $S(x)$ and $\Lambda(x)$s.}
Define $T_i$ as follows.
\begin{align}
    T_i \equiv K_{i,h_0}(x)^2 - \mathbb{E}K_{i,h_0}(x)^2. \nonumber
\end{align}
Then, the definitions of $S(x)$ and $\Lambda(x)'$s are
\begin{align*}
     S(x) &= \Delta_1\mu_{20}^{-1/2} = \frac{1}{(nh_0)^{1/2}}\sum_{i=1}^nS_i\\
     \Lambda_1(x) &= \sqrt{nh_0}\mu_{20}^{-1/2}\left(\frac{\hat{h}-h_0}{h_0}\right)\mathbb{E}\Gamma_{KDE_1}\\
     &\qquad = \frac{C_{PI}h_0^{\frac{2L+1}{2}}}{n^{1/2}\mu_{20}^{1/2}}\sum_{i=1}^nV_i\left(\sum_{l=0}^{L-1}C_{\Gamma,l}(x)h^l\right)+\frac{C_{PI}h_0^{\frac{2L+1}{2}}}{2n^{1/2}(n-1)\mu_{20}^{1/2}}\sum_{i=1}^n\sum_{j\neq i}^nW_{ij}\left(\sum_{l=0}^{L-1}C_{\Gamma,l}(x)h^l\right)\\
     \Lambda_2 (x) &= \sqrt{nh_0}\mu_{20}^{-1/2}\left(\frac{\hat{h}-h_0}{h_0}\right)\{\Gamma_{KDE_1}-\mathbb{E}\Gamma_{KDE_1}\}\\
     &\qquad =\frac{C_{PI}}{n^{3/2}h_0^{1/2}\mu_{20}^{1/2}}\left[\sum_{i=1}^n\sum_{j\neq i}^nV_i\Gamma_j+\sum_{i=1}^nV_i\Gamma_i\right]+\frac{C_{PI}}{2n^{3/2}(n-1)h_0^{1/2}\mu_{20}^{1/2}}\left[\sum_{i=1}^n\sum_{j\neq i}^n\sum_{k\neq i,j}^nW_{ij}\Gamma_k+\sum_{i=1}^n\sum_{j\neq i}^n\{W_{ij}\Gamma_i+W_{ij}\Gamma_j\}\right]\\
     \Lambda_3(x) &=  \frac{1}{2}\mu_{20}^{-1/2}\Delta_1\left(\frac{\hat{h}-h_0}{h_0}\right)\\
     &\qquad = -\frac{C_{PI}}{2n^{3/2}h_0^{1/2}}\left[\sum_{i=1}^n\sum_{j\neq i}^nS_iV_j+\sum_{i=1}^nS_iV_i\right]-\frac{C_{PI}}{4n^{3/2}(n-1)h_0^{1/2}}\left[\sum_{i=1}^n\sum_{j\neq i}^n\sum_{j\neq i,j}^nW_{ij}S_k+\sum_{i=1}^n\sum_{j\neq i}^n\{W_{ij}S_i+W_{ij}S_j\}\right]\\
     \Lambda_4(x) &= -\frac{1}{2}(nh_0)^{-1/2}\mu_{20}^{-3/2}\Delta_1\Delta_2\\
     &\qquad =-\frac{1}{2}(nh_0)^{-3/2}\mu_{20}^{-1}\left[\sum_{i=1}^n\sum_{j\neq i}^nS_iT_j + \sum_{i=1}^nS_iT_i\right]\\
     \Lambda_5(x) &= n^{-1/2}h_0^{1/2}\mu_{20}^{-3/2}f(x)\Delta_1^2\\
     &\qquad = n^{-3/2}h_0^{-1/2}\mu_{20}^{-1/2}f(x)\left[\sum_{i=1}^n\sum_{j\neq i}^nS_iS_j+\sum_{i=1}^nS_i^2\right]\\
     \Lambda_6(x) &= \frac{3}{8}(nh_0)^{-1}\mu_{20}^{-5/2}\Delta_1\Delta_2^2\\
     &\qquad =\frac{3}{8}(nh_0)^{-5/2}\mu_{20}^{-2}\left[\sum_{i=1}^n\sum_{j\neq i}^n\sum_{k\neq  i,j}^nS_iT_jT_k+\sum_{i=1}^n\sum_{j\neq i}^n\{2S_iT_iT_j+S_iT_j^2\}+\sum_{i=1}^nS_iT_i^2\right]\\
     \Lambda_7(x) &=\frac{1}{2}\mu_{20}^{-1/2}\Delta_1\left(\frac{\hat{h}-h_0}{h_0}\right)\\
     &\qquad = -\frac{C_{PI}}{2n^{3/2}h_0^{1/2}}\left[\sum_{i=1}^n\sum_{j\neq i}^nS_iV_j+\sum_{i=1}^nS_iV_i\right]-\frac{C_{PI}}{4n^{3/2}(n-1)h_0^{1/2}}\left[\sum_{i=1}^n\sum_{j\neq i}^n\sum_{j\neq i,j}^nW_{ij}S_k+\sum_{i=1}^n\sum_{j\neq i}^n\{W_{ij}S_i+W_{ij}S_j\}\right]\\
     \Lambda_8(x) &=\mu_{20}^{-3/2}\delta\Delta_1\left(\frac{\hat{h}-h_0}{h_0}\right)\\
     &\qquad = -\frac{C_{PI}\delta}{n^{3/2}h_0^{1/2}\mu_{20}}\left[\sum_{i=1}^n\sum_{j\neq i}^nS_iV_j+\sum_{i=1}^nS_iV_i\right]-\frac{C_{PI}\delta}{2n^{3/2}(n-1)h_0^{1/2}\mu_{20}}\left[\sum_{i=1}^n\sum_{j\neq i}^n\sum_{j\neq i,j}^nW_{ij}S_k+\sum_{i=1}^n\sum_{j\neq i}^n\{W_{ij}S_i+W_{ij}S_j\}\right]
\end{align*}
\begin{align*}
    \Lambda_4(x)^2 &= \frac{1}{4}(nh_0)^{-3}\mu_{20}^{-2}\sum_{i=1}^n\sum_{j=1}^n\sum_{k=1}^n\sum_{l=1}^nS_iT_jS_kT_l\\
    &= \frac{1}{4}(nh_0)^{-3}\mu_{20}^{-2}\Biggl[\sum_{i=1}^n\sum_{j\neq i}^n\sum_{k\neq i,j}^n\sum_{l\neq i,j,k}^nS_iT_jS_kT_l+\sum_{i=1}^n\sum_{j\neq i}^n\sum_{k\neq i,j}^n\{4S_iT_iS_jT_k+S_i^2T_jT_k+S_iS_jT_k^2\}\\
    &\qquad+\sum_{i=1}^n\sum_{j\neq i}^n\{2S_i^2T_iT_j+2S_iT_i^2S_j+S_i^2T_j^22S_iT_iS_jT_j\}+\sum_{i=1}^nS_i^2T_i^2\Biggl]
\end{align*}

\subsection{Expansion of Characteristic Function}
In the following expansion of characteristic function of $S_{pilot}(x)$, we use
\begin{align*}
    & \mathbb{E}[S_1\Gamma_1]=h_0\mu_{20}^{-1/2}\xi_{11},\\
    & \mathbb{E}[S_1V_1]=\mathbb{E}[S_1\mathcal{L}_1]=h_0\mu_{20}^{-1/2}\rho_{11},\\
    & \mathbb{E}[S_1S_2W_{12}] = b^{-2L}h_0\mu_{20}^{-1}\omega_{111},\\
    & \mathbb{E}[S_1\Gamma_2W_{12}] = b^{-2L}h_0\mu_{20}^{-1/2}\psi_{111}.
\end{align*}
Define the characteristic function of $T_{PI}$ as $\chi_{T_{PI}}(t)$ as follows.
\begin{align*}
    \hat{\chi}_{PI}(t) &= \mathbb{E}\left[\exp\left\{it\left(S(x)+\sum_{k=1}^8\Lambda_k(x)\right)\right\}\right]\\
    &= \mathbb{E}\left\{e^{itS(x)}\left(1+\sum_{k=1}^8it\Lambda_k(x)+\frac{1}{2}(it\Lambda_4(x))^2\right)\right\}+o\{(nh_0)^{-1}\}\\
    &= (I)+(II)+(III)+(IV)+(V)+(VI)+(VII)+(VIII)+(IX)+(X)+o\{(nh_0)^{-1}\}
\end{align*}
Define $\gamma(t)=\frac{it}{\sqrt{nh_0}}S_i$. Note that for $m=0,1,2,3$
\begin{align}
    \gamma(t)^{n-m}&= \exp\left(\frac{-t^2}{2}\right)\Biggl[\Biggl\{1+\frac{\mu_{30}\mu_{20}^{-3/2}}{6n^{1/2}h_0^{1/2}}(it)^3+\frac{\mu_{40}\mu_{20}^{-2}}{24nh_0}(it)^4+\frac{\mu_{30}^2\mu_{20}^{-3}}{72nh_0}(it)^6\Biggl\}\Biggl] + o\{(nh_0)^{-1}\}. \nonumber
\end{align}

\begin{align}
    (I)=\gamma(t)^n = \exp\left(\frac{-t^2}{2}\right)\Biggl[\Biggl\{1+\frac{\mu_{30}\mu_{20}^{-3/2}}{6n^{1/2}h_0^{1/2}}(it)^3+\frac{\mu_{40}\mu_{20}^{-2}}{24nh_0}(it)^4+\frac{\mu_{30}^2\mu_{20}^{-3}}{72nh_0}(it)^6\Biggl\}\Biggl] +o\{(nh_0)^{-1}\} \nonumber
\end{align}

Since we have already driven $(II),(III)$ and $(IV)$ in the proof of Theorem \ref{T:Edgeworth Expansion Including Pilot Bandwidth}, we quote them.

Next, we expand $(IV)$ and $(V)$, which is including the effect of studentisation.
\begin{align*}
    (V)&=\mathbb{E}\left[e^{itS(x)}it\Lambda_4(x)\right]\\
    &=-\frac{1}{2}(nh_0)^{-3/2}\mu_{20}^{-1}\mathbb{E}\left[e^{itS(x)}\left\{\sum_{i=1}^n\sum_{j\neq i}^nS_iT_j + \sum_{i=1}^nS_iT_i\right\}\right](it)\\
    &= -\frac{1}{2}n^{1/2}h_0^{-3/2}\mu_{20}^{-1}\gamma(t)^{n-2}\mathbb{E}\left[\left\{1+\frac{it}{\sqrt{nh_0}}(S_1+S_2)+\frac{(it)^2}{2nh_0}(S_1+S_2)^2+\frac{(it)^3}{6\sqrt{nh_0}^3}(S_1+S_2)^3\right\}S_1T_2\right](it)\\
    &\qquad -\frac{1}{2}n^{-1/2}h_0^{-3/2}\mu_{20}^{-1}\gamma(t)^{n-1}\mathbb{E}\left[\left\{1+\frac{it}{\sqrt{nh_0}}S_1\right\}S_1T_1\right](it)+o\{(nh)^{-1}\}\\
    &= -\frac{1}{2}n^{-1/2}h_0^{-5/2}\gamma(t)^{n-2}\mu_{20}^{-1}\mathbb{E}[S_1^2]\mathbb{E}[S_2T_2](it)^3-\frac{1}{4}n^{-1}h_0^{-3}\gamma(t)^{n-2}\mu_{20}^{-1}\left\{\mathbb{E}[S_1^3]\mathbb{E}[S_2T_2]+\mathbb{E}[S_1^2]\mathbb{E}[S_2^2T_2]\right\}(it)^4\\
    &\qquad-\frac{1}{2}n^{-1/2}h_0^{-3/2}\gamma(t)^{n-1}\mu_{20}^{-1}\mathbb{E}[S_1T_1](it)-\frac{1}{2}n^{-1}h_0^{-2}\gamma(t)^{n-1}\mu_{20}^{-1}\mathbb{E}[S_1^2T_1](it)^2 +o\{(nh_0)^{-1}\}\\
    &= -\frac{1}{2}n^{-1/2}h_0^{-1/2}\gamma(t)^{n-2}\mu_{20}^{-3/2}\mu_{11}(it)^3-\frac{1}{4}n^{-1}h_0^{-1}\gamma(t)^{n-2}\left\{\mu_{20}^{-3}\mu_{30}\mu_{11}+\mu_{20}^{-2}\mu_{21}\right\}(it)^4\\
    &\qquad-\frac{1}{2}n^{-1/2}h_0^{-1/2}\gamma(t)^{n-1}\mu_{20}^{-3/2}\mu_{11}(it)-\frac{1}{2}n^{-1}h_0^{-1}\gamma(t)^{n-1}\mu_{20}^{-2}\mu_{21}(it)^2+o\{(nh_0)^{-1}\} +o\{(nh_0)^{-1}\}\\
    &= \exp\left(\frac{-t^2}{2}\right)\Biggl[-\frac{1}{2}(nh_0)^{-1/2}\mu_{20}^{-3/2}\mu_{11}\{(it)^3+(it)\}-\frac{1}{12}(nh_0)^{-1}\mu_{20}^{-3}\mu_{30}\mu_{11}\{(it)^6+(it)^4\}\\
    &\qquad-\frac{1}{4}(nh_0)^{-1}\mu_{20}^{-3}\mu_{30}\mu_{11}(it)^4-\frac{1}{4}(nh_0)^{-1}\mu_{20}^{-2}\mu_{21}\{(it)^4+2(it)^2\}\Biggl]+o\{(nh_0)^{-1}\}\\
    &=\exp\left(\frac{-t^2}{2}\right)\Biggl[-\frac{1}{2}(nh_0)^{-1/2}\mu_{20}^{-3/2}\mu_{11}\{(it)^3+(it)\}-\frac{1}{12}(nh_0)^{-1}\left(\mu_{20}^{-3}\mu_{30}\mu_{11}\{(it)^6+4(it)^4\}\right)\\
    &\qquad-\frac{1}{4}(nh_0)^{-1}\mu_{20}^{-2}\mu_{21}\{(it)^4+2(it)^2\}\Biggl]+o\{(nh_0)^{-1}\}
\end{align*}

\begin{align*}
    (VI)&=\mathbb{E}\left[e^{itS(x)}it\Lambda_5(x)\right]\\
    &= n^{-3/2}h_0^{-1/2}\mu_{20}^{-1/2}f(x)\mathbb{E}\left[e^{itS(x)}\left\{\sum_{i=1}^n\sum_{j\neq i}^nS_iS_j+\sum_{i=1}^nS_i^2\right\}\right](it)\\
    &= n^{1/2}h_0^{-1/2}\mu_{20}^{-1/2}f(x)\gamma(t)^{n-2}\mathbb{E}\left[\left\{1+\frac{it}{\sqrt{nh_0}}(S_1+S_2)+\frac{(it)^2}{2nh_0}(S_1+S_2)^2\right\}S_1S_2\right](it)\\
    &\qquad+n^{-1/2}h_0^{-1/2}\mu_{20}^{-1/2}f(x)\gamma(t)^{n-1}\mathbb{E}\left[S_1^2\right](it)+o\{(nh_0)^{-1}\}\\
    &=n^{-1/2}h_0^{-3/2}f(x)\gamma(t)^{n-2}\mu_{20}^{-1/2}\mathbb{E}[S_1^2]^2(it)^3+n^{-1/2}h_0^{-1/2}f(x)\gamma(t)^{n-1}\mu_{20}^{-1/2}\mathbb{E}[S_1^2](it)+o\{(nh_0)^{-1}\}\\
    &=n^{-1/2}h_0^{1/2}f(x)\mu_{20}^{-1/2}\{\gamma(t)^{n-2}(it)^3+\gamma(t)^{n-1}(it)\}+o\{(nh_0)^{-1}\}\\
    &=\exp\left(\frac{-t^2}{2}\right)\Biggl[n^{-1/2}h_0^{1/2}f(x)\mu_{20}^{-1/2}\{(it)^3+(it)\}\Biggl]+o\{(nh_0)^{-1}\}
\end{align*}

\begin{align*}
    (VII)&=\mathbb{E}\left[e^{itS(x)}it\Lambda_6(x)\right]\\
    &=\frac{3}{8}(nh_0)^{-5/2}\mu_{20}^{-2}\mathbb{E}\left[e^{itS(x)}\left\{\sum_{i=1}^n\sum_{j\neq i}^n\sum_{k\neq  i,j}^nS_iT_jT_k+\sum_{i=1}^n\sum_{j\neq i}^n\{2S_iT_iT_j+S_iT_j^2\}+\sum_{i=1}^nS_iT_i^2\right\}\right](it)\\
    &=\frac{3}{8}n^{1/2}h_0^{-5/2}\gamma(t)^{n-3}\mathbb{E}\left[\left\{1+\cdots+\frac{(it)^3}{(nh_0)^{3/2}}(S_1+S_2+S_3)^3\right\}S_1T_2T_3\right](it)\\
    &\qquad+\frac{3}{8}n^{-1/2}h_0^{-5/2}\gamma(t)^{n-2}\mu_{20}^{-2}\mathbb{E}\left[\left\{1+\frac{it}{\sqrt{nh_0}}(S_1+S_2)\right\}\{2S_1T_1T_2+S_1T_2^2\}\right](it)+o\{(nh_0)^{-1}\}\\
    &=\frac{3}{8}n^{-1}h_0^{-4}\gamma(t)^{n-3}\mu_{20}^{-2}\mathbb{E}[S_1^2]\mathbb{E}[S_2T_2]^2(it)^4\\
    &\qquad+\frac{3}{8}n^{-1}h_0^{-3}\gamma(t)^{n-2}\mu_{20}^{-2}\{2\mathbb{E}[S_1T_1]^2+\mathbb{E}[S_1^2]\mathbb{E}[T_2^2]\}(it)^2+o\{(nh_0)^{-1}\}\\
    &=\frac{3}{8}(nh_0)^{-1}\gamma(t)^{n-3}\mu_{20}^{-3}\mu_{11}^2(it)^4+\frac{3}{8}(nh_0)^{-1}\gamma(t)^{n-2}\mu_{20}^{-2}\{2\mu_{20}^{-1}\mu_{11}^2+\mu_{02}\}(it)^2\\
    &=\exp\left(\frac{-t^2}{2}\right)\Biggl[\frac{3}{8}(nh_0)^{-1}\left(\mu_{20}^{-3}\mu_{11}^2\{(it)^4+2(it)^2\}+\mu_{20}^{-2}\mu_{02}(it)^2\right)\Biggl]+o\{(nh_0)^{-1}\}
\end{align*}

\begin{align*}
    (X)&=\mathbb{E}\left[e^{itS(x)}\frac{1}{2}(it\Lambda_4(x))^2\right]\\
    &=\frac{1}{8}(nh_0)^{-3}\mu_{20}^{-2}\mathbb{E}\Biggl[e^{itS(x)}\sum_{i=1}^n\sum_{j\neq i}^n\sum_{k\neq i,j}^n\sum_{l\neq i,j,k}^nS_iT_jS_kT_l+\sum_{i=1}^n\sum_{j\neq i}^n\sum_{k\neq i,j}^n\{4S_iT_iS_jT_k+S_i^2T_jT_k+S_iS_jT_k^2\}\\
    &\qquad+\sum_{i=1}^n\sum_{j\neq i}^n\{2S_i^2T_iT_j+2S_iT_i^2S_j+S_i^2T_j^2+2S_iT_iS_jT_j\}+\sum_{i=1}^nS_i^2T_i^2\Biggl](it)^2\\
    &= \frac{1}{8}(nh_0)^{-1}\gamma(t)^{n-4}\mu_{20}^{-3}\mu_{11}^2(it)^6\\
    &\qquad + \frac{5}{8}(nh_0)^{-1}\gamma(t)^{n-3}\mu_{20}^{-3}\mu_{11}^2(it)^4+\frac{1}{8}(nh_0)^{-1}\gamma(t)^{n-3}\mu_{20}^{-2}\mu_{02}(it)^4\\
    &\qquad+\frac{1}{8}(nh_0)^{-1}\gamma(t)^{n-2}\mu_{20}^{-2}\mu_{02}(it)^2+\frac{1}{4}(nh_0)^{-1}\gamma(t)^{n-2}\mu_{20}^{-3}\mu_{11}^2(it)^2+o\{(nh_0)^{-1}\}\\
    &=\exp\left(\frac{-t^2}{2}\right)\Biggl[\frac{1}{8}(nh_0)^{-1}\left(\mu_{20}^{-3}\mu_{11}^2\{(it)^6+5(it)^4+2(it)^2\}+\mu_{20}^{-2}\mu_{02}\{(it)^4+(it)^2\}\right)\Biggl]+o\{(nh_0)^{-1}\}
\end{align*}

Finally, we expand $(VIII)$ and $(IX)$, which are including the simultaneous effect of studentisation and bandwidth selection.
\begin{align*}
    (VIII)&=\mathbb{E}\left[e^{itS(x)}it\Lambda_7(x)\right]\\
    &= -\frac{C_{PI}}{2n^{3/2}h_0^{1/2}}\mathbb{E}\left[e^{itS(x)}\left\{\sum_{i=1}^n\sum_{j\neq i}^nS_iV_j+\sum_{i=1}^nS_iV_i\right\}\right](it)\\
    &\qquad-\frac{C_{PI}}{4n^{3/2}(n-1)h_0^{1/2}}\mathbb{E}\left[e^{itS(x)}\left\{\sum_{i=1}^n\sum_{j\neq i}^n\sum_{j\neq i,j}^nW_{ij}S_k+\sum_{i=1}^n\sum_{j\neq i}^n\{W_{ij}S_i+W_{ij}S_j\}\right\}\right](it)\\
    &= \frac{-C_{PI}(n-1)}{2n^{1/2}h_0^{1/2}}\gamma(t)^{n-2}\mathbb{E}\left[\left\{1+\cdots+\frac{(it)^2}{2nh_0}(S_1+S_2)^2\right\}S_1V_2\right](it) -\frac{C_{PI}}{2n^{1/2}h_0^{1/2}}\gamma(t)^{n-1}\mathbb{E}[S_1V_1](it)\\
    &\qquad -\frac{C_{PI}(n-2)}{4n^{1/2}h_0^{1/2}}\gamma(t)^{n-3}\mathbb{E}\left[\left\{1+\cdots+\frac{(it)^3}{6(nh_0)^{3/2}}(S_1+S_2+S_3)^3\right\}W_{12}S_3\right](it)\\
    &\qquad-\frac{C_{PI}}{2n^{1/2}h_0^{1/2}}\gamma(t)^{n-2}\mathbb{E}\left[\left\{1+\frac{it}{(nh_0)^{1/2}}(S_1+S_2)\right\}W_{12}S_1\right](it)+o\{(nh_0)^{-1}\}\\
    &=\frac{-C_{PI}}{2n^{1/2}h_0^{3/2}}\gamma(t)^{n-2}\mathbb{E}[S_1^2]\mathbb{E}[S_2V_2](it)^3-\frac{C_{PI}}{2n^{1/2}h_0^{1/2}}\gamma(t)^{n-1}\mathbb{E}[S_1V_1](it)\\
    &\qquad-\frac{C_PI}{4nh_0^2}\gamma(t)^{n-3}\mathbb{E}[S_1S_2W_{12}]\mathbb{E}[S_3^2](it)^4-\frac{C_{PI}}{2nh_0}\gamma(t)^{n-2}\mathbb{E}[S_1S_2W_{12}](it)^2+o\{(nh_0)^{-1}\}\\
    &=\frac{-C_{PI}h_0^{1/2}}{2n^{1/2}}\gamma(t)^{n-2}\mu_{20}^{-1/2}\rho_{11}(it)^3-\frac{C_{PI}h_0^{1/2}}{2n^{1/2}}\gamma(t)^{n-1}\mu_{20}^{-1/2}\rho_{11}(it)\\
    &\qquad-\frac{C_{PI}}{4nb^{2L}}\gamma(t)^{n-3}\mu_{20}^{-1}\omega_{111}(it)^4-\frac{C_{PI}}{2nb^{2L}}\gamma(t)^{n-2}\mu_{20}^{-1}\omega_{111}(it)^2+o\{(nh_0)^{-1}\}\\
    &=\exp\left(\frac{-t^2}{2}\right)\Biggl[\frac{-C_{PI}h_0^{1/2}}{2n^{1/2}}\mu_{20}^{-1/2}\rho_{11}\{(it)^3+(it)\}-\frac{C_{PI}}{4nb^{2L}}\mu_{20}^{-1}\omega_{111}\{(it)^4+2(it)^2\}\Biggl]+o\{(nh_0)^{-1}\}
\end{align*}

\begin{align*}
    (IX)&=\mathbb{E}\left[e^{itS(x)}it\Lambda_8(x)\right]\\
    &=-\frac{C_{PI}\delta}{n^{3/2}h_0^{1/2}\mu_{20}}\mathbb{E}\left[e^{itS(x)}\left\{\sum_{i=1}^n\sum_{j\neq i}^nS_iV_j+\sum_{i=1}^nS_iV_i\right\}\right](it)\\
    &\qquad-\frac{C_{PI}\delta}{2n^{3/2}(n-1)h_0^{1/2}\mu_{20}}\mathbb{E}\left[e^{itS(x)}\left\{\sum_{i=1}^n\sum_{j\neq i}^n\sum_{j\neq i,j}^nW_{ij}S_k+\sum_{i=1}^n\sum_{j\neq i}^n\{W_{ij}S_i+W_{ij}S_j\}\right\}\right](it)\\
    &= -\frac{C_{PI}\delta (n-1)}{n^{1/2}h_0^{1/2}\mu_{20}}\gamma(t)^{n-2}\mathbb{E}\left[\left\{1+\cdots+\frac{(it)^2}{2nh_0}(S_1+S_2)^2\right\}S_1V_2\right](it)-\frac{C_{PI}\delta}{n^{1/2}h_0^{1/2}\mu_{20}}\gamma(t)^{n-1}\mathbb{E}[S_1V_1](it)\\
    &\qquad-\frac{C_{PI}\delta(n-2)}{2n^{1/2}h_0^{1/2}\mu_{20}}\gamma(t)^{n-3}\mathbb{E}\left[\left\{1+\cdots+\frac{(it)^3}{6(nh_0)^{3/2}}(S_1+S_2+S_3)^3\right\}W_{12}S_3\right](it)\\
    &\qquad-\frac{C_{PI}\delta}{n^{1/2}h_0^{1/2}\mu_{20}}\gamma(t)^{n-2}\mathbb{E}\left[\left\{1+\frac{it}{(nh_0)^{1/2}}(S_1+S_2)\right\}W_{12}S_1\right](it) +o\{(nh_0)^{-1}\}\\
    &=-\frac{C_{PI}\delta}{n^{1/2}h_0^{3/2}\mu_{20}}\gamma(t)^{n-2}\mathbb{E}[S_1^2]\mathbb{E}[S_2V_2](it)^3-\frac{C_{PI}\delta}{n^{1/2}h_0^{1/2}\mu_{20}}\gamma(t)^{n-1}\mathbb{E}[S_1V_1](it)\\
    &\qquad-\frac{C_{PI}\delta}{2nh_0^2\mu_{20}}\gamma(t)^{n-3}\mathbb{E}[S_1S_2W_{12}]\mathbb{E}[S_3^2](it)^4-\frac{C_{PI}\delta}{nh_0\mu_{20}}\gamma(t)^{n-2}\mathbb{E}[S_1S_2W_{12}](it)^2+o\{(nh_0)^{-1}\}\\
    &= -\frac{C_{PI}\delta h_0^{1/2}}{n^{1/2}}\gamma(t)^{n-2}\mu_{20}^{-1}\mu_{20}^{-1/2}\rho_{11}(it)^3-\frac{C_{PI}\delta h_0^{1/2}}{n^{1/2}}\gamma(t)^{n-1}\mu_{20}^{-1}\mu_{20}^{-1/2}\rho_{11}(it)\\
    &\qquad -\frac{C_{PI}\delta }{2nb^{2L}}\gamma(t)^{n-3}\mu_{20}^{-2}\omega_{111}(it)^4-\frac{C_{PI}\delta }{nb^{2L}}\gamma(t)^{n-2}\mu_{20}^{-2}\omega_{111}(it)^2+o\{(nh_0)^{-1}\}\\
    &=\exp\left(\frac{-t^2}{2}\right)\left[-\frac{C_{PI}\delta h_0^{1/2}}{n^{1/2}}\mu_{20}^{-3/2}\rho_{11}\{(it)^3+(it)\}-\frac{C_{PI}\delta }{2nb^{2L}}\mu_{20}^{-2}\omega_{111}\{(it)^4+2(it)^2\}\right]+o\{(nh_0)^{-1}\}
\end{align*}

\subsubsection{Simplifying Terms in the Expansion which Includes the Effect of Studentisation}
$(I),(V),(VI),(VII)$ and $(X)$ includes the effect of studentization. The Edgeworth expansion has been already derived in \cite{Hall91}. 
\begin{align*}
    &(I)+(V)+(VI)+(VII)+(X) \\
    &= \exp\left(\frac{-t^2}{2}\right)\Biggl[\Biggl\{1+\frac{\mu_{30}\mu_{20}^{-3/2}}{6n^{1/2}h_0^{1/2}}(it)^3+\frac{\mu_{40}\mu_{20}^{-2}}{24nh_0}(it)^4+\frac{\mu_{30}^2\mu_{20}^{-3}}{72nh_0}(it)^6\Biggl\}\\
    &-\frac{1}{2}(nh_0)^{-1/2}\mu_{20}^{-3/2}\mu_{11}\{(it)^3+(it)\}-\frac{1}{12}(nh_0)^{-1}\left(\mu_{20}^{-3}\mu_{30}\mu_{11}\{(it)^6+4(it)^4\}\right)-\frac{1}{4}(nh_0)^{-1}\mu_{20}^{-2}\mu_{21}\{(it)^4+2(it)^2\}\\
    &+n^{-1/2}h_0^{1/2}f(x)\mu_{20}^{-1/2}\{(it)^3+(it)\}+\frac{3}{8}(nh_0)^{-1}\left(\mu_{20}^{-3}\mu_{11}^2\{(it)^4+2(it)^2\}+\mu_{20}^{-2}\mu_{02}(it)^2\right)\\
    &+\frac{1}{8}(nh_0)^{-1}\left(\mu_{20}^{-3}\mu_{11}^2\{(it)^6+5(it)^4+2(it)^2\}+\mu_{20}^{-2}\mu_{02}\{(it)^4+(it)^2\}\right)\Biggl]+o\{(nh_0)^{-1}\}\\
    &= \exp\left(\frac{-t^2}{2}\right)\Biggl[(nh_0)^{-1/2}\Biggl\{\frac{-1}{2}\mu_{20}^{-3/2}\mu_{11}(it)+\frac{1}{6}\mu_{20}^{-3/2}(\mu_{30}-3\mu_{11})(it)^3\Biggl\}+n^{-1/2}h_0^{1/2}f(x)\mu_{20}^{-1/2}\{(it)^3+(it)\}\\
    &\qquad+(nh_0)^{-1}\Biggl\{\left(-\frac{1}{2}\mu_{20}^{-2}\mu_{21}+\mu_{20}^{-3}\mu_{11}^2+\frac{1}{2}\mu_{20}^{-2}\mu_{02}\right)(it)^2\\
    &\qquad\qquad+\left(\frac{1}{24}\mu_{40}\mu_{20}^{-2}-\frac{1}{3}\mu_{20}^{-3}\mu_{30}\mu_{11}-\frac{1}{4}\mu_{20}^{-2}\mu_{21}+\mu_{20}^{-3}\mu_{11}^2+\frac{1}{8}\mu_{20}^{-2}\mu_{02}\right)(it)^4\\
    &\qquad\qquad+\left(\frac{1}{72}\mu_{30}^2\mu_{20}^{-3}-\frac{1}{12}\mu_{20}^{-3}\mu_{30}\mu_{11}+\frac{1}{8}\mu_{20}^{-3}\mu_{11}^2\right)(it)^6\Biggl\}\Biggl] +o\{(nh_0)^{-1}\}
\end{align*}
As stated in \cite[p.214]{Hall92book}, if the non-negative integer $i,j,k,l$ satisfy $i+2j=k+2l$ then $\mu_{ij}-\mu_{kl}=O(h)\implies \mu_{ij}=\mu_{kl}+O(h)$. Using this statements, $\mu_{02}=\mu_{40}+O(h_0),\mu_{21}=\mu_{40}+O(h_0)$ and $\mu_{11}=\mu_{30}+O(h_0)$.
\begin{align*}
    &(I)+(V)+(VI)+(VII)+(X) \\
    &= \exp\left(\frac{-t^2}{2}\right)\Biggl[(nh_0)^{-1/2}\Biggl\{\frac{-1}{2}\mu_{20}^{-3/2}\mu_{11}(it)+\frac{1}{6}\mu_{20}^{-3/2}(\mu_{30}-3\mu_{11})(it)^3\Biggl\}+n^{-1/2}h_0^{1/2}f(x)\mu_{20}^{-1/2}\{(it)^3+(it)\}\\
    &\qquad+(nh_0)^{-1}\Biggl\{\mu_{20}^{-3}\mu_{11}^2(it)^2+\left(\frac{2}{3}\mu_{20}^{-2}\mu_{30}^2-\frac{1}{12}\mu_{40}\mu_{20}^{-2}\right)(it)^4+\frac{1}{18}\mu_{20}^{-3}\mu_{30}^2(it)^6\Biggl\}\Biggl] +o\{(nh_0)^{-1}\}
\end{align*}
This is consistent with the expansion in \cite{Hall91,Hall92book}.

\subsubsection{Simplifying Terms in the Expansion which Includes the Effect of Global Plug-In Method}
\begin{align*}
    &(II)+(III)+(IV)\\
    &= \exp\left(\frac{-t^2}{2}\right)\Biggl[C_{PI}\mu_{20}^{-1}\rho_{11}\left(\sum_{l=0}^{L-1}C_{\Gamma,l}(x)h_0^{L+l+1}\right)(it)^2\\
    &\qquad+\frac{C_{PI}}{2}n^{-1/2}h_0^{\frac{2L+1}{2}}b^{-2L}\mu_{20}^{-3/2}\omega_{111}\left(\sum_{l=0}^{L-1}C_{\Gamma,l}(x)h_0^{l}\right)(it)^3\\
    &\qquad+C_{PI}n^{-1/2}h_0^{1/2}\Biggl(\mu_{20}^{-3/2}\rho_{11}\xi_{11}(it)^3 -\frac{\mu_{20}^{-1/2}}{2}\rho_{11}\{(it)^3+(it)\}\Biggl) \nonumber\\
    &\qquad+C_{PI}n^{-1}b^{-2L}\left(\frac{1}{2}\mu_{20}^{-2}\xi_{11}\omega_{111}(it)^4+\mu_{20}^{-1}\psi_{111}(it)^2-\frac{1}{4}\mu_{20}^{-1}\omega_{111}\{(it)^4+2(it)^2\}\right)\Biggl] +o\{(nh_0)^{-1}\}\\
    & =  \exp\left(\frac{-t^2}{2}\right) \Biggl[C_{PI}\mu_{20}^{-1}\rho_{11}\left(\sum_{l=0}^{L-1}C_{\Gamma,l}(x)h_0^{L+l+1}\right)(it)^2 \\
    &\qquad+\frac{C_{PI}}{2}n^{-1/2}h_0^{\frac{2L+1}{2}}b^{-2L}\mu_{20}^{-3/2}\omega_{111}\left(\sum_{l=0}^{L-1}C_{\Gamma,l}(x)h_0^{l}\right)(it)^3\\
    &\qquad+C_{PI}n^{-1/2}h_0^{1/2}\rho_{11}\left(\mu_{20}^{-3/2}\xi_{11}(it)^3-\frac{\mu_{20}^{-1/2}}{2}\{(it)^3+(it)\}\}\right) \nonumber\\
    &\qquad+C_{PI}n^{-1}b^{-2L}\left(\frac{1}{2}\mu_{20}^{-2}\xi_{11}\omega_{111}(it)^4+\mu_{20}^{-1}\psi_{111}(it)^2-\frac{1}{4}\mu_{20}^{-1}\omega_{111}\{(it)^4+2(it)^2\}\right)\Biggl] +o\{(nh_0)^{-1}\}
\end{align*}

\subsubsection{Simplifying Terms in the Expansion which Includes the Simultaneous Effect of Plug-In Method and Studentisation}
\begin{align*}
    &(VIII)+(IX)\\
    &=\exp\left(\frac{-t^2}{2}\right)\Biggl[\frac{-C_{PI}h_0^{1/2}}{2n^{1/2}}\mu_{20}^{-1/2}\rho_{11}\{(it)^3+(it)\}-\frac{C_{PI}}{4nb^{2L}}\mu_{20}^{-1}\omega_{111}\{(it)^4+2(it)^2\}\\
    &\quad -\frac{C_{PI}\delta h_0^{1/2}}{n^{1/2}}\mu_{20}^{-3/2}\rho_{11}\{(it)^3+(it)\}-\frac{C_{PI}\delta }{2nb^{2L}}\mu_{20}^{-2}\omega_{111}\{(it)^4+2(it)^2\}\Biggl] +o\{(nh_0)^{-1}\}\\
    &=\exp\left(\frac{-t^2}{2}\right)\Biggl[\frac{-C_{PI}h_0^{1/2}}{2n^{1/2}}\left\{1+\delta\mu_{20}^{-1}\right\}\mu_{20}^{-1/2}\rho_{11}\{(it)^3+(it)\}-\frac{C_{PI}}{4nb^{2L}}\mu_{20}^{-1}\omega_{111}\left\{1+\delta\mu_{20}^{-1}\right\}\{(it)^4+2(it)^2\}\Biggl] +o\{(nh_0)^{-1}\}
\end{align*}
Then inverting $\chi_{T_{PI}}(t)$, we have Theorem \ref{T:Edgeworth Expansion Including Pilot Bandwidth and the effect of Studentisation}.

\section{Results of Additional Monte Carlo Studies} \label{section:supp-simulation}

\subsection{Simulation Results with $\hat{I}_L^{convo}$}
\subsubsection{Standard Normal}
We adopt $\# 1: N(0,1)$ in \cite{MW92}. For sample size $n=(50, 100, 400, 1000)$, MSE optimal pilot bandiwdths are $b_0=(0.8596,0.8047,0.7052,0.6462)$. We evaluate the accuracy at the point of $x=0, 0.5, 1, 1.5$ and $x=2$.
Simulation results are as follows.

\begin{table}[H]
\small
\centering
    \caption{$x=0, b=\text{MSE optimal},~~$ scaled second derivative$=0.4122$}
    {\begin{tabular}{c||llllllll}
      \hline & \multicolumn{2}{c|}{$n=50$}   &  \multicolumn{2}{c|}{$n=100$} & \multicolumn{2}{c|}{$n=400$} & \multicolumn{2}{c}{$n=1000$}\\
     \hline & \multicolumn{1}{|c|}{CP} & \multicolumn{1}{|c|}{Ave.Length} & \multicolumn{1}{|c|}{CP} & \multicolumn{1}{|c|}{Ave.Length} & \multicolumn{1}{|c|}{CP} & \multicolumn{1}{|c|}{Ave.Length} & \multicolumn{1}{|c|}{CP} & \multicolumn{1}{|c}{Ave.Length} \\
     \hline\hline $N(0,1)$ & 0.7730 & 0.1339 & 0.6930 & 0.1067 & 0.5950 & 0.0672 & 0.5575 & 0.0490 \\
     Hall (1991) & 0.7670 & 0.1341 & 0.6885 & 0.1068 & 0.5525 & 0.0672 & 0.5525 & 0.0490 \\
     Theorem 3.1 & 0.9140$^{**}$ & 0.1868 & 0.8130$^{**}$ & 0.1343 & 0.6630$^{**}$ & 0.0746 & 0.6020$^{**}$ & 0.0521 \\
     Theorem 3.5 & 0.9090$^{*}$ & 0.2035 & 0.8055$^{*}$ & 0.1421 & 0.6575$^{*}$ & 0.0763 & 0.5990$^{*}$ & 0.0526
    \end{tabular}}
    \label{tab:cSNx0}
\end{table}

\begin{table}[H]
\small
\centering
    \caption{$x=0.5, b=\text{MSE optimal},~~$ scaled second derivative$=0.2728$}
    {\begin{tabular}{c||llllllll}
      \hline & \multicolumn{2}{c|}{$n=50$}   &  \multicolumn{2}{c|}{$n=100$} & \multicolumn{2}{c|}{$n=400$} & \multicolumn{2}{c}{$n=1000$}\\
     \hline & \multicolumn{1}{|c|}{CP} & \multicolumn{1}{|c|}{Ave.Length} & \multicolumn{1}{|c|}{CP} & \multicolumn{1}{|c|}{Ave.Length} & \multicolumn{1}{|c|}{CP} & \multicolumn{1}{|c|}{Ave.Length} & \multicolumn{1}{|c|}{CP} & \multicolumn{1}{|c}{Ave.Length} \\
     \hline\hline $N(0,1)$ & 0.8745 & 0.1342 & 0.8250 & 0.1059 & 0.7815 & 0.0656 & 0.7565 & 0.0475 \\
     Hall (1991) & 0.8685 & 0.1343 & 0.8160 & 0.1059 & 0.7730 & 0.0656 & 0.7475 & 0.0475 \\
     Theorem 3.1 & 0.9120$^{*}$ & 0.1520 & 0.8545$^{*}$ & 0.1145 & 0.7915$^{*}$ & 0.0677 & 0.7600$^{**}$ & 0.0482 \\
     Theorem 3.5 & 0.9225$^{**}$ & 0.1664 & 0.8600$^{**}$ & 0.1212 & 0.7925$^{**}$ & 0.0691 & 0.7590$^{*}$ & 0.0487
    \end{tabular}}
    \label{tab:cSNx05}
\end{table}

\begin{table}[H]
\small
\centering
    \caption{$x=1, b=\text{MSE optimal},~~$ scaled second derivative$=0$}
    {\begin{tabular}{c||llllllll}
      \hline & \multicolumn{2}{c|}{$n=50$}   &  \multicolumn{2}{c|}{$n=100$} & \multicolumn{2}{c|}{$n=400$} & \multicolumn{2}{c}{$n=1000$}\\
     \hline & \multicolumn{1}{|c|}{CP} & \multicolumn{1}{|c|}{Ave.Length} & \multicolumn{1}{|c|}{CP} & \multicolumn{1}{|c|}{Ave.Length} & \multicolumn{1}{|c|}{CP} & \multicolumn{1}{|c|}{Ave.Length} & \multicolumn{1}{|c|}{CP} & \multicolumn{1}{|c}{Ave.Length} \\
     \hline\hline $N(0,1)$ & 0.9730 & 0.1281 & 0.9700$^{*}$ & 0.0990 & 0.9680 & 0.0593 & 0.9640$^{**}$ & 0.0422 \\
     Hall (1991) & 0.9680$^{**}$ & 0.1282 & 0.9695$^{**}$ & 0.0990 & 0.9675$^{**}$ & 0.0593 & 0.9655 & 0.0422 \\
     Theorem 3.1 & 0.9705$^{*}$ & 0.1282 & 0.9720 & 0.0990 & 0.9675$^{**}$ & 0.0593 & 0.9650$^{*}$ & 0.0422 \\
     Theorem 3.5 & 0.9805 & 0.1366 & 0.9765 & 0.1027 & 0.9690 & 0.0601 & 0.9660 & 0.0424
    \end{tabular}}
    \label{tab:cSNx1}
\end{table}

\begin{table}[H]
\small
\centering
    \caption{$x=1.5, b=\text{MSE optimal},~~$ scaled second derivative$=0.1673$}
    {\begin{tabular}{c||llllllll}
      \hline & \multicolumn{2}{c|}{$n=50$}   &  \multicolumn{2}{c|}{$n=100$} & \multicolumn{2}{c|}{$n=400$} & \multicolumn{2}{c}{$n=1000$}\\
     \hline & \multicolumn{1}{|c|}{CP} & \multicolumn{1}{|c|}{Ave.Length} & \multicolumn{1}{|c|}{CP} & \multicolumn{1}{|c|}{Ave.Length} & \multicolumn{1}{|c|}{CP} & \multicolumn{1}{|c|}{Ave.Length} & \multicolumn{1}{|c|}{CP} & \multicolumn{1}{|c}{Ave.Length} \\
     \hline\hline $N(0,1)$ & 0.9500$^{**}$ & 0.1083 & 0.9315 & 0.0819 & 0.8700 & 0.0474 & 0.8315 & 0.0331 \\
     Hall (1991) & 0.9575$^{*}$ & 0.1082 & 0.9430$^{**}$ & 0.0818 & 0.8835 & 0.0474 & 0.8400 & 0.0331 \\
     Theorem 3.1 & 0.9755 & 0.1198 & 0.9595$^{*}$ & 0.0880 & 0.8965$^{**}$ & 0.0490 & 0.8495$^{**}$ & 0.0338 \\
     Theorem 3.5 & 0.9755 & 0.1227 & 0.9600 & 0.0892 & 0.8965$^{**}$ & 0.0493 & 0.8495$^{**}$ & 0.0339
    \end{tabular}}
    \label{tab:cSNx15}
\end{table}

\begin{table}[H]
\small
\centering
    \caption{$x=2, b=\text{MSE optimal},~~$ scaled second derivative$=0.1673$}
    {\begin{tabular}{c||llllllll}
      \hline & \multicolumn{2}{c|}{$n=50$}   &  \multicolumn{2}{c|}{$n=100$} & \multicolumn{2}{c|}{$n=400$} & \multicolumn{2}{c}{$n=1000$}\\
     \hline & \multicolumn{1}{|c|}{CP} & \multicolumn{1}{|c|}{Ave.Length} & \multicolumn{1}{|c|}{CP} & \multicolumn{1}{|c|}{Ave.Length} & \multicolumn{1}{|c|}{CP} & \multicolumn{1}{|c|}{Ave.Length} & \multicolumn{1}{|c|}{CP} & \multicolumn{1}{|c}{Ave.Length} \\
     \hline\hline $N(0,1)$ & 0.8195 & 0.0790 & 0.7645 & 0.0585 & 0.6550 & 0.0328 & 0.6135 & 0.0226 \\
     Hall (1991) & 0.8560 & 0.0786 & 0.7945 & 0.0583 & 0.6805 & 0.0328 & 0.6380 & 0.0226 \\
     Theorem 3.1 & 0.9005$^{**}$ & 0.0872 & 0.8305$^{**}$ & 0.0623 & 0.7005$^{**}$ & 0.0336 & 0.6490$^{**}$ & 0.0229 \\
     Theorem 3.5 & 0.8940$^{*}$ & 0.0875 & 0.8285$^{*}$ & 0.0624 & 0.6975$^{*}$ & 0.0336 & 0.6485$^{*}$ & 0.0229
    \end{tabular}}
    \label{tab:cSNx2}
\end{table}

\subsubsection{Skewed Unimodal}
We adopt $\# 2: \frac{1}{5}N(0,1)+\frac{1}{5}N(\frac{1}{2},(\frac{2}{3})^2)+\frac{3}{5}N(\frac{13}{12},(\frac{5}{9})^2)$ in \cite{MW92}. For sample size $n=(50, 100, 400, 1000)$, MSE optimal pilot bandiwdths are $b_0=(0.5318,
0.4978,0.4363,0.3998)$. We evaluate the accuracy at the point of $x=-2, -1.5, -1, -0.5, 0, 0.5, 1, 1.5$ and $x=2$.
Simulation results are as follows.
\begin{table}[H]
\small
\centering
    \caption{$x=-2, b=\text{MSE optimal},~~$ scaled second derivative$=0.0173$}
    {\begin{tabular}{c||llllllll}
      \hline & \multicolumn{2}{c|}{$n=50$}   &  \multicolumn{2}{c|}{$n=100$} & \multicolumn{2}{c|}{$n=400$} & \multicolumn{2}{c}{$n=1000$}\\
     \hline & \multicolumn{1}{|c|}{CP} & \multicolumn{1}{|c|}{Ave.Length} & \multicolumn{1}{|c|}{CP} & \multicolumn{1}{|c|}{Ave.Length} & \multicolumn{1}{|c|}{CP} & \multicolumn{1}{|c|}{Ave.Length} & \multicolumn{1}{|c|}{CP} & \multicolumn{1}{|c}{Ave.Length} \\
     \hline\hline $N(0,1)$ & 0.9140 & 0.0446 & 0.9220 & 0.0328 & 0.9300 & 0.0182 & 0.9300 & 0.0125 \\
     Hall (1991) & 0.9510$^{**}$ & 0.0417 & 0.9455$^{**}$ & 0.0315 & 0.9470$^{**}$ & 0.0180 & 0.9440$^{**}$ & 0.0124 \\
     Theorem 3.1 & 0.9520$^{*}$ & 0.0418 & 0.9455$^{**}$ & 0.0315 & 0.9470$^{**}$ & 0.0180 & 0.9440$^{**}$ & 0.0124 \\
     Theorem 3.5 & 0.9525 & 0.0418 & 0.9455$^{**}$ & 0.0315 & 0.9470$^{**}$ & 0.0180 & 0.9440$^{**}$ & 0.0124
    \end{tabular}}
    \label{tab:cSUmx2}
\end{table}

\begin{table}[H]
\small
\centering
    \caption{$x=-1.5, b=\text{MSE optimal},~~$ scaled second derivative$=0.0278$}
    {\begin{tabular}{c||llllllll}
      \hline & \multicolumn{2}{c|}{$n=50$}   &  \multicolumn{2}{c|}{$n=100$} & \multicolumn{2}{c|}{$n=400$} & \multicolumn{2}{c}{$n=1000$}\\
     \hline & \multicolumn{1}{|c|}{CP} & \multicolumn{1}{|c|}{Ave.Length} & \multicolumn{1}{|c|}{CP} & \multicolumn{1}{|c|}{Ave.Length} & \multicolumn{1}{|c|}{CP} & \multicolumn{1}{|c|}{Ave.Length} & \multicolumn{1}{|c|}{CP} & \multicolumn{1}{|c}{Ave.Length}\\
     \hline\hline $N(0,1)$ & 0.9200 & 0.0688 & 0.9180 & 0.0508 & 0.9390 & 0.0284 & 0.9425 & 0.0196 \\
     Hall (1991) & 0.9280 & 0.0670 & 0.9330 & 0.0500 & 0.9460 & 0.0283 & 0.9540$^{**}$ & 0.0195 \\
     Theorem 3.1 & 0.9285$^{*}$ & 0.0672 & 0.9345$^{*}$ & 0.0283 & 0.9465$^{**}$ & 0.0283 & 0.9545$^{*}$ & 0.0196 \\
     Theorem 3.5 & 0.9310$^{**}$ & 0.0674 & 0.9355$^{**}$ & 0.0502 & 0.9465$^{**}$ & 0.0283 & 0.9545$^{*}$ & 0.0196
    \end{tabular}}
    \label{tab:cSUm05}
\end{table}

\begin{table}[H]
\small
\centering
    \caption{$x=-1, b=\text{MSE optimal},~~$ scaled second derivative$=0.0503$}
    {\begin{tabular}{c||llllllll}
      \hline & \multicolumn{2}{c|}{$n=50$}   &  \multicolumn{2}{c|}{$n=100$} & \multicolumn{2}{c|}{$n=400$} & \multicolumn{2}{c}{$n=1000$}\\
     \hline & \multicolumn{1}{|c|}{CP} & \multicolumn{1}{|c|}{Ave.Length} & \multicolumn{1}{|c|}{CP} & \multicolumn{1}{|c|}{Ave.Length} & \multicolumn{1}{|c|}{CP} & \multicolumn{1}{|c|}{Ave.Length} & \multicolumn{1}{|c|}{CP} & \multicolumn{1}{|c}{Ave.Length}\\
     \hline\hline $N(0,1)$ & 0.9050 & 0.0727 & 0.8965 & 0.0727 & 0.9170 & 0.0410 & 0.9385 & 0.0283 \\
     Hall (1991) & 0.9300 & 0.0970 & 0.9185 & 0.0723 & 0.9285$^{**}$ & 0.0409 & 0.9470$^{*}$ & 0.0283 \\
     Theorem 3.1 & 0.9345$^{*}$ & 0.0975 & 0.9195$^{*}$ & 0.0725 & 0.9285$^{**}$ & 0.0409 & 0.9470$^{*}$ & 0.0283 \\
     Theorem 3.5 & 0.9360$^{**}$ & 0.0982 & 0.9205$^{**}$ & 0.0728 & 0.9285$^{**}$ & 0.0410 & 0.9475$^{**}$ & 0.0284
    \end{tabular}}
    \label{tab:cSUxm1}
\end{table}

\begin{table}[H]
\small
\centering
    \caption{$x=-0.5, b=\text{MSE optimal},~~$ scaled second derivative$=0.1112$}
    {\begin{tabular}{c||llllllll}
      \hline & \multicolumn{2}{c|}{$n=50$}   &  \multicolumn{2}{c|}{$n=100$} & \multicolumn{2}{c|}{$n=400$} & \multicolumn{2}{c}{$n=1000$}\\
     \hline & \multicolumn{1}{|c|}{CP} & \multicolumn{1}{|c|}{Ave.Length} & \multicolumn{1}{|c|}{CP} & \multicolumn{1}{|c|}{Ave.Length} & \multicolumn{1}{|c|}{CP} & \multicolumn{1}{|c|}{Ave.Length} & \multicolumn{1}{|c|}{CP} & \multicolumn{1}{|c}{Ave.Length}\\
     \hline\hline $N(0,1)$ & 0.9090 & 0.1330 & 0.8850 & 0.0993 & 0.8900 & 0.0566 & 0.8995 & 0.0393 \\
     Hall (1991) & 0.9205 & 0.1325 & 0.9020 & 0.0565 & 0.9030 & 0.0565 & 0.9065 & 0.0393 \\
     Theorem 3.1 & 0.9280$^{*}$ & 0.1353 & 0.9075$^{*}$ & 0.1002 & 0.9050$^{**}$ & 0.0567 & 0.9070$^{**}$ & 0.0393 \\
     Theorem 3.5 & 0.9335$^{**}$ & 0.1375 & 0.9095$^{**}$ & 0.1013 & 0.9050$^{**}$ & 0.0570 & 0.9070$^{**}$ & 0.0394
    \end{tabular}}
    \label{tab:cSUxm05}
\end{table}

\begin{table}[H]
\small
\centering
    \caption{$x=0, b=\text{MSE optimal},~~$ scaled second derivative$=0.2029$}
    {\begin{tabular}{c||llllllll}
      \hline & \multicolumn{2}{c|}{$n=50$}   &  \multicolumn{2}{c|}{$n=100$} & \multicolumn{2}{c|}{$n=400$} & \multicolumn{2}{c}{$n=1000$}\\
     \hline & \multicolumn{1}{|c|}{CP} & \multicolumn{1}{|c|}{Ave.Length} & \multicolumn{1}{|c|}{CP} & \multicolumn{1}{|c|}{Ave.Length} & \multicolumn{1}{|c|}{CP} & \multicolumn{1}{|c|}{Ave.Length} & \multicolumn{1}{|c|}{CP} & \multicolumn{1}{|c}{Ave.Length}\\
     \hline\hline $N(0,1)$ & 0.9620$^{**}$ & 0.1709 & 0.9580$^{**}$ & 0.1299 & 0.9315 & 0.0759 & 0.9205 & 0.0533 \\
     Hall (1991) & 0.9660$^{*}$ & 0.1709 & 0.9635$^{*}$ & 0.1298 & 0.9365 & 0.0758 & 0.9205 & 0.0533 \\
     Theorem 3.1 & 0.9820 & 0.1868 & 0.9735 & 0.1378 & 0.9435$^{*}$ & 0.0777 & 0.9330$^{**}$ & 0.0540 \\
     Theorem 3.5 & 0.9865 & 0.1955 & 0.9760 & 0.1416 & 0.9440$^{**}$ & 0.0785 & 0.9330$^{**}$ & 0.0543    
     \end{tabular}}
    \label{tab:cSUx0}
\end{table}

\begin{table}[H]
\small
\centering
    \caption{$x=0.5, b=\text{MSE optimal},~~$ scaled second derivative$=0.1170$}
    {\begin{tabular}{c||llllllll}
      \hline & \multicolumn{2}{c|}{$n=50$}   &  \multicolumn{2}{c|}{$n=100$} & \multicolumn{2}{c|}{$n=400$} & \multicolumn{2}{c}{$n=1000$}\\
     \hline & \multicolumn{1}{|c|}{CP} & \multicolumn{1}{|c|}{Ave.Length} & \multicolumn{1}{|c|}{CP} & \multicolumn{1}{|c|}{Ave.Length} & \multicolumn{1}{|c|}{CP} & \multicolumn{1}{|c|}{Ave.Length} & \multicolumn{1}{|c|}{CP} & \multicolumn{1}{|c}{Ave.Length}\\
     \hline\hline $N(0,1)$ & 0.9405$^{**}$ & 0.1963 & 0.9290$^{*}$ & 0.1532 & 0.9445$^{**}$ & 0.0933 & 0.9555 & 0.0669 \\
     Hall (1991) & 0.9370 & 0.1964 & 0.9230 & 0.1532 & 0.9395 & 0.0933 & 0.9500$^{**}$ & 0.0669 \\
     Theorem 3.1 & 0.9350 & 0.1936 & 0.9220 & 0.1511 & 0.9380 & 0.0924 & 0.9485 & 0.0665 \\
     Theorem 3.5 & 0.9605$^{*}$ & 0.2176 & 0.9380$^{**}$ & 0.1621 & 0.9435$^{*}$ & 0.0947 & 0.9505$^{*}$ & 0.0673    
     \end{tabular}}
    \label{tab:cSUx05}
\end{table}

\begin{table}[H]
\small
\centering
    \caption{$x=1, b=\text{MSE optimal},~~$ scaled second derivative$=0.7019$}
    {\begin{tabular}{c||llllllll}
      \hline & \multicolumn{2}{c|}{$n=50$}   &  \multicolumn{2}{c|}{$n=100$} & \multicolumn{2}{c|}{$n=400$} & \multicolumn{2}{c}{$n=1000$}\\
     \hline & \multicolumn{1}{|c|}{CP} & \multicolumn{1}{|c|}{Ave.Length} & \multicolumn{1}{|c|}{CP} & \multicolumn{1}{|c|}{Ave.Length} & \multicolumn{1}{|c|}{CP} & \multicolumn{1}{|c|}{Ave.Length} & \multicolumn{1}{|c|}{CP} & \multicolumn{1}{|c}{Ave.Length}\\
     \hline\hline $N(0,1)$ & 0.7605 & 0.2022 & 0.6760 & 0.1595 & 0.5520 & 0.0992 & 0.5275 & 0.0720 \\
     Hall (1991) & 0.7550 & 0.2024 & 0.6680 & 0.1596 & 0.5455 & 0.0992 & 0.5220 & 0.0720 \\
     Theorem 3.1 & 0.9390$^{**}$ & 0.3080 & 0.8370$^{**}$ & 0.2164 & 0.6440$^{**}$ & 0.1152 & 0.5860$^{**}$ & 0.0788 \\
     Theorem 3.5 & 0.9355$^{*}$ & 0.3378 & 0.8270$^{*}$ & 0.2308 & 0.6395$^{*}$ & 0.1184 & 0.5810$^{*}$ & 0.0799    
     \end{tabular}}
    \label{tab:cSUx1}
\end{table}

\begin{table}[H]
\small
\centering
    \caption{$x=1.5, b=\text{MSE optimal},~~$ scaled second derivative$=0.1559$}
    {\begin{tabular}{c||llllllll}
      \hline & \multicolumn{2}{c|}{$n=50$}   &  \multicolumn{2}{c|}{$n=100$} & \multicolumn{2}{c|}{$n=400$} & \multicolumn{2}{c}{$n=1000$}\\
     \hline & \multicolumn{1}{|c|}{CP} & \multicolumn{1}{|c|}{Ave.Length} & \multicolumn{1}{|c|}{CP} & \multicolumn{1}{|c|}{Ave.Length} & \multicolumn{1}{|c|}{CP} & \multicolumn{1}{|c|}{Ave.Length} & \multicolumn{1}{|c|}{CP} & \multicolumn{1}{|c}{Ave.Length}\\
     \hline\hline $N(0,1)$ & 0.8965$^{*}$ & 0.1944 & 0.8855$^{*}$ & 0.1503 & 0.8810$^{**}$ & 0.0903 & 0.8930$^{**}$ & 0.0644 \\
     Hall (1991) & 0.8875 & 0.1945 & 0.8805 & 0.1503 & 0.8735 & 0.0903 & 0.8880 & 0.0644 \\
     Theorem 3.1 & 0.8965$^{*}$ & 0.1990 & 0.8830 & 0.1518 & 0.8735 & 0.0903 & 0.8875 & 0.0643 \\
     Theorem 3.5 & 0.9245$^{**}$ & 0.2170 & 0.8960$^{**}$ & 0.1601 & 0.8795$^{*}$ & 0.0920 & 0.8895$^{*}$ & 0.0649    
     \end{tabular}}
    \label{tab:cSUx15}
\end{table}

\begin{table}[H]
\small
\centering
    \caption{$x=2, b=\text{MSE optimal},~~$ scaled second derivative$=0.3591$}
    {\begin{tabular}{c||llllllll}
      \hline & \multicolumn{2}{c|}{$n=50$}   &  \multicolumn{2}{c|}{$n=100$} & \multicolumn{2}{c|}{$n=400$} & \multicolumn{2}{c}{$n=1000$}\\
     \hline & \multicolumn{1}{|c|}{CP} & \multicolumn{1}{|c|}{Ave.Length} & \multicolumn{1}{|c|}{CP} & \multicolumn{1}{|c|}{Ave.Length} & \multicolumn{1}{|c|}{CP} & \multicolumn{1}{|c|}{Ave.Length} & \multicolumn{1}{|c|}{CP} & \multicolumn{1}{|c}{Ave.Length}\\
     \hline\hline $N(0,1)$ & 0.8600 & 0.1419 & 0.7995 & 0.1060 & 0.6765 & 0.0603 & 0.6145 & 0.0418 \\
     Hall (1991) & 0.8730 & 0.1415 & 0.8195 & 0.1058 & 0.6965 & 0.0602 & 0.6275 & 0.0418 \\
     Theorem 3.1 & 0.9190$^{**}$ & 0.1681 & 0.8670$^{**}$ & 0.1202 & 0.7365$^{**}$ & 0.0641 & 0.6505$^{**}$ & 0.0433 \\
     Theorem 3.5 & 0.9130$^{*}$ & 0.1710 & 0.8630$^{*}$ & 0.1214 & 0.7335$^{*}$ & 0.0643 & 0.6505$^{**}$ & 0.0434    
     \end{tabular}}
    \label{tab:cSUx2}
\end{table}

\subsection{Distributions of $\hat{I}_L$}
In this section, we provide the distributions of $\hat{I}_L$ estimated with $\hat{h}$ and $h_{convo}$ respectively for sample size $n=(50, 100, 400, 1000, 5000)$. Whereas For the simulation to see coverage probabilities, simulation with $n=5000$ is not carried out, we show the result of estimation of $I_L$ at the expense of reducing the iteration form $2000$ to $300$, in order to confirm the consistency of $\hat{I}_L$ to $I_L$. On the following figures, blue histograms represent the distribution of $\hat{I}_L$ and red lines $I_L$.
\subsubsection{Standard Normal}

\subsubsection{Skewed Unimodal}

\subsection{Distributions of KDEs at relatively high curvature points}
On the following figures, the blue histograms represent the distributions of KDEs with plug-in bandwidth ($\hat{f}_{\hat{h}}(x)$ or $\hat{f}_{h_{convo}}(x)$), the orange histograms KDE with optimal bandwidth $\hat{f}_{h_0}(x)$, the red lines $\mathbb{E}[\hat{f}_{h_0}(x)]$ and the black lines $f(x)$.
\subsubsection{Standard Normal at $x=0$}

\subsubsection{Standard Normal at $x=1.5$}

\subsubsection{Skewed Unimodal at $x=1$}

    \newpage
    \bibliographystyle{apalike}
    \bibliography{References}

\end{document}